\documentclass[12pt]{amsart}

\usepackage{amsmath,amssymb}
\usepackage{mathpazo}
\usepackage{amsthm}
\usepackage{comment}
\usepackage{amscd}
\usepackage[all]{xy}
\usepackage{cite}

\usepackage{graphicx}

\usepackage{eq}

\theoremstyle{plain}
\newtheorem{thm}{Theorem}[section]
\newtheorem{lem}[thm]{Lemma}
\newtheorem{prop}[thm]{Proposition}
\newtheorem{cor}[thm]{Corollary}

\theoremstyle{definition}
\newtheorem{mydef}[thm]{Definition}

\newcommand{\cG}{{\mathcal G}}
\newcommand{\cH}{{\mathcal H}}
\newcommand{\cS}{{\mathcal S}}
\newcommand{\cV}{{\mathcal V}}
\newcommand{\cW}{{\mathcal W}}

\newcommand{\ch}{\operatorname{ch}}

\newcommand{\R}{\mathbb{R}}
\newcommand{\C}{\mathbb{C}}

\newcommand{\Proj}{\mathbb{P}}
\newcommand{\Oct}{\mathbb{O}}
\newcommand{\Aut}{\mathrm{Aut}}
\newcommand{\Ker}{\mathrm{Ker}}
\newcommand{\Aff}{\mathrm{Aff}}
\newcommand{\tr}{\mathrm{tr}}
\newcommand{\J}{\mathcal{J}}
\newcommand{\GL}{\mathrm{GL}}
\newcommand{\SL}{\mathrm{SL}}
\newcommand{\SO}{\mathrm{SO}}
\newcommand{\M}{\mathrm{M}}

\newcommand{\Zero}{\mathrm{Zero}}
\newcommand{\Lg}{\mathfrak{g}}
\newcommand{\Lge}{\mathfrak{g_1}}

\newcommand{\Lh}{\mathfrak{h}}
\newcommand{\Le}{\mathfrak{h_1}}
\newcommand{\Lso}{\mathfrak{so}}
\newcommand{\Lgl}{\mathfrak{gl}}
\newcommand{\Spin}{\mathrm{Spin}}

\newcommand{\al}{\alpha}
\newcommand{\be}{\beta}
\newcommand{\ga}{\gamma}
\newcommand{\minimat}{( \begin{smallmatrix} a&b\\ c&d
\end{smallmatrix} \bigr)}
\newcommand{\uppmat}{( \begin{smallmatrix} a&b\\ 0&d
\end{smallmatrix} \bigr)}

\newcommand{\dispfr}[2]{\displaystyle \frac{#1}{#2}}
\newcommand{\oline}[1]{\overline{#1}}

\newcommand{\Sss}{\mathfrak{S}_3}

\newcommand{\ks}{k^{\mathrm{sep}}} 
\newcommand{\HGw}{\mathrm{H}^1(k, G_w)} 
\newcommand{\HGE}{\mathrm{H}^1(k, G_1)}

\newcommand{\KerHGw}{\Ker \left(\mathrm{H}^1(k, G_w) \to \mathrm{H}^1(k, G)\right)}
\newcommand{\Mh}{\mathcal{M}(h)} 
 
\newcommand{\nh}{\mathfrak{n}(h)} 
\newcommand{\at}{\mathfrak{t}}

\newcommand{\Mj}{\mathcal{M}} 
\newcommand{\an}{\mathfrak{n}}

\newcommand{\Gal}{\mathrm{Gal}(k^{\mathrm{sep}}/k)} 

\newcommand{\Autt}{\mathrm{Aut(\J, \at)}}
\newcommand{\Lautt}{\mathfrak{a}(\J,\at)}
\newcommand{\Laut}{\mathfrak{f}}

\newcommand{\bbW}{\mathbb{W}}

\newcommand{\Der}{\mathrm{Der}_k(\J,\J)}
\newcommand{\Lie}{\mathrm{Lie}}
\newcommand{\hj}{h(s_1,s_2,s_3, x_1, x_2, x_3)}
\newcommand{\hji}{h(s_{i1},s_{i2},s_{i3}, x_{i1}, x_{i2}, x_{i3})}
\newcommand{\Trace}{\mathrm{Tr}}
\newcommand{\Vss}{V^{\mathrm{ss}}}

\begin{document}

\address[R. Kato]
{Department of Mathematics, Graduate School of Science, 
Kyoto University, Kyoto 606-8502, Japan}
\email{rkato@math.kyoto-u.ac.jp}

\address[A. Yukie]
{Department of Mathematics, Graduate School of Science, 
Kyoto University, Kyoto 606-8502, Japan}
\email{yukie@math.kyoto-u.ac.jp}
\thanks{The first author was partially supported by 
Grant-in-Aid (B) (24340001)\\}


\title{Rational orbits of the space of pairs of exceptional Jordan algebras}
\author{Ryo Kato}
\author{Akihiko Yukie}
\maketitle


\begin{abstract}
Let $k$ be a field of characteristic not equal to $2,3$, 
$\Oct$ an octonion over $k$ and $\J$ the exceptional Jordan 
algebra defined by $\Oct$. 
We consider the prehomogeneous vector space $(G,V)$
where $G=GE_6\times \GL(2)$ and $V=\J\oplus \J$. 
We prove that  generic rational orbits
of this prehomogeneous vector space are in bijective
correspondence with $k$-isomorphism classes of pairs $(\Mj,\an)$
where $\Mj$'s are isotopes of $\J$ and $\an$'s are 
cubic \'etale subalgebras of $\Mj$. Also we prove that 
if $\Oct$ is split, then generic rational orbits
 are in bijective correspondence with isomorphism classes
 of separable extensions of $k$ of degrees up to $3$. 
\end{abstract}

\section{Introduction}
\label{sec:introduction}

Part of the results in this paper is taken from
the first author's master thesis. 

Let $k$ be a field of characteristic not equal to $2, 3$, 
$k^{\times}=k\setminus \{0\}$,
$\ks$ the separable closure of $k$ and $\oline{k}$ 
the algebraic closure of $k$. We use the notation
$\text{ch}(k)$ for the characteristic of $k$.
If $X$ is a finite set, then let $| X | $ denote its cardinality.
We denote the algebra of $n\times n$ matrices by $\M(n)$ and
the group of $n\times n$ invertible matrices by $\GL(n)$.
Let $\SL(n)=\{g\in\GL(n)\,|\, \det(g)=1\}$.
If $V$ is a finite dimensional vector space over $k$, 
then $\GL(V)$ is the group of invertible
$k$-linear maps from $V$ to $V$ itself.
For $g\in \GL(V)$, $\det (g)$ is well-defined.
Let $\SL(V)=\{g\in\GL(V)\,|\, \det(g)=1\}$.
We denote the Lie algebras of $\GL(n),\GL(V)$
by $\Lgl(n),\Lgl(V)$ respectively.

Let $\widetilde\Oct$ be the split octonion over $k$.  It is the normed algebra
over $k$ obtained by the Cayley--Dickson process
(see \cite[pp.101--110]{harvey}). If $A=\M(2)$ and the norm is the
determinant, then
$\widetilde\Oct$ is $A(+)$ with the notation 
of \cite{harvey}.  An {\it octonion} is, by definition, a normed algebra
which is a $k$-form of $\widetilde\Oct$.
Let $\Oct$ be an octonion. We use the notation $\|x\|$ 
for the norm of $x\in \Oct$.  If $a\in k$, $\|ax\|=a^2\|x\|$. 
For $x, y \in \Oct$, let
\begin{equation*}
Q(x,y)= \tfrac 12(\|x+y\|-\|x\|-\|y\|). 
\end{equation*}
This is a non-degenerate symmetric bilinear form such that $Q(x,x)=\|x\|$. 
Note that $2Q(x,y)$ is denoted by $\langle x,y\rangle$ in \cite[p.1]{Springer}.  
Let $\bbW\subset \Oct$ be the orthogonal complement of $k\cdot 1$
with respect to $Q$.  If $x=x_1+x_2$  where $x_1\in k\cdot 1$,
$x_2\in \bbW$, then we define
\begin{equation}
\label{eq:iota-defn}
\iota(x) = \oline{x} = x_1-x_2
\end{equation}
and call it the {\it conjugate} of $x$. Note that $\|x\|=x\oline{x}$.
The map $\Oct\ni x\mapsto \iota(x)\in\Oct$ is an element of 
$\GL(\Oct)$ where $\Oct$ is regarded as a $k$-vector space. 
For $x \in \Oct$, we define the trace tr$(x)$ by tr$(x) = x + \oline{x}$.
Note that
\begin{equation*}
\mathrm{tr}(x\oline{y}) =
\mathrm{tr}(\oline{y}x) =
2Q(x, y),\; \text{tr}(xy) = \text{tr}(yx). 
\end{equation*}
Also the following properties:
\begin{equation}
\label{eq:trace-relation}
\text{tr}((xy)z) = \text{tr}(x(yz))
\end{equation}
are satisfied for all $x,y,z\in \mathbb O$
(see \cite[p.8, Lemma 1.3.2]{Springer}). 
Therefore, we may write 
${\rm{tr}}(xyz)$ instead of 
${\rm{tr}}((xy)z)$ or ${\rm{tr}}(x(yz))$.

Let  $\J$ be the exceptional Jordan algebra over $k$.
Any element $X \in \J$ is of the form:
\begin{equation}
\label{eq:X-defn}
X = 
\begin{pmatrix}
s_1&x_3&\oline{x_2}\\
\oline{x_3}&s_2&x_1\\
x_2&\oline{x_1}&s_3
\end{pmatrix},
\ \ s_i \in  k, \ x_i \in \Oct \ \ (i =1, 2, 3).
\end{equation}
We sometimes denote this element by $\hj$. 
For elements of $\J$, 
the notion of the determinant is well-defined and is given by 
\begin{equation}
\label{eq:Jdet-defn}
\mathrm{det}(X) = s_1s_2s_3 +\mathrm{tr}(x_1 x_2 x_3)
-s_1 \| x_1 \| -s_2 \| x_2 \| -s_3 \|x_3\|.
\end{equation}

The multiplication  of $\J$ is defined as follows:
\begin{equation*}
X \circ Y = \frac{1}{2} (XY + YX),
\end{equation*}
where the multiplication used on the right-hand side is 
the multiplication of matrices.

Let $\mathfrak{S}_n$ be the symmetric group of $\{1, \cdots, n\}$. 
We define the multiplication of $\sigma,\tau\in \mathfrak{S}_n$
by 
\begin{equation}
(\sigma \tau)(i) = \tau(\sigma(i)) 
\end{equation}
for $i \in \{1,\cdots, n\}$. 
Let $(i \ j)$ denote the transposition of $i\not=j \in \{1,\cdots, n\}$.

The algebraic groups $\SO(Q)$, $E_6$ and $GE_6$ are given by
\begin{align*}
\SO(Q) & = 
 \{ \al \in \SL(\Oct) \,|\, {}^{\forall} x, y \in \Oct, \
 Q(\al (x), \al (y) ) = Q(x, y )
 \}, \\
E_6 & = \{ L \in  \mathrm{GL}(\J) \,|\,
{}^{\forall} X \in \J, \  \det(LX) = \det(X) \}, \\
GE_6 & = \{ L \in  \mathrm{GL}(\J)  \,|\,  
{}^{\forall} X\in \J,   \det(LX) = c(L)\det(X) \
\text{for some} \  c(L) \in \GL (1)  \}
\end{align*}
respectively.
Then $c: GE_6 \to \GL (1)$ is a character and there exists an exact sequence
\begin{equation}\label{E6exact}
0 \to E_6 \hookrightarrow GE_6 \overset{c}{\to} \GL (1) \to 0. 
\end{equation}

It is known that $E_6$ is a smooth connected quasi-simple
simply-connected algebraic group of type ${\mathrm E}_6$ (see \cite[p.181,
Theorem 7.3.2]{Springer}).
The terminology ``quasi-simple'' means that
its inner automorphism group is simple (see
\cite[p.136]{Springer-LAG}). 
Since the dimension of $E_6$ as a variety and the dimension of the 
Lie algebra of $E_6$ coincide (see the proof of \cite[p.181,
Theorem 7.3.2]{Springer}), 
the smoothness of the group follows.


Let $H_1 = E_6$, $G_1 = GE_6$ and  
$H = H_1 \times \GL(2)$  respectively.
Let
\begin{equation}
\label{eq:GV-defn} 
G = G_1 \times \GL(2),\;
V=\J\otimes \Aff^2  
\end{equation}
where $\Aff^2$ is the
2-dimensional affine space regarded as a vector space.  
Then $V$ is a representation of $G$.

We define a character $c'$ of $G$ by the composition of homomorphisms 
\begin{equation}
\label{eq:c'-defn}
G = G_1 \times \GL(2) \overset{\mathrm{pr}_1}\longrightarrow G_1
 \overset{c}\to \GL(1)
\end{equation}
where ``$\mathrm{pr}_1$'' is the natural projection.
The representation $\J$ of $G_1$ is irreducible
(see the proof of \cite[p.181, Theorem 7.3.2]{Springer}) and so $V$ is an irreducible 
representation of $G$.  
The pair $(G,V)$ is what we call a prehomogeneous vector space.

We review the definition of
prehomogeneous vector spaces.
In the following definition, $k$ is an arbitrary field
and $G,V$ are not
necessarily the above $G,V$. 

\begin{mydef}
\label{PVS-def}
Let $G$ be a connected reductive group, $V$ a representation and 
$\chi$ a non-trivial primitive character of $G$, all defined over $k$.
Then,  $(G, V, \chi)$ is called a
{\it prehomogeneous vector space} if it satisfies
the following properties.
\begin{enumerate}
\item There exists a Zariski open orbit. 
\item There exists a non-constant polynomial $\Delta(x) \in k[V]$ such
that $\Delta(gx) = \chi(g)^a \Delta(x)$ for a positive integer $a$.
\end{enumerate}
The polynomial $\Delta$ is called a {\it relative invariant polynomial}.
\end{mydef}

If $(G,V,\chi)$ is an irreducible (as a representation)
prehomogeneous vector space,  
then the choices of $\chi,\Delta$ are essentially unique and we may
write $(G,V)$ instead of $(G,V,\chi)$. We define
$V_k^{\mathrm{ss}} = \{ x \in V_k \,|\, \Delta(x) \neq 0 \}$ 
and call it the set of {\it semi-stable points}.

Now we assume that $\text{ch}(k)\not=2,3$ and 
that $G,V$ are as in (\ref{eq:GV-defn}) again. 
We shall show the existence of an open orbit in 
Proposition \ref{prop:identity-comp}(2).

We identify Aff$^2$ as the space of linear forms in 
two variables $v=(v_1, v_2)$ and 
regard elements of $V$ as the set of $x =x_1 v_1 + x_2 v_2$ where
$x_1, x_2 \in \J$.
The action of $g = (g_1, g_2) \in G$ where $g_2 = \minimat$ on $V$ is 
given by 
\begin{equation}
g(x_1 v_1 + x_2 v_2) = g_1(x_1)(av_1 + cv_2) +  g_1(x_2)(bv_1 + dv_2).
\end{equation}

For $x \in V$, let 
\begin{equation}
\label{eq:Fx-defn} 
F_x(v) = F_x(v_1,v_2) = \det(x_1 v_1 + x_2 v_2). 
\end{equation}
This is a homogeneous polynomial  of degree $3$
of $v=(v_1,v_2)$. 
We define an action of $G$ on $\mathrm{Sym}^3 \Aff^2$ by
\begin{equation}
G \times \mathrm{Sym}^3\Aff^2
\ni \left((g_1,g_2), f(v)\right)
\mapsto c(g_1)f(vg_2)
\in \mathrm{Sym}^3 \Aff^2, 
\end{equation}
where we are regarding $v$ as a row vector. 
Then the map
$V \ni x \mapsto F_x \in  \mathrm{Sym}^3 \Aff^2$ 
is $G$-equivariant.

For $f \in \mathrm{Sym}^3 \Aff^2$, we denote the discriminant of $f$ by
$\Delta(f)$. It is easy to see that $\Delta(F_x)$ 
is a homogeneous polynomial of degree $12$ in $k[V]$, 
and that
\begin{equation*}
\Delta(F_{gx}) = c(g_1)^4 (\det(g_2))^6\Delta(F_x).
\end{equation*}
So $(G, V)$ is an irreducible prehomogeneous vector space and 
$\Delta(F_x)$ is a relative invariant polynomial of degree $12$
of $x$.  We shall discuss the details in Section \ref{sec:stabilizer}.

Let
\begin{equation}
\label{eq:w-defn}
w = w_1v_1 + w_2 v_2 =
\begin{pmatrix}
1&0&0\\
0&-1&0\\
0&0&0
\end{pmatrix}
v_1 +
\begin{pmatrix}
0&0&0\\
0&1&0\\
0&0&-1
\end{pmatrix}
v_2
\ \in \J \otimes \mathrm{Aff}^2.
\end{equation}
Let $G_w$ denote the stabilizer of $w$ in $G$ 
and $H_w = G_w \cap H$.

Suppose that $(G,V)$ is an arbitrary
irreducible prehomogeneous vector space over $k$.
If there exists $w \in \Vss_k$ such that  
$G_w$ is reductive (we are assuming that reductive groups are smooth
over the ground field), then $(G,V)$ is said to be {\it regular}. 
If $(G,V)$ is irreducible and regular, it is known that
$V_{\ks}^{\mathrm{ss}}$ is a single $G_{\ks}$-orbit (see \cite{regularity-of-PVS}).
Note that it is proved in \cite{regularity-of-PVS} that
$V_{\ks}\setminus G_{\ks} w$ is a hypersurface. However,
since the representation is irreducible, the $k$-algebra
generated by relative invariant polynomials has a single generator,
say $\Delta(x)$.  Then $V_{\ks}\setminus G_{\ks} w$ must coincide 
with $\{x\in V_{\ks}\,|\, \Delta(x)=0\}$. Therefore,
$\{x\in V_{\ks}\,|\,\Delta(x)\not=0\}=G_{\ks}w$ and so 
$V^{\mathrm{ss}}_{\ks}$ is a single $G_{\ks}$-orbit. 

This implies that for any $x \in V_{\ks}^{\mathrm{ss}}$, 
there  exists $g\in G_{\ks}$ such that 
$x = g w$. We define a map $c_x$ by
\begin{equation}
\label{eq:cx-defn}
c_x : \Gal \ni \sigma \mapsto g^{-1} g^{\sigma} \in G_{w\,\ks}.
\end{equation}
Then $c_x$ is a $1$-cocycle and defines an element
of the first Galois cohomology set
$\HGw$. For the definition of the Galois action and the Galois cohomology 
set, see Section \ref{sec:preliminaries}.
For arbitrary algebraic groups $G, H$ defined over $k$,
 a homomorphism $f: G \to H$ defined over $k$ and 
a map   $h_f^1 : \mathrm{H}^1(k, G) \to \mathrm{H}^1(k, H)$
which is induced by $f$,
we denote the inverse image of $1 \in  \mathrm{H}^1(k, H)$ 
under $h_f^1$  by 
Ker$\left(\mathrm{H}^1(k, G) \to \mathrm{H}^1(k, H)\right)$.
The following theorem is well-known (see \cite[pp.268,269]{Igusa} 
for example). 
\begin{thm}
\label{thm:galois-coh}
If $(G,V)$ is an irreducible prehomogeneous vector space, 
the orbit of $w\in V$ is open in $V$ and $G_w$ is reductive,
then the map
\begin{equation}
\ 
G_{k} \backslash V_k^{\mathrm{ss}}  \ni  x \mapsto c_x \in \KerHGw 
\end{equation}
is well-defined and bijective.
\end{thm}
Note that it is assumed in \cite{Igusa} that $\text{ch}(k)=0$. 
The proof of the above theorem works
as long as $V_{\ks}^{\mathrm{ss}}$ is a single $G_{\ks}$-orbit.
As we pointed out above, with the assumption of the theorem,
$V_{\ks}^{\mathrm{ss}}$ is a single $G_{\ks}$-orbit even if
the characteristic of $k$ is positive.
So the above theorem holds in all characteristics.

Now we go back to our situation and assume that 
$G = G_1 \times \GL(2)$, 
$V=\J\otimes \Aff^2$.  
When $k$ is the field of complex numbers,  
Sato and Kimura have shown in \cite[pp.138,139]{SatoKimura}
that $(G, V)$ is a prehomogeneous vector space having 
a relative invariant polynomial of degree $12$ 
and for a generic point $w \in V$,
the Lie algebra of $G_w$ is 
isomorphic to the Lie algebra of $\SO(Q)$.
However, the structure of $G_w$ has not been determined completely over 
an arbitrary ground field and the
interpretation of $G_k \backslash V_k^{\mathrm{ss}}$ is unknown.

In Section \ref{sec:stabilizer}, we determine the structure of the stabilizer as an algebraic 
group over $k$ and obtain the following proposition.
\begin{prop}
\label{prop:stabilizer}
$G_w$ is isomorphic to $\GL(1) \times (\Spin(Q) \rtimes \Sss)$.
In particular, $G_w$ is a smooth reductive algebraic group.
\end{prop}

Let $\text{JIC}(k)$ be the set of equivalence classes
of pairs $(\Mj,\an)$ where $\Mj$'s are isotopes of $\J$ 
and $\an$'s are cubic \'etale subalgebras of $\Mj$
(see Section \ref{sec:Rational Orbits} or \cite[pp.154--158]{Springer}
for the definition of isotopes of $\J$). 
For the details of the equivalence relation, see Section  
\ref{sec:Rational Orbits}.
In Section \ref{sec:Rational Orbits}, we associate a pair 
$(\Mj,\an)\in \text{JIC}(k)$ 
to each point in $V^{\mathrm{ss}}_k$. 
Our main theorem concerns a correspondence between 
the set $G_k\backslash V^{\mathrm{ss}}_k$ of rational 
orbits and the set of equivalence classes of pairs
$(\Mj,\an)$ as above. 
The following theorem is our main theorem.
\begin{thm}
\label{thm:main}
The set $G_k\backslash V^{\mathrm{ss}}_k$ corresponds 
bijectively with the set  
\begin{math}
\mathrm{JIC}(k).
\end{math}
\end{thm}

In Section \ref{sec:equivariant-map}, to each element of $V_{k}^{\mathrm{ss}}$, 
we associate an isotope of $\J$ 
and its cubic \' etale subalgebra
explicitly by constructing an equivariant map from $\J \otimes \Aff^2$
to $\J$. 

If $\alpha\in\overline k$ and $F_x(\alpha,1)=0$, 
we call $\alpha$ a root of $F_x(v)$. 
In Section \ref{sec:split-case}, 
we prove the following theorem.
\begin{thm}
\label{thm:main-split-intro}
If $k$ is a finite field or $\Oct$ is the split octonion, then 
there is a bijective correspondence
between $G_k\backslash V^{\rm{ss}}_k$ 
and  $\text{H}^1(k,\mathfrak S_3)$. Moreover, 
if $x\in V^{\rm{ss}}_k$, then the corresponding 
cohomology class in $\text{H}^1(k,\mathfrak S_3)$
is the element determined by the action of the Galois group 
on the set of roots of $F_x(v)$. 
\end{thm} 

If $n>0$, then it is well-known that 
$\mathfrak S_n$ is the automorphism group 
of the $k$-algebra $k^n$.  Therefore, 
the set $\mathrm{H}^1(k,\mathfrak S_n)$ can
be identified with $k$-isomorphism classes of \'etale $k$-algebras
of degree $n$. In the situation of Theorem \ref{thm:main-split-intro},
considering the case $n=3$, 
rational orbits in $V^{\rm{ss}}_k$ are in bijective correspondence with
separable extensions of $k$ of degrees up to $3$. 
Also if $x\in V^{\rm{ss}}_k$ corresponds to a separable cubic extension,
then the identity component $G_w^{\circ}$ is 
isomorphic to a $k$-form of $\Spin(Q)$ which comes from the triality
and so $G_k\backslash V^{\rm{ss}}_k$ parametrizes such objects.

The problem of finding arithmetic interpretations of rational orbits
of prehomogeneous vector spaces is classical and goes back to the
work of Gauss.  Some cases which are in some sense similar
to our case have been considered. For example, the prehomogeneous
vector spaces considered in \cite{wryu},
\cite{kayu}, \cite{taniguchi-orbits} have close relations
with field extensions of the ground field of degrees up to $5$.
Among the cases considered in these papers,
consider the following prehomogenous vector spaces:

(a) $G={\rm GL}(3)\times {\rm GL}(2)$, 
$V={\rm Sym}^2 {\rm Aff}^3 \otimes {\rm Aff}^2$ (see \cite{wryu}), 
 
(b) $G=\operatorname{Res}_{k_1/k}{\rm GL}(3)\times {\rm GL}(2)$, 
$V={\rm H}_3(k_1) \otimes {\rm Aff}^2$ (see \cite{kayu}), 

(c) $G={\rm GL}_3(D)\times {\rm GL}(2)$,
$V={\rm H}_3(D) \otimes {\rm Aff}^2$ (see \cite{taniguchi-orbits})

\noindent
where $k_1$ is a quadratic extension of $k$,
$D$ is a quaternion algebra over $k$ and
${\mathrm H}_3(k_1),{\mathrm H}_3(D)$ are the spaces of
$3\times 3$ Hermitian matrices with entries in $k_1,D$
respectively. These cases have the same set of weights
with respect to their maximal split tori. If the ground field
is $\R$, these cases and our case correspond to
$\R,\C,{\mathbb H}$ (the Hamiltonian), $\Oct$.  
So our case is the last of these series.

The problem of finding arithmetic interpretations of rational orbits
of our case over rings is more difficult. We may consider such problems
in the future. 

In Section \ref{sec:preliminaries}, we define and review basic notions
such as the Galois cohomology, the Lie algebra of groups of type
${\mathrm E}_6$ and the triality concerning ${\mathrm{Spin}}(8)$. 
In Section \ref{sec:stabilizer}, we determine the structure
of the stabilizer $G_w$. 
In Section \ref{sec:Rational Orbits}, we
give an intrinsic arithmetic interpretation
of the set of rational orbits $G_k\backslash V^{\mathrm{ss}}_k$. 
In Section \ref{sec:equivariant-map}, we construct
the isotopes corresponding to points in $V^{\mathrm{ss}}_k$
explicitly by using an equivariant map from $V$ to $\J$.  
In Section \ref{sec:split-case}, we prove Theorem
\ref{thm:main-split-intro}. 


\section{Preliminaries}
\label{sec:preliminaries}

In this section, we define basic notations and prove  
fundamental facts used in the subsequent sections.

We define an action of the Galois group on an algebraic group
as follows.
Suppose that $K$ is a Galois extension of $k$.
Then we define the multiplication 
in the Galois group $\mathrm{Gal}(K/k)$ by 
\begin{equation}
\label{eq:galois-product}
(\sigma\tau)(a) = \tau\left(\sigma(a)\right).
\end{equation}
for $\sigma, \tau \in \mathrm{Gal}(K/k)$ and $a \in K$.
We denote the action of $\sigma\in \mathrm{Gal}(K/k)$
on $x\in K$ by $x^{\sigma}$. This is a right action.

Let $G$ be an arbitrary algebraic group over $k$.
The Lie algebra of $G$ is denoted by $\text{Lie}(G)$. 
We denote the identity component (the connected component
containing the unit element) by $G^{\circ}$. 
If $G$ acts on a variety $X$ and $x\in X$, 
then we denote the stabilizer and the orbit of $x$ 
by $G_x$, $G(x)$ respectively. 

If $F$ is any extension field of $k$, then
$G$ can be regarded as an algebraic group over $F$. 
If $R$ is a $k$-algebra, we denote the
set of $R$-rational points by $G_R$. 
If $K/k$ is a field extension, then
we denote the $K$-algebra of regular functions
in the sense of \cite[p.15, 6.3]{Borel} on
$G_K$ by $K[G]$. Then, $\ks [G] \cong \ks \otimes_{k} k[G]$.
We define a right action of $\Gal$ on $\ks[G]$ by 
\begin{equation*}
\nu(\sigma): 
\ks \otimes_{k} k[G]
\ni \sum a_i \otimes v_i
\mapsto
\sum a_i^{\sigma} \otimes v_i
\in \ks \otimes_{k} k[G].
\end{equation*} 

For any $p \in G_{\ks}$, let $p^{\#}$ be the
associated $\ks$-algebra 
homomorphism from $\ks[G]$ to $\ks$.  Then
$\sigma \circ p^{\#}\circ \nu(\sigma)^{-1}$ is also 
a $\ks$-algebra homomorphism. We denote the element in 
$G_{\ks}$ corresponding to 
$\sigma \circ p^{\#}\circ \nu(\sigma)^{-1}$ 
by $p^{\sigma}$.
Then, $p^{\sigma \tau} = (p^{\sigma})^{\tau}$ for any $\sigma, \tau \in \Gal$.
Therefore, we can define a right action of $\Gal$ on $G_{\ks}$ by
\begin{equation*}
G_{\ks} \times \Gal
\ni (p, \sigma)
\mapsto
p^{\sigma}\in G_{\ks}.    
\end{equation*} 
%


We next define notations regarding the Galois cohomology.
\begin{mydef}
\label{def:galcoh}
For an algebraic group $G$ over $k$, we consider the discrete topology
on $G_{\ks}$.
A continuous function $h : \Gal \to G_{\ks}$ is called a {\it $1$-cocycle} if 
$h(\sigma\tau) = h(\tau)h(\sigma)^{\tau}$ for any $\sigma, \tau \in \Gal$.
Two 1-cocycles $h$ and $h'$ are equivalent if there exists $g \in G_{\ks}$
such that 
\begin{equation*}
h(\sigma) = g^{-1} h'(\sigma) g^{\sigma}
\end{equation*}
for any $\sigma \in \Gal$.
The {\it first Galois cohomology set} $\mathrm{H}^1(k, G)$
is the set of equivalent classes of $1$-cocycles 
by the above equivalence relation. 
\end{mydef}

We shall define several notions regarding 
the exceptional Jordan algebra $\J$ (see (\ref{eq:X-defn}))
and review its fundamental properties. 
For the rest of this paper $G,G_1,H,H_1$ are the groups
defined in Introduction.  

Let $I_n$ be the $n \times n$ unit matrix.
We denote the diagonal matrix with diagonal entries
$\alpha_1,\ldots,\alpha_n$ by
$\text{diag}(\alpha_1,\ldots,\alpha_n)$.

Let $\Oct$ be an octonion and $\J$
the exceptional Jordan algebra as in 
Introduction.  Let $\Oct^{\times}=\{x\in \Oct\,|\, \|x\|\not=0\}$.
This set is closed under multiplication, but may not be
a group since the multiplication is not associative.
However, if $x\in \Oct^{\times}$, then $x^{-1}$ exists in 
$\Oct$. 
If $X,Y\in\J$, then we denote the usual 
matrix multiplication by $XY$. If $X,Y,Z\in\J$
and $(XY)Z=X(YZ)$, we may write $XYZ$.

For $X \in \J$, we define an endomorphism $R_X \in \mathrm{End}(\J)$ by 
\begin{equation*}
R_X(W) = W \circ X = \frac{1}{2}(WX + XW)
\end{equation*}
for $W \in \J$.   
We define a symmetric trilinear form 
$D$ on $\J$ by
\begin{equation*}
\begin{split}
6D(X, Y , Z)  = 
&\det(X + Y  + Z) - \det(X + Y) - \det(Y  + Z) - \det(Z + X) \\
&+ \det(X) + \det(Y) +\det(Z).
\end{split}
\end{equation*}
We denote the trace (the sum of diagonal entries) on $\J$  by $\Trace$. 
For $X, Y \in \J$, 
we define a symmetric bilinear form $\langle \ ,\ \rangle $ on $\J$ by
\begin{equation*}
\langle  X, Y\rangle  = \Trace(X\circ Y).
\end{equation*}
One can verify by direct computation that the symmetric bilinear form 
$\langle \ ,\ \rangle $ satisfies the following equation:
\begin{equation*}
\label{bi-formula}
\langle X \circ Y , Z \rangle  = \langle X, Y \circ Z\rangle, 
\quad {}^{\forall} X, Y, Z \in \J. 
\end{equation*}
%

For $X, Y \in \J$,  $X \times Y$ is, by definition,
the element satisfying 
the following equation:
\begin{equation*}
\langle  X \times Y,  Z \rangle  = 3D(X, Y, Z),
\quad {}^{\forall} Z \in \J. 
\end{equation*}
Let $e=I_3$. Then the following equations are satisfied 
(see \cite[p.122, Lemma 5.2.1]{Springer}). 
\begin{equation}
\label{cross-formula}
{}^{\forall} X\in \J, \quad X \circ (X\times X) = \det(X)e, 
\quad e \times e = e.
\end{equation}

Let
\begin{equation}
\label{elements-defn}
\begin{split}
&E_1 =
\begin{pmatrix}
1&0&0\\
0&0&0\\
0&0&0
\end{pmatrix}
, \ \ \ \ \ 
(\al)_1 = 
\begin{pmatrix}
0&0&0\\
0&0&\al \\
0&\oline{\al}&0 
\end{pmatrix},\\
&E_2 =
\begin{pmatrix}
0&0&0\\
0&1&0\\
0&0&0
\end{pmatrix}
, \ \ \ \ \ 
(\be)_2 = 
\begin{pmatrix}
0&0&\oline{\be}\\
0&0&0 \\
\be&0&0 
\end{pmatrix},\\
&E_3 =
\begin{pmatrix}
0&0&0\\
0&0&0\\
0&0&1
\end{pmatrix}
, \ \ \ \ \ 
(\ga)_3 = 
\begin{pmatrix}
0&\ga&0\\
\oline{\ga}&0&0 \\
0&0&0 
\end{pmatrix}
\end{split}
\end{equation}
for $\al, \be, \ga \in \Oct$. 
Then,
\begin{equation}
\label{cal}
\begin{array}{ccc}
(\al)_1 \circ (E_i) =
\begin{cases}
\vspace{5pt}
0&i = 1\\
\vspace{5pt}
\dispfr{1}{2}(\al)_1&i = 2\\
\dispfr{1}{2}(\al)_1&i = 3
\end{cases},
&
(\be)_2 \circ (E_i) =
\begin{cases}
\vspace{5pt}
\dispfr{1}{2}(\be)_2		&i = 1\\
\vspace{5pt}
0						&i = 2\\
\dispfr{1}{2}(\be)_2	&i = 3
\end{cases},\\
\\
(\ga)_3 \circ (E_i) =
\begin{cases}
\vspace{5pt}
\dispfr{1}{2}(\ga)_3		&i = 1\\
\vspace{5pt}
\dispfr{1}{2}(\ga)_3		&i = 2\\
0						&i = 3
\end{cases}.
&
\end{array}
\end{equation}

We define an injective $k$-linear map $\iota_i : \Oct \to \J$ 
for $i = 1, 2,3$ by
\begin{align}
\label{eq:oct-i-defn}
\iota_i(x) = (x)_i, \quad { }^{\forall} x \in \Oct.
\end{align}
We denote $\iota_i(\Oct) \subset \J$ by $\Oct_ i$.

The algebraic group $F_4$ is defined as follows:
\begin{equation}
\label{eq:F4-defn}
F_4 = \Aut(\J) = \{L \in G_1 \,|\, L(X) \circ L(Y) = L(X \circ Y ), \quad
{}^{\forall}     X, Y \in \J \}.
\end{equation}
It is known that
\begin{equation}
\label{eq:F4-e}
F_4 = \{L \in G_1 \,|\, L(e) = e \}
\end{equation}
(see \cite[p.159, Proposition 5.9.4]{Springer}) and 
$F_4$ is a connected simple algebraic group of type ${\mathrm F}_4$
which is defined over $k$ (see \cite[p.178, Theorem 7.2.1]{Springer}).
Note that if $g\in G_1$ and $g(e)=e$, then by taking the determinant,
$c(g)=1$ and so $g\in H_1$.

Let $\Lh_1$ and $\Laut$ be 
the Lie algebras of $H_1$ and $F_4$ respectively. 
Let $k[\varepsilon]/(\varepsilon^2)$ be the ring of dual numbers.
Since ch$(k) \neq 3$, it is easy to see that
\begin{equation}
\label{eq:Lie-AutJ}  
\begin{aligned}
\Lh_1  & = \{ t \in  \Lgl(\J)  \, | \,  
\det \left((1+\varepsilon t)(X)\right)
= \det(X), \quad {}^{\forall} X \in \J\} \\ 
& =\{ t \in \Lgl(\J) \,|\, D(t(X), X, X)= 0, \quad {}^{\forall} X \in \J\},  \\ 
\Laut  & = \{ t \in  \Lh_1 \,|\,  
(1+\varepsilon t)(e) = e \} 
=\{ t \in \Lh_1 \,|\, t(e) = 0 \}.
\end{aligned}
\end{equation}
Moreover, it is known that dim $\Lh_1 = 78$ and dim $\Laut = 52$
(see \cite[p.181, Theorem 7.3.2]{Springer} 
and  \cite[p.180,Corollary 7.2.2]{Springer}). 

Let $\mathrm{Der}_ {k} (\J, \J)$ be the Lie algebra of
$k$-derivations of $\J$.
The following fact is known (see \cite[p.180, Corollary 7.2.2]{Springer}). 
\begin{lem}
\label{lem:52}  
\begin{math}
\Laut = \Der
\end{math} 
and $\mathrm{dim}_k \Der = 52$.  
 \end{lem}

We define a vector subspace $\mathfrak{D}_0  \subset  \Le$  by
\begin{equation*}
\mathfrak{D}_0 = \{ X \in \mathrm{Der}_ {k} (\J, \J) 
\,|\, X(E_i) = 0, i = 1, 2, 3\}.
\end{equation*}
Let $\Lso(Q)$
denote the Lie algebra of $\SO(Q)$. 
Then
\begin{equation*}
\Lso(Q) = \{ t \in \Lgl(\Oct) \,|\, Q( t(x), y) +  
Q( x, t(y) ) = 0, \quad {}^{\forall} x, y \in \Oct \}. 
\end{equation*} 
If $L:\Oct\to \Oct$ is a $k$-linear map, then
we define $\widehat L=\iota \, L \,\iota$ (see (\ref{eq:iota-defn})). 
Obviously, if $L\in\SO(Q),\Lso(Q)$, then
$\widehat L\in \SO(Q),\Lso(Q)$ respectively. 
 
The following proposition is probably known
(if $k=\C$, it is proved in
\cite[pp.24--27]{SatoKimura}).  
However, we could not find 
a good reference over arbitrary fields of 
characteristic not equal to $2,3$ and so
we include the proof.  
\begin{prop}
\label{prop:Lie-alg-G1}
\begin{itemize}
\item[(1)] 
\begin{math}
\Lh_1 = \Der \oplus \{R_Y \,|\, \Trace(Y) = 0\}.
\end{math}
\item[(2)] 
\begin{math}
\Der \cong \mathfrak{D}_0 \oplus \Oct \oplus \Oct \oplus \Oct. 
\end{math}
\end{itemize}
\end{prop}
\begin{proof}
We shall prove this proposition in the
following order. 

(I) The proof of 
$\mathfrak{D}_0 \cong \Lso(Q)$. 

(II) The proof of (2).
 
(III) The proof of (1).


(I)  
We use the following theorem (see \cite[p.53, Theorem 3.5.5]{Springer}
and {\cite[p.56, Lemma 3.5.9]{Springer}}).
\begin{thm}
\label{Pri-Loc-Tri}
If $\mathrm{ch}(k) \neq 2$, then for any $t_1 \in 
\Lso(Q)$,
there exist unique $t_2, t_3 \in \Lso(Q)$ such that 
\begin{equation}
\label{Local-tri}
t_1(xy) = t_2(x)y + xt_3(y) 
\end{equation} 
for any $x, y \in \Oct$.
If $(t_1, t_2, t_3)$ satisfies the equation $\eqref{Local-tri}$, then
$(t_2, t_1,\widehat t_3)$ and $(t_3, \widehat t_2, t_1)$ also satisfy 
the equation $\eqref{Local-tri}$.
\end{thm}

For any $D \in \mathfrak{D}_0$, by the equations of $\eqref{cal}$,
\begin{equation*}
D((\al)_1 \circ (E_i)) =  D((\al)_1)\circ E_i =
\begin{cases}
\vspace{5pt}
0&i = 1\\
\vspace{5pt}
\dispfr{1}{2}D((\al)_1)&i = 2\\
\dispfr{1}{2}D((\al)_1)&i = 3
\end{cases}.
\end{equation*} 
It follows immediately that $D(\Oct_1) \subset \Oct_1$.  Similarly, we
have $D(\Oct_i) \subset \Oct_i$ for $i = 2, 3$.
For $\alpha\in \Oct$, we define $D_1(\alpha)$ to be the element
of $\Oct$ such that $D((\al)_1)=(D_1(\al))_1$.
We define $D_2,D_3$ similarly.

Since $(\al)_1^2 = \|\al\|(E_2 + E_3)$ and
\begin{align*}
D\left((\al)_1^2\right) & = \|\al\|D(E_2 + E_3) = 0, \\
D\left((\al)_1^2\right) & =
2(D_1 \left(\al\right))_1 \circ (\al)_1 =
2Q(D_1(\al), \al)(E_2 + E_3), 
\end{align*}   
we have $Q(D_1(\al), \al) = 0$. 
It follows that $D_1 \in \Lso(Q)$.
Similarly, we have $D_2, D_3 \in \Lso(Q)$.
We define a $k$-linear map $r$ by
\begin{equation*}
r : \mathfrak{D}_0\ni D \mapsto \widehat D_3
\in \Lso(Q). 
\end{equation*}
We prove that $r$ is bijective.

Suppose that $t_1 \in \Lso(Q)$. 
We show that $t_1$ is in the image of $r$. 
There exist unique $t_2, t_3$ satisfying
$\eqref{Local-tri}$ by Theorem \ref{Pri-Loc-Tri}. 
We define a linear map 
$D_{t_1}: \J \to \J$ by
\begin{equation*}
D_{t_1}:
X = \begin{pmatrix}
s_1&x_3&\oline{x_2}\\
\oline{x_3}&s_2&x_1\\
x_2&\oline{x_1}&s_3
\end{pmatrix}
\mapsto 
\begin{pmatrix}
0&\widehat t_1(x_3)&\oline{t_3(x_2)}\\
\oline{\widehat t_1(x_3)}&0&t_2(x_1)\\
t_3(x_2)&\oline{t_2(x_1)}&0
\end{pmatrix}.
\end{equation*}
By Theorem \ref{Pri-Loc-Tri}, 
\begin{equation*}
\begin{array}{rcl}
\widehat t_1(\oline{x_1 x_2})&=&\oline{t_2(x_1)x_2 + x_1 t_3(x_2)}, \\
t_3(\oline{x_3 x_1})&=&\oline{\widehat t_1 (x_3)x_1 + x_3t_2(x_1)}, \\
t_2(\oline{x_2x_3})&=&\oline{t_3(x_2)x_3 + x_2 \widehat t_1 (x_3)}.
\end{array}
\end{equation*}

By computation, 
\begin{equation*}
X^2   
= \begin{pmatrix}
s_1^2 + \|x_2\| + \|x_3\|
& (s_1+s_2)x_3 + \oline{x_2} \, \oline{x_1}
& (s_1+s_3)\oline{x_2} + x_3x_1 \\
(s_1+s_2)\oline{x_3} + x_1x_2 
& s_2^2 + \|x_1\| + \|x_3\|
& (s_2+s_3)x_1 + \oline{x_3} \, \oline{x_2} \\
(s_1+s_3)x_2 + \oline{x_1} \, \oline{x_3}
& (s_2+s_3)\oline{x_1} + x_2x_3 
& s_3^2 + \|x_1\| + \|x_2\|
\end{pmatrix}.
\end{equation*}
Therefore, 
\begin{equation*}
\begin{array}{l}
D_{t_1}(X^2) = \\[7pt]
\!\!\!
\begin{pmatrix}
\!0\!&\!\!(s_1 + s_2)\widehat t_1(x_3) + \widehat t_1(\oline{x_1 x_2})\!&\!
(s_1 + s_3)\oline{t_3(x_2)} + \oline{t_3(\oline{x_3 x_1})}\\
\!\!(s_1 + s_2)\oline{\widehat t_1(x_3)} 
+ \oline{\widehat t_1(\oline{x_1 x_2})}
\!\!&\!0\!\!&\!(s_2 + s_3)t_2(x_1) + t_2(\oline{x_2x_3})\\
\!(s_1 + s_3)t_3(x_2) + t_3(\oline{x_3 x_1})
\!\!&\!\!(s_2 + s_3)\oline{t_2(x_1)} + \oline{t_2(\oline{x_2x_3})}\!&\!\!0
\end{pmatrix} \vspace{5pt} \\
= 2D_{t_1}(X) \circ X 
\end{array}
\vspace{-4pt}
\end{equation*}
for any $X \in \J$.
This implies that $D(X\circ Y)=D(X)\circ Y+X\circ D(Y)$
for all $X,Y\in \J$
by considering $D((X+Y)\circ (X+Y))$. 
Thus, $D_{t_1} \in \mathfrak{D}_0$
and so $r(D_{t_1})=t_1$. Therefore,  
$r$ is surjective.

Suppose that $D\in \mathfrak{D}_0$. 
Since $\left(\,\overline{\al\be}\,\right)_3=2(\al)_1\circ (\be)_2$, 
we have
\begin{equation*}
\begin{array}{rcl}
\vspace{2pt}
D\left(\left(\,\overline {\al\be}\,\right)_3\right)
&=& 2D((\al)_1 \circ (\be)_2) =
2\left((D_1(\al))_1 \circ (\be)_2 + (\al)_1 \circ 
(D_2(\be))_2\right)  \\
\vspace{2pt}
&=&  \left(\,\overline{D_1 (\al)\be + \al D_2(\be)}\,\right)_3
\end{array}
\end{equation*}   
for any $\al, \be \in \Oct$.
So $\widehat{D}_3(\al \be) = D_1 (\al)\be + \al D_2(\be)$. 
Therefore, $(\widehat D_3, D_1, D_2)$ satisfies the equation
$\eqref{Local-tri}$.
So the uniqueness property of Theorem
\ref{Pri-Loc-Tri} implies that  
$r$ is injective. 
Therefore, we have $\mathfrak{D}_0 \cong \Lso(Q)$ and 
dim $\mathfrak{D}_0 = 28$.

(II)
We next investigate the structure of $\Der$ in detail.
For $\al,  \be, \ga \in \Oct$, we define elements
$(\al)'_1,(\be)'_2,(\ga)'_3$ of $\Le$ by 
\begin{align*}
(\al)'_1 & = [R_{E_2}, R_{(\al)_1}] \, 
(=  R_{E_2}R_{(\al)_1} - R_{(\al)_1}R_{E_2}), \\
(\be)'_2 & = [R_{E_3}, R_{(\be)_2}], \\
(\ga)'_3 & = [R_{E_1}, R_{(\ga)_3}].
\end{align*}
We show that $(\al)'_1 \in \mathrm{Der}_k(\J, \J)$.

For $X \in \J$ as in (\ref{eq:X-defn}), 
\begin{equation}
\label{eq:al1'}
(\al)'_1(X) = \dispfr{1}{4}
\begin{pmatrix}
0&\oline{\al x_2}&-x_3 \al\\
\al x_2 & 2Q(\al, x_1)& (-s_2 + s_3)\al\\
-\oline{x_3 \al}&(-s_2 + s_3)\oline{\al}& -2Q(\al, x_1)
\end{pmatrix}. 
\end{equation} 
Note that the following equation is satisfied
(see \cite[p.8, Lemma 1.3.3]{Springer}).
\begin{equation}
\label{eq:associator}
(xy)\oline{z} + (xz)\oline{y} = 2Q(y, z)x,  \quad {}^{\forall} x, y, z \in \Oct.
\end{equation} 
We have
\begin{equation*}
\begin{array}{l}
2(\al)'_1(X) \circ X = \\[7pt] \dispfr{1}{4} \!\!
\vspace{3pt}
\begin{pmatrix}
\!0&\!\!(s_1 + s_3)\oline{\al x_2} + (x_3x_1)\oline{\al}
&\!\!-(s_1 + s_2)x_3\al - \oline{x_1x_2}\al\\
\!(s_1 + s_3)\al x_2 + \al \overline{x_3x_1}
&\!\!2Q((s_2 + s_3)x_1 + \oline{x_2x_3},\al)
&\!\!(s_3^2 - \|x_3\| -s_2^2 + \|x_2\|)\al \\
\!-(s_1 + s_2)\oline{x_3\al} - \oline{\oline{x_1x_2}\al}
&\!\!(s_3^2 - \|x_3\| -s_2^2 + \|x_2\|)\oline{\al}
&\!\!-2Q((s_2 + s_3)x_1 + \oline{x_2x_3},\al)
\end{pmatrix} \vspace{5pt}\\
\vspace{3pt}
=  (\al)'_1(X^2)
\end{array}
\vspace{-4pt}
\end{equation*}
((\ref{eq:associator}) is used for the computations of
$(i,j)$-entries ($i\not=j$)). 
So, $(\al)'_1 \in \mathrm{Der}_k(\J, \J)$.

Similarly, we have 
$(\be)'_2, (\ga)'_3 \in \mathrm{Der}_k(\J, \J)$.
Also by (\ref{eq:al1'}) and similar calculations, we have
\begin{equation}
\label{eq:albegam'} 
\begin{array}{ccc}
(\al)'_1 (E_i) =
\begin{cases}
\vspace{5pt}
0							&i = 1\\
-\dispfr{1}{4}(\al)_1		&i = 2\\
\vspace{5pt}
\dispfr{1}{4}(\al)_1			&i = 3
\end{cases},
&
(\be)'_2  (E_i) =
\begin{cases}
\vspace{5pt}
\dispfr{1}{4}(\be)_2		&i = 1\\
\vspace{5pt}
0						&i = 2\\
- \dispfr{1}{4}(\be)_2	&i = 3
\end{cases},\\
\\
(\ga)'_3  (E_i) =
\begin{cases}
\vspace{5pt}
-\dispfr{1}{4}(\ga)_3 & i = 1\\
\vspace{5pt}
\dispfr{1}{4}(\ga)_3	 & i = 2\\
0 & i = 3
\end{cases}.
&
\end{array}
\end{equation}

We define a $k$-linear map $l$ by
\begin{equation*}
\begin{array}{cccc}
l:&\mathfrak{D}_0 \oplus \Oct \oplus \Oct \oplus \Oct &\to&
\mathrm{Der}_k(\J,\J) \\
&\rotatebox{90}{$\in $}&&\rotatebox{90}{$\in$} \\
&(X_0, \al, \be, \ga) &\mapsto&X_0 + (\al)'_1 +  (\be)'_2 +  (\ga)'_3 
\end{array}.
\end{equation*}
We show that $l$ is bijective.  For that purpose,
it is enough to prove that $l$ is surjective since
the dimensions of both sides are equal to $52$
(see Lemma \ref{lem:52} and (I)).
 
For $D \in \Der$, let $X_i = D(E_i) = \hji \in \J$ where $i = 1, 2, 3$
(see (\ref{eq:X-defn})). 
Then,
\begin{equation*}
D(E_i) = D(E_i^2) = 2E_i \circ D(E_i).
\end{equation*}
We consider the case $i = 1$. Since
\begin{align*}
X_1 = D(E_1) =2E_1 \circ D(E_1) =
\begin{pmatrix}
2s_{11}&x_{13}&\oline{x_{12}}\\
\oline{x_{13}}&0&0 \\
x_{12}&0&0
\end{pmatrix},
\end{align*}
we have $s_{11}=s_{12}=s_{13}=x_{11}=0$.
Therefore,  
\begin{math}
X_1 = (x_{12})_2 + (x_{13})_3.
\end{math}
By similar calculations,
we have $s_{2i}=s_{3i}=0$ for $i=1,2,3$,
$x_{22}=x_{33}=0$ and 
\begin{math}
X_2 = (x_{23})_3 + (x_{21})_1, \;
X_3 = (x_{31})_1 + (x_{32})_2.
\end{math}
Moreover, since
\begin{equation*}
0 = D(E_i\circ E_j) = D(E_i) \circ E_j + E_i \circ D(E_j) = 
X_i \circ E_j + E_i \circ X_j 
\end{equation*}
for any $i\not=j \in \{1, 2, 3\}$,
\begin{equation}
\label{eq:x12etc}
\begin{aligned}
(x_{13})_3 + (x_{23})_3 = 0  \quad \text{from}\; (i,j)=(1,2), \\
(x_{12})_2 + (x_{32})_2 = 0  \quad \text{from}\; (i,j)=(1,3), \\
(x_{21})_1 + (x_{31})_1 = 0  \quad \text{from}\; (i,j)=(2,3). 
\end{aligned}
\end{equation}
Hence, using (\ref{eq:albegam'}), (\ref{eq:x12etc}),
we have
\begin{equation*}
\left(D - 4(x_{31})'_1 - 4(x_{12})'_2 - 4(x_{23})'_3\right)(E_i) = 0.
\end{equation*}

Let
\begin{equation*}
D_0 = D - 4(x_{31})'_1 - 4(x_{12})'_2 - 4(x_{23})'_3.
\end{equation*}
Then $D_0 \in \mathfrak{D}_0$ and
\begin{equation*}
l((D_0, 4x_{31}, 4x_{12},  4x_{23})) = D.
\end{equation*}
Therefore, $l$ is surjective. So we conclude that $\Der \cong 
\mathfrak{D}_0 \oplus \Oct \oplus \Oct \oplus \Oct$.

(III)
Finally, we prove (1). 
For any $R_Y$ where $\Trace(Y) = \langle Y, e \rangle = 0$,
\begin{align*}
D(R_Y(X), X, X)  
& = \frac 13 \langle Y \circ X, X \times X\rangle  
= \dispfr{1}{3} \langle Y, X \circ (X\times X)\rangle \\
& = \dispfr{1}{3}\det(X)\langle Y, e\rangle  = 0
\end{align*}
for any $X \in \J$.
Thus, $R_Y \in \Le$.

For $t \in \Le$,
let $Y = t(e)$. Then,
since $D(t(e),e,e) = 0$, by $\eqref{cross-formula}$, we have
\begin{equation*}
D(t(e),e,e)  = \dispfr{1}{3}\langle t(e), 
e \times e\rangle  = \dispfr{1}{3}\langle Y, e\rangle = 0,
\end{equation*}
It follows that Tr$(Y) = 0$.
Since $(t - R_Y)(e) = 0$, by 
$\eqref{eq:Lie-AutJ}$ , 
$t - R_Y \in \mathrm{Der}_{k} (\J, \J)$.
Thus,
\begin{equation*}
\Lh_1 = \Der + \{R_Y \,|\, \Trace(Y) = 0\}.
\end{equation*}
Since $\dim \Lh_1=78=\dim \Der + \dim \{R_Y \,|\, \Trace(Y) = 0\}$,  
the above sum is a direct sum. 
This completes the proof of (1).
\end{proof}

We define varieties $M$ and $RT(\Oct)$ by 
\begin{equation*}
 M \! = \! \{(A, B, C) \in \SO(Q) \times \SO(Q) \times \SO(Q)
\,|\, A(x)B(y) = \widehat C(xy), {}^{\forall} x, y \in \Oct \},
\end{equation*}
\begin{equation*}
 RT(\Oct) \! = \! \{(A, B, C) \in \SO(Q) \times \SO(Q) \times \SO(Q)
\,|\, A(xy) = B(x)C(y), {}^{\forall} x, y \in \Oct \}
\end{equation*}
respectively.
Then $M$ and $RT(\Oct)$  are defined over $k$.

The variety $RT(\Oct)$ is given in \cite[p.59]{Springer}.
It is known that $RT(\Oct)$ is a closed subgroup of the algebraic group 
$\SO(Q)\times \SO(Q)\times \SO(Q)$ (see \cite[p.59]{Springer}).
We show that $M$ is also a
subgroup of $\SO(Q)\times \SO(Q)\times \SO(Q)$
and is isomorphic to $RT(\Oct)$ over $k$.

For any $(A, B, C), (A',B',C') \in M$ and $x, y \in \Oct$, we have
\begin{align*}
(A'A)(x)(B'B)(y) 
& = A'(A(x))B'(B(y)) 
= \widehat C'(A(x)B(y)) \\ 
& = (\widehat C' \widehat C)(xy)
= (\widehat {C'C})(xy), \\
(A^{-1})(x)(B^{-1})(y) 
& = (\widehat C^{-1} \widehat C)(A^{-1}(x)B^{-1}(y)) \\
& = \widehat C^{-1}(A(A^{-1}(x))(B(B^{-1}(y)))) \\
& = \widehat C^{-1}(xy) 
= \widehat {C^{-1}}(xy).
\end{align*}
It follows that $M$ is closed under multiplication 
and taking the inverse. 
Obviously, the unit element $(I_8, I_8, I_8)$ of 
$\SO(Q)\times \SO(Q)\times \SO(Q)$ belongs to $M$.
Therefore, $M$ is a subgroup of $\SO(Q)\times \SO(Q)\times \SO(Q)$.

We define a homomorphism $f$ by
\begin{equation*}
f : M \ni (A, B, C)
\mapsto (\widehat{C}, A, B) \in RT(\Oct). 
\end{equation*}
Then, $f$ is obviously an isomorphism defined over $k$. 
So $M$ is isomorphic to $RT(\Oct)$ over $k$. 
Since $RT(\Oct)$ is smooth over $k$ (see \cite[p.60]{Springer}), 
$M$ is also smooth over $k$.

It is known that the algebraic group $RT(\Oct)$ satisfies
the following properties. Note that we are assuming 
that $\mathrm{ch}(k) \neq 2$. 
\begin{prop}
\label{Prop-RT}
\ \!
\begin{enumerate}
\item Let ${\mathrm{pr}}_i : RT(\Oct) \to \SO(Q)$ be the projection
\begin{math}
{\mathrm{pr}}_i((t_1,t_2,t_3)) = t_i
\end{math}
($i=1,2,3$).
Then ${\mathrm{pr}}_1$ is a surjective homomorphism from $RT(\Oct)$ to $\SO(Q)$ 
with kernel of order 2  and
the representations ${\mathrm{pr}}_i$ on $\Oct$ are irreducible and 
pairwise inequivalent.
\item
Let $\eta_1, \eta_2, \eta_3$ be the following
automorphisms of $RT(\Oct)$
\begin{equation*}
\begin{array}{rcl}
\eta_1 &:&(t_1, t_2, t_3) \to
(\widehat t_1, \widehat t_3, \widehat t_2) \\
\eta_2&:&(t_1, t_2, t_3) \to (t_3,\widehat  t_2, t_1) \\
\eta_3&:&(t_1, t_2, t_3) \to (t_2, t_1, \widehat t_3)
\end{array}.
\end{equation*} 
Then $\eta_i^2=1$ $(i = 1, 2, 3)$ and 
$\{\eta_1,\eta_2,\eta_3\}$  generates a group $\tilde S$ of outer 
automorphisms of $RT(\Oct)$ which is isomorphic to  $\Sss$.
\item
$RT(\Oct)$ is isomorphic to $\Spin(Q)$ as algebraic groups. 
\end{enumerate}
 \end{prop}
\begin{proof}
See \cite[p.59, Proposition 3.6.1]{Springer}, 
\cite[p.60, Proposition 3.6.3]{Springer} and
\cite[p.64, Proposition 3.6.6,]{Springer}.
\end{proof}
For
the definition of the algebraic group $\Spin(Q)$,  
see \cite[pp.38,39]{Springer}.
It is known that $\Spin(Q)$ is connected (see \cite[p.40]{Springer}). 
Moreover, there exists a surjective homomorphism from 
$\Spin(Q)$ to $\SO(Q)$ with kernel of 
order $2$ if ch$(k) \neq 2$ (see  \cite[p.40]{Springer}). 
Since $\SO(Q)$ is semi-simple, by \cite[p.192, 14.11 Corollary]{Borel},
$\Spin(Q)$ is semi-simple.
So $RT(\Oct)$ is also connected and semi-simple.

Since $M$ is isomorphic to $RT(\Oct)$, 
we have the following proposition
(we are still assuming that $\mathrm{ch}(k) \neq 2$). 
\begin{prop}
\label{Prop-M}
\
\begin{enumerate}
\item Let ${\mathrm{pr}}_i : M \to \SO(Q)$ be the projection
\begin{math}
{\mathrm{pr}}_i((t_1,t_2,t_3)) = t_i
\end{math}
($i=1,2,3$).
Then ${\mathrm{pr}}_3$ is a surjective homomorphism from $M$ to $\SO(Q)$ 
with kernel of order 2  and
the representations ${\mathrm{pr}}_i$ on $\Oct$ are irreducible and 
pairwise inequivalent.
\item
Let $\eta'_1, \eta'_2, \eta'_3$ be the following
automorphisms of $M$
\begin{equation*}
\begin{array}{rcl}
\eta'_1 &:&(t_1, t_2, t_3) \to (\widehat t_1, \widehat t_3, \widehat t_2) \\
\eta'_2&:&(t_1, t_2, t_3) \to (\widehat t_3,\widehat  t_2, \widehat t_1) \\
\eta'_3&:&(t_1, t_2, t_3) \to (\widehat t_2, \widehat t_1, \widehat t_3)
\end{array}.
\end{equation*} 
Then $\eta'_i{}^2=1$ $(i = 1, 2, 3)$ and
$\{\eta'_1,\eta'_2,\eta'_3\}$
generates a group of outer 
automorphisms $\tilde S'$ of $M$ which is isomorphic to  $\Sss$.
\item $M$ is isomorphic to $\Spin(Q)$ as algebraic groups.
\end{enumerate}
 \end{prop}

Note that in the above proposition, 
$\eta_1',\eta_2',\eta_3'$ correspond to 
$\eta_2,\eta_3,\eta_1$ respectively. 
For example, by the isomorphism $M\cong RT(\Oct)$
and $\eta_1$, we have a map 
\begin{equation*}
(t_1,t_2,t_3)
\to (\widehat t_3,t_1,t_2)
\to (t_3,\widehat t_2,\widehat t_1) 
\to (\widehat t_2,\widehat t_1,\widehat t_3), 
\end{equation*}
which is $\eta_3'$. 

From now on, we denote $M$ by $\Spin(Q)$.

\section{Stabilizer}
\label{sec:stabilizer}

In this section, we shall determine the structure of 
the stabilizer $G_w$ as an algebraic group.


%
\begin{thm}
\label{thm:stabilizer}
$G_w$ is isomorphic to $\GL(1) \times (\Spin(Q) \rtimes \Sss)$.
In particular, $G_w$ is a smooth reductive algebraic group.
\end{thm}
\begin{proof}
 To prove this proposition, we first determine
 the identity component $G_w^{\circ}$ of the stabilizer mainly
 by Lie algebra computations. Next, given an arbitrary element 
$g\in G_w$, 
 we replace $g$ by the action of $\Spin(Q) \rtimes \Sss$
so that it commutes with all elements of $\Spin(Q)$. 
Then Schur's lemma enables us to simplify the situation 
enough to prove the proposition. 

We first prove the following proposition. 
\begin{prop}
\label{prop:identity-comp}  
\begin{itemize}
\item[(1)] 
$G_w^{\circ}$ is isomorphic to  $\GL(1)\times \Spin(Q)$
over $k$.   
\item[(2)] 
$(G,V)$ is an irreducible prehomogeneous vector space. 
\end{itemize}
\end{prop}
\begin{proof}
We shall prove this proposition in the following two steps.

(I) \, We show that dim $G_w = 29$ and $G_w$ is smooth.
In this process, we prove (2).

(II) \, We construct an injective homomorphism 
$\xi : \GL(1)\times \Spin(Q) \to G_w$ which 
is defined over $k$ and induces an isomorphism between
$\text{Lie}(\GL(1)\times \Spin(Q))$ and $\text{Lie}(G_w)$.

If (I) and (II) are proved, since
the dimensions of both $\GL(1)\times \Spin(Q)$ and $G_w$ are equal 
to $29$ and both algebraic groups are smooth,  we have
$\GL(1)\times \Spin(Q) \cong G_w^{\circ}$ by $\xi$.

(I)
Let $\Lg_w$ and $\Lh_w$ 
denote the Lie algebras of
$G_w$ and $H_w$ respectively. 
We shall prove that $\Lh_w = \mathfrak{D}_0 $ and dim $H_w = 28$.
Since obviously $\mathfrak{D}_0 \subset \Lh_w$, 
it is enough to show that $\Lh_w \subset \mathfrak{D}_0.$

Let $Z = \left( X,  
\bigl( \begin{smallmatrix}
a&b\\ c&d
\end{smallmatrix} \bigr) \right)
\in \Lh_w \subset \Le \oplus \mathfrak{gl}(2)$, 
where $X \in \Le$ and $ 
\bigl( \begin{smallmatrix}
a&b\\ c&d
\end{smallmatrix} \bigr)
\in \Lgl(2)$.
Then by Proposition \ref{prop:Lie-alg-G1}, 
there exists
$X_0  \in \mathfrak{D}_{0}, \ \al,\ \be,\ \ga \in \Oct$ and  
$Y \in \J$ such that
\begin{equation*}
X = X_0 + (\al)'_1 + (\be)'_2 + (\ga)'_3 + R_Y
\end{equation*}
and  $\mathrm{Tr}(Y) = 0$.

Let $T =  (\al)'_1 + (\be)'_2 + (\ga)'_3 $.
Then
\begin{equation}\label{Tw_1}
\begin{split}
T&
\begin{pmatrix}
1&0&0\\
0&-1&0\\
0&0&0
\end{pmatrix}
= \frac{(\al)_1}{4} + \frac{(\be)_2}{4} - \frac{(\ga)_3}{2}, \\
&\\
T&
\begin{pmatrix}
0&0&0\\
0&1&0\\
0&0&-1
\end{pmatrix}
= -\frac{(\al)_1}{2} + \frac{(\be)_2}{4} + \frac{(\ga)_3}{4}. 
\end{split}
\end{equation}

Let $Y = 
\begin{pmatrix}
t_1 & y_3 & \oline{y_2}\\
\oline{y_3} & t_2 & y_1\\
y_2 & \oline{y_1} & t_3
\end{pmatrix},
$
where $t_i \in  k$, $y_i \in \Oct \ (i = 1, 2, 3)$ and $t_1 + t_2 +t_3 = 0$.
Then
\begin{align*}
\begin{pmatrix}
1 & 0 & 0 \\
0 & -1 & 0 \\
0 & 0 & 0
\end{pmatrix}
\circ Y 
& = \frac{1}{2}\left(
\begin{pmatrix}
t_1 & y_3 & \oline{y_2}\\
- \oline{y_3} & -t_2 & -y_1\\
0 & 0 & 0
\end{pmatrix}
+
\begin{pmatrix}
t_1 & -y_3 & 0 \\
\oline{y_3} & -t_2 & 0 \\
y_2 & -\oline{y_1} & 0
\end{pmatrix}\right)\\
&= 
\begin{pmatrix}
t_1 & 0 &  \frac{1}{2}\oline{y_2}\\
0 & -t_2 &  -\frac{1}{2}y_1\\
 \frac{1}{2} y_2 & -  \frac{1}{2} \oline{y_1}
 & 0
\end{pmatrix}, \\
\begin{pmatrix}
0&0&0\\
0&1&0\\
0&0&-1
\end{pmatrix}
\circ Y 
&= 
\begin{pmatrix}
0 & \frac{1}{2}y_3 & - \frac{1}{2}\oline{y_2} \\ 
\frac{1}{2} \oline{y_3} & t_2 & 0 \\
- \frac{1}{2} y_2 & 0 & -t_3
\end{pmatrix}.
\end{align*}

Since $Z w = \left(X, \minimat\right)w = 0$,
we have
\begin{equation}
\label{TwY-eq1}
T
\begin{pmatrix}
1&0&0\\
0&-1&0\\
0&0&0
\end{pmatrix}
+
\begin{pmatrix}
a&0&0\\
0&-a+b&0\\
0&0&-b
\end{pmatrix}
+
\begin{pmatrix}
1&0&0\\
0&-1&0\\
0&0&0
\end{pmatrix}
\circ Y = 0
\end{equation}
and
\begin{equation}
\label{TwY-eq2}
T
\begin{pmatrix}
0&0&0\\
0&1&0\\
0&0&-1
\end{pmatrix}
+
\begin{pmatrix}
c&0&0\\
0&-c + d&0\\
0&0&-d
\end{pmatrix}
+
\begin{pmatrix}
0&0&0\\
0&1&0\\
0&0&-1
\end{pmatrix}
\circ Y = 0.
\end{equation}

By \eqref{Tw_1}, 
\eqref{TwY-eq1} and \eqref{TwY-eq2}, we have
\begin{equation}
\label{TwY-cal1}
\begin{pmatrix}
a +t_1 &-\frac 12\ga		
&\frac 12\oline{y_2} + \frac 14\oline{\be} \\
-\frac 12\oline{\ga}	&-a +b -t_2				
&-\frac 12y_1 +\frac 14\al \\
\frac 12y_2 + \frac 14\be
&-\frac 12\oline{y_1} +\frac 14\oline{\al}&-b
\end{pmatrix}
= 0
\end{equation}
and
\begin{equation}\label{TwY-cal2}
\begin{pmatrix}
c&\frac 12y_3 +\frac 14\ga
&-\frac 12\oline{y_2} +\frac 14\oline{\be} \\
\frac 12\oline{y_3} +\frac 14\oline{\ga}&-c+d+t_2
&-\frac 12\al \\
-\frac 12y_2 +\frac 14\be
&-\frac 12\oline{\al}&-d - t_3
\end{pmatrix}
= 0.
\end{equation}
By \eqref{TwY-cal1}, we have
\begin{equation*}
 b = 0, \ t_1 = -a, \ t_2 = -a, \ \ga = 0,
 \ y_1 = \frac {\al}2, \ y_2 = - \frac {\be}2.
\end{equation*}
By \eqref{TwY-cal2}, we have
\begin{equation*}
c =0,\ t_2 = -d,\ t_3 = -d, \ \al = 0, \ y_3 = 0, \ y_2 = \frac {\be}2.
\end{equation*}
Hence, we have
\begin{equation*}
\begin{aligned}
&t_1 = t_2 = t_3 = -a = -d, \\
&y_1 = y_2 =y_3 = \al = \be = \ga = b = c = 0.
\end{aligned}
\end{equation*}
Since $\mathrm{Tr}(Y) = 0$, $3t_1 = 0 \ (\mathrm{ch}(k) \neq 3)$.
Therefore,
$Y = 0, \minimat = 0$ and $\al = \be = \ga = 0$.
Thus, $Z = X_0$, and so $\Lh_w  = \mathfrak{D}_0$. 
Then, we also have dim $H_w \le \dim \ \Lh_w = \dim \ \mathfrak{D}_0 = 28$.

Since $\dim H_w + \dim H (w) = \dim H$ ($H(w)$ is the orbit),   
we have
\begin{equation*}
\dim H_w = \dim H - \dim H (w) 
\ge \dim H - \dim V 
\ge 82 -54 = 28. 
\end{equation*}
It follows that $\dim H_w = 28$. 
Since dim $H = 82$ and $H_w = 28$, $\dim H(w) = 54$. Moreover,
it is a dense subset of $V$.

We show that dim $G_w = 29$.
Since  $\dim H_1 = \dim \Le =78$, by the exact sequence 
$\eqref{E6exact}$, we have $\dim G_1 = 79$.
Thus, since the dimension of the $G$-orbit of $w$ is also $54$ 
and dim $G = 83$, we have dim $G_w = 29$.
Since $G(w)$ is a constructible subset which is irreducible,
it is locally closed. 
Moreover, since dim $G(w) = 54 = \dim V$, $G(w)$ is open in $V$.
It follows that 
the condition of Definition \ref{PVS-def} (i)
is satisfied. We pointed out in Introduction that
 $V$ is an irreducible representation and that there exists
a relative invariant polynomial.  So
 we have (2) of the proposition.

We prove that $G_w$ is smooth. For that purpose, 
it is enough to show that $\dim \Lg_w = \dim G_w = 29$. 
Let $\Lg_1$ denote the Lie algebra of $G_1$.
We define $dc$ and $dc'$ by the Lie algebra homomorphisms induced by 
the characters $c : G_1 \to \GL(1)$ and $c' : G \to \GL(1) $
respectively (see (\ref{eq:c'-defn})).
Then,
\begin{equation*}
\Lge = \left\{ L \in  \Lgl(\J) \ {\vrule \
\begin{matrix}
{}^{\forall} X \in \J, \ \det \left((1+\varepsilon L)(X)\right)
& = c(1+ \varepsilon L)\det(X) \hfill \\ 
& = \left(1+\varepsilon dc(L)\right)\det(X) \hfill
\end{matrix}}
\right\}.
\end{equation*}

For $a \in k$,
\begin{equation*}
\det((1+a\varepsilon)I_{27}(X)) 
= (1+a\varepsilon)^3 \det(X)
= (1+3a\varepsilon)\det(X).
\end{equation*}
So $dc(a I_{27}) = 3a$.  Since ch$(k) \neq 3$, 
$dc : \Lge \to \Lgl(1)$ is surjective. Hence, the sequence
\begin{equation*}
0 \to \Le \to \Lge \to \Lgl (1) \to 0
\end{equation*}
is exact.
For $L = (L_1, L_2)  \in \Lg_w \subset \Lge \oplus \Lgl (2)$,
if $dc'(L) = 0$, then $dc(L_1) = 0$. 
It follows that $L_1 \in \Le$ and $L \in \Lh_w$.   
Since $dc'\left((aI_{27}, -aI_2)\right) = 3a \in \Lgl(1)$, 
the restriction of $dc'$ on $\Lg_w$ is also surjective.   
Hence, the following sequence 
\begin{equation*}
0 \to \Lh_w \to \Lg_w \to \Lgl (1) \to 0
\end{equation*}
is exact.
Therefore, $\dim \Lg_w = 29$ and $G_w$ is smooth.


(II)
We shall construct an injective homomorphism 
$\xi : \GL(1) \times $Spin$(Q) \to G_w$.

We identify $\J$ with 
$k^3 \oplus \Oct \oplus \Oct \oplus \Oct$ by the 
following isomorphism:
\begin{equation}
\label{identify-J}
\begin{array}{ccc}
\J &\cong & k^3 \oplus \Oct \oplus \Oct \oplus \Oct \\
\rotatebox{90}{$\in$} & &\rotatebox{90}{$\in$} \\
\begin{pmatrix}
s_1&x_3&\oline{x_2}\\
\oline{x_3}&s_2&x_1\\
x_2&\oline{x_1}&s_3
\end{pmatrix} &\leftrightarrow&
\left((s_1, s_2, s_3), x_1, x_2, x_3\right)
\end{array}.
\end{equation}
We choose 
$\{ (1,0,0), (0,1,0), (0,0,1) \}$ as the basis of $k^3$.
By the isomorphism $\J\cong k^3 \oplus \Oct \oplus \Oct \oplus \Oct$, 
we express elements of $GL(\J)$ in $4\times 4$ block form. 

We define a homomorphism 
$\xi_0 : \mathrm{GL}(1) \times \mathrm{Spin}(Q) \rightarrow \GL(\J)$ by 
\begin{equation}
\label{eq:ABC} 
\begin{array}{cc}
\xi_0 ((t, (A, B, C))) = 
\begin{pmatrix}
t I_3	&	0	&	0	&0\\
0		&tA 	&0		&0\\
0		&0		&tB		&0\\
0		&0		&0		&tC
\end{pmatrix}
\end{array}
\end{equation}
for
\begin{math}
(t,A,B,C)\in \mathrm{GL}(1) \times \mathrm{Spin}(Q).
\end{math}

Let $Z =(t, (A, B, C)) \in \GL (1) \times \mathrm{Spin}(Q)$.
Note that we are identifying $\mathrm{Spin}(Q)$ with $M$. 
We prove that $\xi_0 (Z) \in G_1$.
For any $X$ in the form (\ref{eq:X-defn}), 
we have
\begin{equation*}
\xi_0 (Z)X = 
\begin{pmatrix}
ts_1 & tC(x_3) & t\overline{B(x_2)} \\
t\overline{C(x_3)} & ts_2 & tA(x_1)\\
tB(x_2)	& t\overline{A(x_1)} & ts_3
\end{pmatrix}
\end{equation*}
and 
\begin{equation*}
\begin{split}
\mathrm{det}(\xi_0(Z)X)
& = t^3 s_1s_2s_3 
+\mathrm{tr}(tA(x_1) tB(x_2) tC(x_3)) \\
& \quad -ts_1 \| tA(x_1) \| -ts_2 \| tB(x_2) \| -ts_3 \|tC(x_3)\|.
\end{split}
\end{equation*}
Since $A, B, C \in \SO(Q)$, they preserve the inner product.
Thus, we have $\| tA(x_1) \| = t^2 \| x_1 \|, \ \| tB(x_2) \| 
= t^2 \| x_2 \|,\ \| tC(x_3) \| = t^2 \| x_3 \|$ and
\begin{equation*}
\begin{split}
\mathrm{tr}\left(tA(x_1) tB(x_2) tC(x_3)\right)
&= t^3
\mathrm{tr}\left(\overline{C(\overline{x_1 x_2})} C(x_3)\right)
=2t^3Q\left( \overline{x_1 x_2}, x_3 \right) \\
&=t^3 \mathrm{tr}(x_1 x_2 x_3).
\end{split}
\end{equation*}
Therefore, $\mathrm{det}(\xi_0 (Z)X) = t^3 \mathrm{det}(X)$,
and so $\xi_0 (Z) \in G_1$.

We define a homomorphism $\xi$ by
\begin{equation*}
\xi \ : \GL (1) \times \mathrm{Spin}(Q) 
\ni Z  = (t, (A, B, C)) \mapsto
(\xi_0 (Z), t^{-1}I_2)
\in G_1 \times \GL(2). 
\end{equation*}
By (\ref{eq:ABC}), $\text{Im}(\xi)\subset G_w$.
Since $\GL (1) \times \mathrm{Spin}(Q)$ is connected, we have
$\mathrm{Im} (\xi) \subset G_w^{\circ}$.
Moreover, $\xi$ is obviously injective and 
defined over $k$.
Thus, since dim $\GL (1) \times \Spin (Q) = 29$, 
we have dim $\xi(\GL (1) \times\Spin (Q)) = 29$.

The induced  homomorphism of 
Lie algebras from Lie$(\GL(1) \times \Spin (Q))$ to $\Lg_w$ is as follows:
\begin{equation*}
\begin{array}{ccccc}
&\! \! \! \! \! \!\Lgl (1) \oplus \Lso (Q)\oplus \Lso (Q)\oplus \Lso (Q)
&& \Lge \oplus \Lgl(2) &
\\

&\rotatebox{90}{$\subset$} &&\rotatebox{90}{$\subset$}&\\

d\xi  :  \! \! \! &\! \! \! \! \! \!  
\mathrm{Lie}(\GL (1) \times \mathrm{Spin}(Q))           
&\! \! \! \! \! \!  \longrightarrow &\Lg_w & \\
		 &\! \! \! \! \! \! \rotatebox{90}{$\in$} &                 		  
& \rotatebox{90}{$\in$} &\\[-4pt]
	         &\! \! \! \! \! \! (t,U_1,U_2,U_3)         &\! \! \! \! \! \!  \longmapsto     
& \left(tI_{27} +
\begin{pmatrix}
0&0&0&0\\
0&U_1&0&0\\
0&0&U_2&0\\
0&0&0&U_3
\end{pmatrix},  
\begin{pmatrix}
-t&0\\
0&-t
\end{pmatrix}\right)
&
\end{array}.
\end{equation*}
Note that $(\widehat U_3, U_1,  U_2)$ satisfies the equation 
$\eqref{Local-tri}$.
Since $d\xi$ is injective and 
dim Lie$(\GL(1) \times \Spin (Q)) =$ dim $\Lg_w = 29$,
  $d\xi$ is an isomorphism.
Thus,  $G^{\circ}_w \cong \GL (1) \times\mathrm{Spin}(Q)$ over $k$.
This completes the proof of Proposition \ref{prop:identity-comp}
\end{proof}

From now on, we identify $G_w^{\circ}$ 
with $\GL (1) \times\mathrm{Spin}(Q)$ by the above isomorphism.


We next determine the structure of $G_w / G_w^{\circ}$.
We  prove that $G_w / G_w^{\circ}$ $\cong$ $\mathfrak{S}_3 $.
Since $G_w$ and $G_w^{\circ}$ are defined over $k$, $G_w / G_w^{\circ}$  
and the natural homomorphism  $\pi: G_w \to G_w / G_w^{\circ}$ are 
defined over $k$ 
(see \cite[p.]98, 6.8 Theorem]{Borel}). 
We prove that there is a finite subgroup of $G_w$ which is isomorphic to $\Sss$ and is
mapped bijectively to  $G_w / G_w^{\circ}$ by $\pi$.

Let 
\begin{equation}
\label{eq:tau-defn} 
\begin{array}{rl}
\tau_1 &= \left(\left(
\begin{array}{@{\,}ccc|ccccc@{\,}}
0 & 1 & 0 & & & & \\
1 & 0 & 0 & \multicolumn{4}{c}{\raisebox{-5pt}[0pt][0pt]{\Huge $0$}}\\
0 & 0 & 1 & & & & \\

\hline
\vspace{-8pt}
& & & & & & \\
& & & &\hspace{-8pt} 0 & \iota & 0 \\
\multicolumn{2}{c}{\kern.8em {\raisebox{-5pt}[0pt][0pt]{\Huge $0$}}}
& & &\hspace{-8pt} \iota & 0 & 0 \\
& & & &\hspace{-8pt}0 & 0 & \iota \\
\end{array}
\right),
\begin{pmatrix}
-1 & 0\\
1 & 1 
\end{pmatrix}
\right),
\\
\tau _2 & = 
\left(\left(
\begin{array}{@{\,}ccc|ccccc@{\,}}
1 & 0 & 0 & & & & \\
0 & 0 & 1 & \multicolumn{4}{c}{\raisebox{-5pt}[0pt][0pt]{\Huge $0$}}\\
0 & 1 & 0 & & & & \\
\hline
\vspace{-8pt}
 & & & & & & \\
 & & & & \hspace{-8pt} \iota & 0 & 0 \\
\multicolumn{2}{c}{\kern.8em {\raisebox{-5pt}[0pt][0pt]{\Huge $0$}}}&& &\hspace{-8pt} 0 & 0 & \iota \\
&&&&\hspace{-8pt} 0 & \iota & 0 \\
\end{array}
\right),
\begin{pmatrix}
1 & 1\\
0 & -1
\end{pmatrix}
\right),
\vspace{8pt}
\\
e_G & = (I_{27}, I_2)
\end{array}
\end{equation}
(see (\ref{eq:iota-defn})).  Then, by easy computations, we have
 $\tau_1^2 = \tau_2^2 = e_G$, $(\tau_1 \tau_2)^2 = \tau_2 \tau_1$ 
and $\tau_1, \tau_2 \in G_w$.
Let $\langle \tau_1, \tau_2 \rangle$ denote the finite subgroup of $G_w$
generated by $\tau_1, \tau_2$.
It is easy to see that $|\langle \tau_1, \tau_2 \rangle| = 6$.

For $g \in G_w$, 
we define $I_g \in \mathrm{Aut}(G_w)$ by
\begin{equation*}
I_g: G_w\ni x \mapsto gxg^{-1}\in G_w.
\end{equation*}
Since $G_w^{\circ}$ is a normal subgroup of $G_w$, 
by restricting $I_g$ to $G_w^{\circ}$,
 $I_g$ induces an element of Aut$(G_w^{\circ})$. 
Moreover, since $\GL (1)$ is contained 
in the center of $G_w$, $I_g$ also induces  an element of
Aut$(\Spin (Q))$. 
So, we regard $I_g \in \mathrm{Aut}(\Spin (Q))$.
Then, the map 
\begin{equation*}
G_w\ni g \mapsto I_g\in \mathrm{Aut}(\Spin(Q)) 
\end{equation*}
is a homomorphism. We denote this homomorphism by $I$.

We define a subgroup $\mathrm{Inn}(\Spin (Q))$ of 
$\mathrm{Aut}(\mathrm{Spin}(Q))$ by
\begin{equation*}
\text{Inn(Spin$(Q))$} = \{ I_h \in \mathrm{Aut}(\mathrm{Spin}(Q)) 
\,|\, h \in \mathrm{Spin}(Q) \}.
\end{equation*}
Since $\text{Spin}(Q)$ is a group of type $\mathrm{D}_4$
and the automorphism group of the Dynkin diagram of
$\mathrm{D}_4$ is $\Sss$, by Proposition \ref{Prop-M}(ii) and
\cite[p.190, Proposition]{Borel}, 
the following sequence
\begin{equation}\label{AutSpin}
0 \to \text{Inn(Spin$(Q))$} \to \mathrm{Aut}(\mathrm{Spin}(Q))  
\to \mathfrak{S}_3 \to 0
\end{equation}
is exact.
We denote the above homomorphism from
$\mathrm{Aut}(\mathrm{Spin}(Q)) $ to $\mathfrak{S}_3$ by $f$.
The definition of $f$ is as follows.
For $Z = (Z_1, Z_2, Z_3) \in \Spin(Q)$, we define 
${\mathrm{pr}}_i(Z) = Z_i \ (i =1, 2, 3)$. 
For any $x \in \mathrm{Aut}(\Spin(Q))$ and $i \in \{1, 2, 3 \}$, 
${\mathrm{pr}}_i\circ x$ is an irreducible representation 
of $\Spin(Q)$ on $\Oct$ and equivalent to
some ${\mathrm{pr}}_j$ where $j \in \{1, 2, 3 \}$.
Then we define $f(x)(i) = j$.

We denote the composition of homomorphisms 
\begin{equation*}
G_w\overset{I}\to \mathrm{Aut}(\Spin(Q)) \overset{f}\to \Sss
\end{equation*}
by $\varphi$.
It is easy to see that the representation $Z =  (Z_1, Z_2, Z_3) \mapsto$  
$\widehat Z_i$ is equivalent to ${\mathrm{pr}}_i$ for $i = 1,2,3$. 
Since $I_{\tau_1}(Z) = (\widehat Z_2, \widehat Z_1, \widehat Z_3)$,
$I_{\tau_2}(Z) = (\widehat Z_1, \widehat Z_3, \widehat Z_2)$,
we have
\begin{equation}
\label{varphi}
\varphi(\tau_1) = (1 \ 2), \;
\varphi(\tau_2) = (2 \ 3).
\end{equation}
It follows that
$\varphi( \langle \tau_1, \tau_2 \rangle ) = \Sss$. 
Since $|\langle \tau_1, \tau_2 \rangle|=6$, 
$\langle \tau_1, \tau_2 \rangle\cong \Sss$ and 
$\varphi : G_w \to \Sss$ is surjective.

We define  another natural homomorphism from $G_w$ to $\Sss$ 
as follows.
For $x = x_1 v_1 + x_2 v_2 \in V$, let $F_x$ be the cubic form 
defined in (\ref{eq:Fx-defn}). 
For $x \in V$, 
we define 
%
\begin{equation}
\label{eq:zero-p1}
\Zero(x) = \{ q \in \Proj^1_k \,|\, F_x (q) = 0 \}.
\end{equation}
We call $\Zero(x)$ the {\it zero set} of $x$.
Note that this is well-defined since $F_x$ is homogeneous.
We express elements of $\Proj^1_k$ by row vectors
as $q=(q_1,q_2)$.

We shall define an action of $G_x$ on $\Zero(x)$.  
For any $g = (g_1, g_2) \in G_x$, 
\begin{equation*}
F_x(v) = F_{gx}(v)  
= c(g_1)F_{x}(vg_2).
\end{equation*}
It follows that if $q \in \Zero(x)$, then $qg_2$
also belongs to $\Zero(x)$.
Hence, we can define a right action of  $G_x$ on $\Zero(x)$ by 
\begin{equation*}
G_x \times \Zero(x)
\ni \left((g_1, g_2), q \right)
\mapsto qg_2 \in \Zero(x). 
\end{equation*}
Let $\Zero(x) = \{ q_{x,1}, q_{x,2}, q_{x,3}\}$. For $g = (g_1, g_2)\in G_x$, 
let $\eta_x(g) \in \Sss$ be the element such that
$q_{x,i}g_2  = q_{x,\eta_x(g)(i)}$ for $i = 1,2,3$. 

We define a map $\eta_x$  by
\begin{equation}
\label{eq:psix-defn}
\eta_x : G_x \ni g \mapsto \eta_x(g) \in \Sss.
\end{equation}
For any $q_{x,i} \in \Zero(x)$,
$g=(g_1, g_2)$ and $h = (h_1, h_2) \in G_x$,  
\begin{equation}
\label{eq:q-eta-cal}
q_{x,\eta(gh)(i)} = 
q_{x,i}g_2 h_2 = q_{x, \eta(g)(i)} h_2 = q_{x, \eta(h)(\eta(g)(i))}.
\end{equation}
Thus,
\begin{equation*}
\eta_x(gh)(i) = \eta_x(h)(\eta_x(g)(i)). 
\end{equation*}
Therefore, the map  $\eta_x :G_x \to \Sss$ is a homomorphism.
The definition of $\eta_x$ depends on the order of 
$\{q_{x,1},q_{x,2},q_{x,3}\}$ but its conjugacy class 
depends only on $x$.

We remind the reader that $w\in V_k$ is the element defined
in (\ref{eq:w-defn}). 
By direct computation, we have $F_w(v_1, v_2) = v_1 v_2 (v_1 - v_2)$.
It follows that
\begin{equation*}
\Zero(w) = \{ (0,1), (1,1), (1,0) \} \subset \Proj^k_1.
\end{equation*}
By choosing $q_{w,1}=(0,1),q_{w,2}=(1,1),q_{w,3}=(1,0)$, 
we obtain  a homomorphism  $\eta_w:G_w \to \Sss$.
We  shall show at the end of this section that
the following diagram is commutative.
\begin{equation}
\label{later on}
\xymatrix{
G_w \ar[d]_{I} \ar[r]^{\eta_w}&\Sss \\
\mathrm{Aut}(\Spin(Q)) 
\ar[ur]^{\!\!\!\!\!\!\!\!\!\!\!\!\!\!\circlearrowleft}_{f} \\
}
\end{equation}

Next, we show that the following sequence
\begin{equation}
\label{gw/gw0}
0 \to G_w^{\circ} \to G_w \overset{\varphi}\to \Sss \to 0
\end{equation}
is exact and split.  We may assume that $k$ is algebraically closed.
Since $\varphi$ is surjective, 
we only have to show that $G_w^{\circ} = \Ker(\varphi)$.
Since $G_w^{\circ}$ is connected,
$G_w^{\circ} \subset \Ker(\varphi)$.
So we shall show that $G_w^{\circ} \supset \Ker(\varphi)$.

Let $g =(g_1, g_2) \in \Ker (\varphi) \subset G_1 \times \GL (2)$. 
Since $g \in \Ker (\varphi)$, $f(I_g)=1$. 
Hence, by the exact sequence 
$\eqref{AutSpin}$, $I_g \in \mathrm{Inn}(\Spin(Q))$.
It follows that there exists $h \in \Spin(Q)$ such that $I_{gh}$ 
is trivial on $\Spin(Q)$ (since $k=\overline k$). 
Thus, we may assume that $I_g$ is trivial on $\Spin(Q)$.
 
Since $g$ commutes with all elements of $\Spin(Q)$,
$g_1$ preserves each irreducible non-equivalent 
representation of $\Spin(Q)$ in $\J$.
Thus, $g_1$ is in the block form: 
\begin{equation*}
g_1 = 
\begin{pmatrix}
A_{0}	&0		&0		&0	 	\\
0		&A_1	 &0		&0 		\\
0		&0		&A_2	&0		 \\
0		&0		&0		&A_3
\end{pmatrix},
\end{equation*}
where $A_0 \in \GL(3), A_1, A_2, A_3 \in \GL(\Oct)$.
Moreover, since $A_i Z_i = Z_i A_i $ for all $i \in \{1, 2, 3 \}$ 
and $Z =(Z_1, Z_2, Z_3) \in \Spin (Q)$, 
by Schur's lemma \cite[p.7]{Fulton-Harris} ($k=\overline k$),
there exists $c_i \in k^{\times}$ for $i = 1, 2, 3$ such that
$A_{i} 	= c_i I_8$. 
Thus, 
\begin{equation*}
g_1 = 
\begin{pmatrix}
A_{0}	&0		&0		&0	 	\\
0		&c_1 I_8 &0		&0 		\\
0		&0		&c_2 I_8	&0		 \\
0		&0		&0		&c_3 I_8
\end{pmatrix}.
\end{equation*}

Let $A_{0} = (a_{ij} )_{i,j =1, 2, 3}$, where $a_{ij} \in k$.
Since $g_1 \in G_1$, there exists $c \in k^{\times}$ such that 
$\det(g_1 X) = c\det(X) $ for any $X \in \J$.
Let
\begin{equation*}
X_1 =
\begin{pmatrix}
1		&t_3		&t_2		\\
t_3		&0			&t_1		\\
t_2		&t_1		&0
\end{pmatrix}, 
\end{equation*}
where $t_1, t_2, t_3 \in k$.
Then,
\begin{equation*}
\det(X_1) = 2t_1t_2t_3 - t_1^2.
\end{equation*}
We have 
\begin{equation*}
\begin{aligned}
g_1 X_1 & =
\begin{pmatrix}
a_{11} & c_3 t_3	& c_2 t_2		\\
c_3 t_3 & a_{21}	& c_1 t_1		\\
c_ 2 t_2 & c_1 t_1 & a_{31}	
\end{pmatrix},
\vspace{-2pt}\\
\det (g_1 X_1) & = a_{11}a_{21}a_{31} + 2c_1c_2c_3t_1t_2t_3
- \sum_{i=1}^3 a_{i1}c_i^2 t_i^2.
\end{aligned}
\end{equation*}
Since $\det(g_1 X_1)  =c\det(X_1)$ for any $t_1, t_2, t_3 \in k$ and
ch$(k) \neq 2$ , we have 
\begin{equation}
\label{c_k}
\begin{aligned}
c_1c_2c_3 &= c \neq 0,\; 
a_{11} c_1^2 = c, \;
a_{21} c_2^2 = a_{31} c_3^2 = 0\\
\end{aligned}.
\end{equation}
So $a_{21} = a_{31} = 0$ since $c_2, c_3 \neq 0$.

Let
\begin{equation*}
X_2 =
\begin{pmatrix}
0		&t_3		&t_2		\\
t_3		&1			&t_1		\\
t_2		&t_1		&0
\end{pmatrix},\; 
X_3 =
\begin{pmatrix}
0		&t_3		&t_2		\\
t_3		&0			&t_1		\\
t_2		&t_1		&1
\end{pmatrix}.
\end{equation*}
Similarly as above, 
we have
\begin{equation}
\label{A00}
a_{ij} = 0 \; (i \neq j), \; a_{11} c_1^2 = a_{22} c_2^2 = a_{33} c_3^2 = c.
\end{equation}

Let $g_2 = \bigl( 
\begin{smallmatrix}
p & q \\ r & s
\end{smallmatrix} \bigr)$
where $p, q, r, s \in k$. 
Since $g$ fixes $w$, we have
\begin{equation}
\begin{aligned}
\label{fixw1}
a_{11}p = 1, \;
-a_{22}p + a_{22}q = -1, \;
-a_{33}q = 0, \\
a_{11}r = 0, \;
-a_{22}r + a_{22}s = 1, \;
-a_{33}s = -1.
\end{aligned}
\end{equation}
By $\eqref{fixw1}$ and $\eqref{A00}$, we have
\begin{equation*}
\begin{aligned}
a_{11} = a_{22} = a_{33}, \;
p = s = a_{11}^{-1}, \;
q = r = 0
\end{aligned}		
\end{equation*}
and
\begin{equation}\label{ca}
c_1^2 = c_2^2 = c_3^2 = a_{11}^{-1} c.
\end{equation}
By  $\eqref{c_k}$ and $\eqref{ca}$, 
we have $c_i = \pm a_{11}$ for $i = 1, 2, 3$. 
Let $\varepsilon_i = c_i / a_{11} \in \{ \pm 1 \}$ for $i = 1, 2, 3$.
Since $c_1^2 = a_{11}^{-1} c = a_{11}^{-1} c_1 c_2 c_3$, 
we have $\varepsilon_1\varepsilon_2\varepsilon_3 = 1$.
Thus, 
\begin{equation*}
(\varepsilon_1, \varepsilon_2, \varepsilon_3) 
= (1, 1, 1), (1, -1, -1), (-1, 1, -1) \;\text{or}\; (-1, -1, 1). 
\end{equation*}
Therefore,
\begin{equation*}
g = \left(a_{11} \left(
\begin{array}{@{\,}c|ccc@{\,}}
I_3&\multicolumn{3}{c}{\raisebox{0pt}[0pt][0pt]{\Large $0$}}\\
\hline
\vspace{-10pt}
&&&\\
&\varepsilon_1 I_8&0&0\\
\text{\Large $0$} &0&\varepsilon_2 I_8&0\\
&0&0&\varepsilon_3 I_8\\
\end{array}
\right),
a_{11}^{-1} I_2
\right).
\end{equation*}

Obviously,
 $ (I_8, I_8, I_8), \ $  $(-I_8, -I_8, I_8) \ $, $(-I_8, I_8, -I_8) \ $ 
and $(I_8, -I_8, -I_8) \in  M (= \Spin(Q))$. 
Hence, we have $g \in G_w^{\circ}$.
Thus, $G_w^{\circ} = \Ker (\varphi)$ and 
the sequence $\eqref{gw/gw0}$ is exact. So 
$G_w / G_w^{\circ} \cong  \Sss$.
We denote the restriction of $\varphi$ on 
$\langle \tau_1, \tau_2 \rangle$ by 
$\varphi|_{\langle \tau_1, \tau_2 \rangle}$.
Since 
$\varphi|_{\langle \tau_1, \tau_2 \rangle}$
is an isomorphism, 
$(\varphi|_{\langle \tau_1, \tau_2 \rangle})^{-1}$
is a section of $\varphi$.
So the exact sequence $\eqref{gw/gw0}$ is split.

Therefore, the restriction of the natural homomorphism 
$\pi : G_w \to G_w /  G_w^{\circ}$ on 
$\langle \tau_1, \tau_2 \rangle$ is an isomorphism 
 defined over $k$.
This completes the proof of Theorem \ref{thm:stabilizer}. 
\end{proof}

Theorem \ref{thm:stabilizer} implies that the assumption on
the stabilizer in Theorem \ref{thm:galois-coh} is satisfied.

Finally,
we shall show that the diagram $\eqref{later on}$ is commutative.
Since 
\begin{equation*}
\begin{array}{rl}
q_{w,i}\tau_1&=
\begin{cases}
\vspace{3pt}
(0, 1)\!
\begin{pmatrix}
-1 & 0\\
1 & 1
\end{pmatrix}
= (1, 1) & i = 1\\
(1, 1)\!
\vspace{3pt}
\begin{pmatrix}
-1 & 0\\
1 & 1
\end{pmatrix}
= (0, 1) & i = 2
\\
(1, 0)\!
\begin{pmatrix}
-1 & 0\\
1 & 1
\end{pmatrix}
= (-1, 0) & i = 3
\end{cases},
\vspace{2pt} \\
q_{w,i}\tau_2&=
\begin{cases}
\vspace{2pt}
(0, 1)\!
\begin{pmatrix}
1&1\\
0&-1
\end{pmatrix}
= (0, -1) &i = 1
\\
\vspace{2pt}
(1, 1)\!
\begin{pmatrix}
1&1\\
0&-1
\end{pmatrix}
= (1, 0)&i = 2
\\
(1, 0)\!
\begin{pmatrix}
1&1\\
0&-1
\end{pmatrix}
= (1, 1)&i = 3
\end{cases}, 
\end{array}
\end{equation*}
by (\ref{varphi}), we have
\begin{equation}
\label{eq:eta-cal}
\eta_w(\tau_1) = (1 \ 2) = \varphi(\tau_1), \;
\eta_w(\tau_2) = (2 \ 3) = \varphi(\tau_2).
\end{equation}
It follows that 
$\eta_w(\tau) = \varphi(\tau)$ for all
$\tau \in \langle \tau_1, \tau_2 \rangle$. 
For any $g \in G_w$, there exists $\tau \in \langle \tau_1,
\tau_2 \rangle$ such that $\varphi(g) = \varphi(\tau)$. This implies that
$g \tau^{-1} \in G_w^{\circ}$. So the $\GL(2)$ part of $g \tau^{-1}$
is a scalar matrix. 
Thus, $\eta_w(g \tau^{-1})$ is the unit element of $\Sss$.
So we have $\eta_w(g) = \eta_w(\tau) = \varphi(\tau) = \varphi(g)$.
Therefore, the diagram $\eqref{later on}$ is commutative.


\section{Rational Orbits}
\label{sec:Rational Orbits}	

In this section, we determine the arithmetic
interpretation of the set of rational orbits.
We have determined the stabilizer in the
previous section. This and Theorem \ref{thm:galois-coh}
enable us to reduce the consideration to that of Galois
cohomology.  Since $\text{H}^1(k,\GL(1))=\{1\}$,
we have to interpret the set
$\text{H}^1(k,\Spin(Q)\rtimes \Sss)$. 
We first present an alternative description of
$\Spin(Q)\rtimes \Sss$.

We denote the cubic subalgebra $\sum_{i = 1}^{3} k E_i $  of $\J$ by 
$\at$.
We define an algebraic group $\Autt$ by 
\begin{equation*}
\Autt = \{L \in F_4 \,|\, L(\at) = \at \}
\end{equation*}
(see (\ref{eq:F4-defn}) for the definition of $F_4$)
and denote the Lie algebra of  $\Autt$ by $\Lautt$. Then
\begin{equation*}
\Lautt = \{L \in \Laut\,|\, L(\at) \subset \at \}.
\end{equation*}

\begin{lem}
\label{lem:at=D0}  
$\Lautt = \mathfrak{D}_0$. 
\end{lem}
\begin{proof}
Let $X = X_0 + (\al)'_1 + (\be)'_2 + (\ga)'_3 \in \Laut = \Der$ 
where $X_0 \in \mathfrak{D}_0,
\al, \be, \ga \in \Oct$.  
Since 
\begin{equation*}
X(E_i) =
\begin{cases}
\vspace{2pt}
\dispfr{1}{4}\left((\be)_2 
-(\ga)_3\right) &i = 1 \\
\vspace{2pt}
\dispfr{1}{4}\left((\ga)_3
-(\al)_1\right)&i = 2 \\
\vspace{2pt}
\dispfr{1}{4}\left((\al)_1
-(\be)_2\right)&i = 3
\end{cases},
\end{equation*}
if $X \in \Lautt$, then
we have $\al = \be = \ga = 0$.
Hence, we have
$\Lautt \subset \mathfrak{D}_0$.
The opposite inclusion is trivial and so 
$\Lautt = \mathfrak{D}_0$.
\end{proof}

Let ${\mathrm{pr}}_1$ be the natural projection from $G$ to $G_1$ as in
Introduction. 
Let $P={\mathrm{pr}}_1(\Sss)$. Then it is easy to see that ${\mathrm{pr}}_1$ induces 
an isomorphism $\Sss \cong P$. 
Since all elements in $\Spin(Q)$ preserve three  
octonions $\Oct_i$ where $i = 1,2,3$ (see (\ref{eq:oct-i-defn})) 
in $\J$ and such an element in $P$ 
is only $1$, $\Spin(Q)\cap P=\{1\}$. So 
${\mathrm{pr}}_1$ induces an isomorphism 
\begin{math}
\Spin(Q)\rtimes \Sss \cong \Spin(Q)\rtimes P. 
\end{math}
The following lemma is 
an alternative description of
$\Spin(Q)\rtimes \Sss$. 

\begin{lem}
\label{lem:Autt}
\begin{math}
\Spin (Q) \rtimes P = \Autt.
 \end{math}
\end{lem}
\begin{proof}
 Since $\Spin(Q)$ acts trivially on $\at$
(elements of $\Spin(Q)$ are in the form (\ref{eq:ABC}) where $t=1$)
and  all elements in $P$ preserve $\at$ (see (\ref{eq:tau-defn})), 
we have $(\Spin (Q) \rtimes P)  \subset
\Autt$. It follows that 
\begin{equation*}
28 = \dim \Spin(Q) \rtimes P \le \dim \Autt \le \dim \Lautt = 28.
\end{equation*}
In particular, $\dim \Autt = \dim \Lautt = 28$ and $\Autt$ is smooth.
We have shown in the step (I) of the proof of Proposition \ref{prop:Lie-alg-G1}
that $\mathfrak D_0=\mathfrak{so}(Q)$. So by Lemma \ref{lem:at=D0},
$\Lie(\Spin (Q) \rtimes P) = \Lautt$. Since both  
$\Spin(Q)\subset \Autt$ are smooth, the identity component of  
$\Autt$ coincides with $\Spin(Q)$.
Therefore, we may assume that $k=\overline k$
and prove that 
$\Autt_{\oline{k}} \subset (\Spin (Q) \rtimes P)_{\oline{k}}$
set-theoretically.  

By the above consideration, $\Spin(Q)$ is a normal subgroup of $\Autt$.
Hence, we can define a homomorphism from $\Autt$ to $\Aut(\Spin(Q))$ by
\begin{equation*}
\begin{array}{cccc}
I : 	&\Autt& \to &\Aut (\Spin(Q)) \\
	&\rotatebox{90}{$\in$}&&\rotatebox{90}{$\in$}\\
	&g&\mapsto&I_g : X \mapsto g Xg^{-1}  
\end{array}
\end{equation*}
for $X \in \Spin(Q)$.

Since $\Spin(Q) \rtimes P \subset \Autt$ and the restriction of $I$
to $\Spin(Q) \rtimes P$ is surjective,
$I$ is also surjective.
Hence, for $g \in \Autt$, by multiplying an element of 
$\Spin(Q) \rtimes P$, we may 
assume that $I_g$ is trivial. 
Then, by the same argument as in Section 
\ref{sec:stabilizer},
we have  
\begin{equation*}
g = \left(
\begin{array}{@{\,}ccc|ccccc@{\,}}
a_1&0&0&&&&\\
0&a_2&0&\multicolumn{4}{c}{\raisebox{-5pt}[0pt][0pt]{\Huge $0$}}\\
0&0&a_3&&&&\\
\hline
\vspace{-8pt}
&&&&&&\\
&&&&\hspace{-8pt}c_1 I_8&0&0\\
\multicolumn{2}{c}{\kern.1em {\raisebox{-5pt}[0pt][0pt]
{\quad \quad \quad \Huge $0$}}}&& &\hspace{-8pt}0&c_2 I_8&0\\
&&&&\hspace{-8pt}0&0&c_3 I_8\\\end{array}
\right),
\end{equation*}
where $a_i, c_i \in k^{\times}$ $(i =1, 2, 3)$ and
\begin{equation}
\label{eq:ai-ci}
a_1 c_1^2 = a_2 c_2^2 = a_3 c_3^2 = c_1 c_2 c_3
\end{equation}
(see $\eqref{c_k}, \eqref{A00}$).

Since $g \in \Autt \subset F_4 = \Aut(\J)$, we have $g(e) = e$. So
$a_1 = a_2 = a_3 = 1$. By (\ref{eq:ai-ci}),
we have $c_1^2 = c_2^2 = c_3^2 = c_1 c_2 c_3$.
It follows that
\begin{equation*}
(c_1, c_2, c_3) = (1,1,1), (1,-1,-1), (-1, 1,-1) \ \text{or} \ (-1,-1, 1).
\end{equation*}
Thus, $g \in \Spin(Q)$. 
Therefore, 
$\Autt_{\oline{k}} \subset (\Spin (Q) \rtimes P)_{\oline{k}}$.
\end{proof}

We shall consider the interpretation of $\mathrm{H}^1(k, G_w)$.  
We have shown in Section \ref{sec:stabilizer} that 
$G_w \cong \GL(1) \times (\Spin(Q) \rtimes \Sss)$.
Note that $\mathrm{H}^1(k, \GL(n))=\{1\}$ for all $n$. 
So $\mathrm{H}^1(k, G_w),\mathrm{H}^1(k,G)$
can be identified with 
$\mathrm{H}^1(k,\Spin(Q) \rtimes \Sss),\mathrm{H}^1(k,G_1)$
respectively.  So by Lemma \ref{lem:Autt}, 
we may consider
$\Ker \big(\mathrm{H}^1(k,\Autt) \to \HGE \big)$.

Let $\mathcal{M}$ be a Jordan algebra over $k$ and 
$\mathfrak{n}$ be a cubic \'etale subalgebra of $\mathcal{M}$ over $k$.
If there exists an isomorphism between $\Mj$ and $\J$ over $\ks$ 
which induces an isomorphism between $\an$ and $\at$ over $\ks$,
we  call $(\Mj, \an)$ a $k$-form of $(\J, \at)$.
Two $k$-forms $(\J_1, \at_1),(\J_2, \at_2)$
are said to be equivalent if
there exists an isomorphism between $\J_1$ and $\J_2$ over $k$ 
which induces an isomorphism between $\at_1$ and $\at_2$ over $k$.
Note that even if two cubic \' etale subalgebras $\at_1$ and $\at_2$ are 
isomorphic over $k$, two $k$-forms
$(\J, \at_1)$ and $(\J, \at_2)$ may not be equivalent.

We shall give a correspondence between
$\mathrm{H}^1(k,\Autt)$ and the set of $k$-forms
of $(\J, \at)$.

Suppose that $\Mj$ is a $k$-form of $\J$. 
For $\sigma \in \Gal$, let $v_{\Mj}(\sigma)$ be the semi-linear homomorphism
\begin{align*}
v_{\Mj}(\sigma) : \Mj_{\ks} =	\ks \otimes \Mj  \ni
\sum a_i\otimes m_i \mapsto \sum a_i^{\sigma}\otimes m_i 
\in \ks \otimes \Mj.
\end{align*}
We define a map 
\begin{equation}
 A: \mathrm{H}^1(k,\Autt)
  \to\{ \text{equivalence classes of k-forms of $(\J, \at)$} \}
\end{equation}
as follows.

Let $c \in \mathrm{H}^1(k,\Autt)$ be the element determined by
a 1-cocycle $h$ with coefficients in $\Autt_{\ks}$.
We define 
\begin{equation*}
\Mh = \{ X \in \J_{\ks} \,|\, h_1(\sigma)   v_{\J} (\sigma)(X) = X \},
\end{equation*}
\begin{equation*}
\nh = \{ Y \in \at_{\ks} \,|\, h_1(\sigma)   v_{\J} (\sigma)(Y) = Y \}.
\end{equation*}
Then, $\Mh$, $\nh$ are closed under
the addition and 
the multiplication of $\J_{\ks}$ 
since $h(\sigma) \in F_{4\,\ks}$ and 
the action of the Galois group commutes with
the addition and the multiplication of $\J_{\ks}$. 
So $(\Mh, \nh)$ is a $k$-form of $(\J, \at)$.
We define $A(h) = (\Mh, \nh)$.
One can show by the standard argument that
$A(h)$ depends only on the cohomology class $c$.  
By abuse of notation, we use the notation $A(c)$
also. 


Next, we define a map
\begin{equation*}
B: 
\{ \text{equivalence classes of k-forms of $(\J, \at)$} \} 
\to \mathrm{H}^1(k,\Autt)
\end{equation*}
as follows.

Let $(\Mj, \an)$ be a $k$-form of $(\J, \at)$.
Let $\phi : \Mj_{\ks} \to \J_{\ks} $ be an isomorphism of 
Jordan algebras which induces an isomorphism from $\an_{\ks}$
to $\at_{\ks} $.
Then we define a map $h_{\Mj}$ from $\Gal$ to $\GL(\J)_{\ks}$ by  
\begin{equation*}
h_{\Mj} (\sigma) 
= \phi \,  v_{\Mj}(\sigma)  \, \phi^{-1}   v_{\J}(\sigma)^{-1}.  
\end{equation*}
Then one can show by an easy computation,
that $h_{\Mj}$ is 1-cocycle with coefficients in
$\Autt_{\ks}$.  Let $B((\Mj,\an))$ be the element
of $\mathrm{H}^1(k,\Autt)$ determined by $h_{\Mj}$. 
This definition is also well-defined.

The proof of the following proposition is 
standard and is left to the reader. 

\begin{prop}
The maps $A,B$ are the inverses of each other and so 
\begin{equation*}
A : \HGw \to \{ \text{equivalence classes of k-forms of $(\J, \at)$} \}
\end{equation*}
is bijective.
\end{prop}
%
%

We next consider the interpretation of 
the set $\Ker \left(\HGw \to \HGE \right)$.

We first define the notion of {\it isotopes} of $\J$.
For $m \in \J$ such that $\det(m) \neq 0$, we can  
associate a Jordan algebra $\J_m$ related to $m$ as follows.
The underlying vector space of $\J_m$ is $\J$. 
We define
\begin{equation*}
\begin{array}{rcl}
\vspace{2pt}
\langle  x, y\rangle _m &=& -6\det(m)^{-1}D(x, y, m)  
+ 9\det(m)^{-2}D(x, m, m)D(y, m, m),  \\
\vspace{2pt}
x \circ_m y &=& 4\det(m)^{-1}(x \times m) \times (y \times m) +
\dispfr{1}{2}
(\langle  x, y\rangle _m - \langle  x, m\rangle _m \langle  y, m\rangle
_m)m 
\end{array}
\end{equation*}
for $x, y \in \J$. 
It is known that $\J_m$ is a Jordan algebra
with $x\circ_m y$, $\langle  x, x\rangle _m$ and $m$
as the product, the norm and the unit element. 
Moreover, there is an isomorphism of Jordan algebras
$\J_m\otimes \ks\to \J\otimes \ks$. 
This construction is given in \cite[p.155]{Springer}.  

Note that if the octonion is split, then all isotopes of
$\J$ are isomorphic to $\J$ over $k$
(it follows from the three facts [p. 155, Proposition 5.9.2], 
[p. 147, Proposition 5.6.2] and [p. 153, Proposition 5.8.2]).

For an arbitrary $k$-from $\Mj$ of $\J$,
the notion of the determinant on $\Mj$
is defined as follows. 
Let $\phi : \Mj_{\ks} \to \J_{\ks}$ be an isomorphism.
For any $X \in \Mj$, we define $\det_{\Mj}$ on $\Mj$
by $\det_{\Mj}(X) = \det(\phi(X))$.
Since automorphisms of $\J$ preserve  the determinant on $\J$
(see \cite[p.133, Proposition 5.3.10]{Springer}),
$\det_{\Mj}$ is independent of the choice of
the isomorphism $\phi$ and is defined over $k$.

It is also known that the determinant $\det_m$ on $\J_m$ satisfies 
the following equation:%
\begin{equation}
\label{eq:det-isotope}
\text{det}_m(x) = \det(m)^{-1}\det(x)
\end{equation}
for $x \in \J$ (see 
\cite[p.155, Proposition 5.9.2]{Springer}).  

Let $\text{JIC}(k)$ be the set of equivalence classes of pairs
$(\J_m, \an)$ of isotopes of $\J$ and their cubic \' etale subalgebras.
We show that elements in $\text{JIC}(k)$ correspond 
to elements in $\Ker \big(\mathrm{H}^1(k, G_w) \to  \mathrm{H}^1(k, G)\big)$
bijectively.

Let $h$ be a 1-cocycle with coefficients in
$G_{w\,\ks}$ which corresponds to $(\Mj,\an)$.
We write $h=(h_1,h_2)$ where
$h(\sigma)=(h_1(\sigma),h_2(\sigma))$
and $h_1(\sigma)\in G_{1\, \ks},h_2(\sigma)\in \GL(2)_{\ks}$.  

Let $c\in \mathrm{H}^1(k, G_w)$ be the cohomology 
class determined by $h$ and 
$A(c)=(\Mj,\an)$. 
Suppose that 
$c\in \Ker \big(\mathrm{H}^1(k, G_w) \to  \mathrm{H}^1(k, G)\big)$.
Since $h$ is trivial in $ \mathrm{H}^1(k, G)$, 
there exists $g = (g_1, g_2)\in G_{\ks}$ such that
\begin{equation*}
h(\sigma) = g^{-1} g^{\sigma} 
\end{equation*}
for any $\sigma \in \Gal$.

For $g = (g_1, g_2)\in G$, 
we define an element of $G_1$ by
\begin{equation}
\label{eq:mu-g-defn}
\mu_g = c(g_1) \det(g_2)^2 g_1 \in G_1.
\end{equation}
Note that the map $g \to \mu_g$ is a homomorphism.
If $g \in G_{w\,\ks}$, then by Theorem \ref{thm:stabilizer} and its proof,
there exist $t \in k, u\in \Spin(Q), \tau \in \langle \tau_1, \tau_2\rangle$
 such that
$g = (tu, t^{-1}I_2)\tau$.
So by computation,
\begin{equation}
\mu_g = c(tu\mathrm{pr}_1(\tau))(t^{-4})tu\mathrm{pr}_1(\tau) = u\mathrm{pr}_1(\tau)
\end{equation}
where $\mathrm{pr}_1$ is the natural projection from $G$ to $G_1$.
It follows that $\mu_g \in (\Spin(Q) \rtimes P)_{\ks} = \Autt_{\ks}$.

For the above 1-cocycle $h$,
since $a = \left( c(g_1) \det(g_2)^2 I_{27}, 
(c(g_1)\det(g_2)^2)^{-1}I_2\right)\in G_{w\,\ks}$,
we may substitute  $a^{-1}h(\sigma)a^{\sigma}$ for $h(\sigma)$
for any $\sigma \in \Gal$.
Then $h_1(\sigma) = \mu_g^{-1}\mu_g^{\sigma} = \mu_{g^{-1}g^{\sigma}}$.

We show that there exists $m \in \J$ with $\det(m)\neq 0$ such that 
$\Mj$ is isomorphic to $\J_m$.
Let $\mu_g(e) = m$.
Since $g^{-1}g^{\sigma}\in G_{w\,\ks}$,  
$\mu_{g^{-1}g^{\sigma}} \in  \Autt_{\ks}\subset F_{4\, \ks}$. 
Therefore, we have $\mu_g^{-1}\mu_g^{\sigma}(e) = e$. 
So 
\begin{align*}
v_{\J}(\sigma)(m) &  = v_{\J}(\sigma)(\mu_g(e)) 
= v_{\J}(\sigma)(\mu_g(v_{\J}(\sigma)^{-1}(v_{\J}(\sigma)(e)))) \\
& = \mu_g^{\sigma}(v_{\J}(\sigma)(e)) =
\mu_g^{\sigma}(e) = \mu_g(e) = m 
\end{align*}
for any $\sigma \in \Gal$.
Hence, we have $m \in \J$.

It is known that all isotopes are isomorphic over $\ks$ (it follows from 
\cite[p.158, Proposition 5.9.3]{Springer} and the proof of 
\cite[p.181, Theorem 7.3.2]{Springer}).
So there exists an isomorphism $\phi :\J_{\ks} \to \J_{m\, \ks}$.
Since $\phi^{-1} \mu_g(e) = e$, by (\ref{eq:F4-e}), 
$\phi^{-1}  \mu_g \in F_{4\,\ks}$.
Thus, $\mu_g \, (= \phi  \phi^{-1}  \mu_g)$ is an isomorphism
from $\J_{\ks}$ to $\J_{m\,\ks}$. Since 
\begin{equation}
\label{eq:Mj-defn}
\begin{array}{rcl}
\Mj &=& \{ X \in \J_{\ks} \,|\, \mu_g^{\sigma}(v_{\J}(\sigma)(X)) 
= \mu_g(X), 
\quad {}^{\forall} \sigma \in \Gal \} \\
&=& \{ X \in \J_{\ks} \,|\, v_{\J}(\sigma)(\mu_g(X)) = \mu_g(X), 
\quad {}^{\forall} \sigma \in \Gal \} \\
&=& \mu_g^{-1}(\J)
\end{array},
\end{equation}
the $k$-linear map $\mu_g : \Mj \to \J_m$ is a $k$-linear isomorphism.
Therefore, $\Mj \cong \J_m$ as Jordan algebras over $k$. 

Conversely, suppose that $m\in \J$, $\det(m)\not=0$
and $\an$ is a cubic \' etale subalgebra of $\J_m$. 
We show that $B((\J_m, \an)) \in \KerHGw$.

If $k = \ks$, then all cubic \'etale algebras are isomorphic to 
$(\ks)^3$.
So there exist three  idempotents $u_1, u_2, u_3$ such that
$u_i$'s are pairwise orthogonal and
$u_1 + u_2 + u_3$ is the unit element.
So for any cubic \'etale subalgebras $\an_1, \an_2$ of $\J$, 
there exist $u_{i1}, u_{i2}, u_{i3} \in \an_{i\,\ks}$ for $i = 1, 2$
which satisfy the above conditions.
Since $\Oct$ is split if $k = \ks$, 
the map $\Oct^{\times}_{\ks} 
\ni x \mapsto
\|x\|  \in (\ks)^{\times}$ is surjective.
In the situation where $k = \ks$, by  
\cite[p.413, T{\tiny HEOREM} 9]{Albert},
there exists $\phi \in F_{4\,\ks}$ such that 
$\phi(u_{1j}) = u_{2j}$ for $j = 1,2,3$.
Then $\phi(\an_{1\,\ks}) = \an_{2\,\ks}$.

For any $(\J_m, \an)$, there exists an isomorphism
$\phi_1: \J_{m\,\ks} \to \J_{\ks}$.
By the above argument, there exists $\phi_2 \in F_{4\ks}$ such that
$\phi_2(\phi_1(\an_{\ks})) = \at_{\ks}$.
Then $\phi_2\phi_1 : \J_{m\, \ks} \to \J_{\ks}$ is an isomorphism 
which induces an isomorphism 
from $\an_{\ks}$ to $\at_{\ks}$. 
So $(\J_m, \an)$ is a $k$-form of $(\J, \at)$.
Moreover, the map 
$\sigma \mapsto \phi \, v_{\J_{m}}(\sigma)\, 
\phi^{-1} v_{\J}(\sigma)^{-1}$ 
is a 1-cocycle which defines the cohomology class
$B((\J_m, \an))$. 

Since the 
determinant on $\J_{m\, \ks} $ satisfies $\eqref{eq:det-isotope}$ and 
$\det(\phi(X)) = \det_m(X)$ for any $X \in \J_{m\, \ks}$
, we have $\phi \in G_{1 \,\ks}$. 
It follows that $ B((\J_m, \an)) \in 
\Ker$ $\big(\mathrm{H}^1(k, G_w)$ $\to$ $\mathrm{H}^1(k, G)\big)$. 
Thus, 
elements in $\text{JIC}(k)$ 
correspond to elements in $\Ker \big(\HGw \to$ $\HGE \big)$ bijectively.

Therefore, by the above argument and Theorem \ref{thm:galois-coh},
we have the following interpretation of the set of rational orbits.
\begin{thm}
The following map 
\label{thm:main-theorem}
\begin{equation*}
G_k \backslash V_k^{\mathrm{ss}} \ni x \mapsto	A(c_x) \in \mathrm{JIC}(k)
\end{equation*}
(see (\ref{thm:galois-coh}) for the definition of $c_x$)
is bijective.
\end{thm}
%


\section{Equivariant map}
\label{sec:equivariant-map}

In Section \ref{sec:Rational Orbits}, we have shown
that elements of
the set  $\KerHGw$ are in bijective correspondence with
equivalence classes of 
pairs $(\J_m,\an)$ of isotopes of $\J$ and
their cubic subalgebras.
In this section, for any $x \in V_k^{\mathrm{ss}}$, 
we construct explicitly the 
isotope of $\J$ and its cubic subalgebra corresponding to $x$
in Theorem \ref{thm:main-theorem} 
by an equivariant map from $V$ to $\J$.

Before defining an equivariant map, 
we review some properties of the cross product.
The following formulas are proved in 
\cite[p.154]{Springer} and 
\cite[pp.122-123, Lemma 5.1.2]{Springer}.

\begin{lem}
\label{cross-pro}
The following equations are satisfied for all $x,y\in \J$.

\begin{enumerate}

\item $\displaystyle
\vspace{3pt}
x \times y = x \circ y - \frac{1}{2}\langle  y,e\rangle x
- \frac{1}{2}\langle  x,e\rangle y - 
\frac{1}{2}\langle  x,y\rangle e +
\frac{1}{2}\langle  x,e\rangle \langle  y,e\rangle e \displaystyle  \displaystyle$.
\item $x\circ (x \times x) = \det(x)e$ (see (\ref{cross-formula}) also).  
\vspace{2pt}
\item $(x \times x) \times (x \times x) = \det(x)x \displaystyle$. 
\vspace{2pt}
\item $4x \times \left(y \times(x \times x) \right) = \det(x)y
+ \langle  x, y\rangle x \times x$. 
\vspace{2pt}
\item $4(x \times x) \times (x \times y) = 
\det(x)y + 3D(x, x, y) x$.
\vspace{2pt}
\item $4(x \times y) \times (x \times y) = 
-2(x \times x) \times (y \times y) 
+ 3D(x, y, y)x + 3D(x, x, y)y 
\displaystyle $.
 \end{enumerate}
 \end{lem}

For any $g \in \GL(\J)$, we define $\widetilde{g} \in \GL(\J)$ by
\begin{equation*}
\langle  g(x), \widetilde{g}(y)\rangle = \langle x, y\rangle,
\quad {}^{\forall}  x, y \in \J.
\end{equation*}
The following lemma is proved in 
\cite[p.180, Proposition 7.3.1] {Springer}.
\begin{lem}
\label{lem:tilde}
The map $g \mapsto \widetilde{g}$ is an automorphism of $H_1$ 
with order $2$ and the following equations:
\begin{align*}
g(x \times y) = \widetilde{g}(x) \times \widetilde{g}(y), \ \,
\widetilde{g}(x \times y)  = g(x) \times g(y) \quad
{}^{\forall}  x, y \in \J
\end{align*}
are satisfied.
\end{lem}



We define a map $m:V=\J \otimes \Aff^2\to \J$ by
\begin{equation*}
m(x_1v_1 + x_2v_2)
= 6(x_1 \times x_1) \times (x_2 \times x_2) - 3D(x_1,x_2,x_2)x_1 
-3D(x_1,x_1,x_2)x_2.
\end{equation*}
We define an action of $G$ on $\J$ by 
\begin{equation}
G \times \J \ni ((g_1, g_2), x) \mapsto c(g_1)(\det(g_2))^2g_1(x) \in \J.
\end{equation}
\begin{prop}
\label{prop:G-eq}
For any $(g_1, g_2) \in G$ and $x_1v_1 + x_2v_2 \in V$,
\begin{equation*}
m\left((g_1, g_2)(x_1v_1 +  x_2v_2)\right) 
= c(g_1)\det(g_2)^2g_1(m\left(x_1v_1 +  x_2v_2)\right).
\end{equation*}
\end{prop}
\begin{proof}
We first prove that the map $m$ is $\GL(2)$-equivariant.
Since $\GL(2)$ is generated by upper triangular matrices and 
$(
\begin{smallmatrix} 
0&1\\ 1&0
\end{smallmatrix}
)$, it is enough to consider these two cases.
We denote $m(x_1v_1 + x_2v_2 )$ by $m(x_1, x_2)$. 
For any upper triangular matrix $g_2 = \uppmat$,
\begin{equation*}
\renewcommand{\arraystretch}{1.15}
\begin{array}{rl}
\vspace{2pt}
\!m(ax_1 + bx_2, dx_2) =&\!\!\!
 6\left((ax_1 + bx_2) \times (ax_1 + bx_2)\right) 
\times (dx_2 \times dx_2) \\
\vspace{2pt}
&\!\!\!\!-3D(ax_1 + bx_2, dx_2, dx_2)(ax_1 + bx_2)  \\
&\!\!\!\!-3D(ax_1+bx_2, ax_1 +bx_2, dx_2)dx_2 \\
\vspace{2pt}
=&\!\!\! a^2d^2 \{6(x_1\times x_1)\times (x_2\times x_2) \\ 
&\!\!\!\! -3D(x_1,x_2,x_2)x_1 -3D(x_1,x_1,x_2)x_2\} \\
\vspace{2pt}
&\!\!\!\! + 12abd^2(x_1 \times x_2) \times (x_2 \times x_2) 
+6b^2d^2(x_2 \times x_2) \times (x_2 \times x_2)  \\
&\!\!\!\! -3abd^2D(x_2,x_2,x_2)x_1 -9abd^2D(x_1,x_2,x_2)x_2 \\
&\!\!\!\!  -6b^2d^2D(x_2,x_2,x_2)x_2 \\
\vspace{2pt}
=&\!\!\! a^2d^2m(x_1, x_2)  \\
&\!\!\!\!+ 3abd^2\{4(x_1 \times x_2) \!\times\! (x_2 \times x_2) \!
-\det(x_2)x_1\! -3D(x_1, x_2, x_2)x_2\}\\
\vspace{2pt}
=&\!\!\! a^2d^2m(x_1,x_2) \; \text{(by (v))} \vspace{2pt}\\
=&\!\!\! (\det(g_2))^2m(x_1, x_2)
\end{array}.
\end{equation*}
For $g_2 = 
(
\begin{smallmatrix} 
0&1\\ 1&0
\end{smallmatrix}
)$, by definition,
$m(g_2(x_1v_1 + x_2v_2)) = m(x_2,x_1) = m(x_1, x_2)$.
Therefore, $m$ is $\GL(2)$-equivariant.

We next prove that $m$ is $G_1$-equivariant.
For any $g \in H_1$, by Lemma \ref{lem:tilde},
\begin{align*}
g\left((x \times y) \times (z \times w)\right)
& = \left(\widetilde g(x\times y)\right) \times
\left(\widetilde g(z \times w)\right) \\
& = \left(g(x) \times g(y)\right) \times \left(g(z) \times g(w)\right)
\end{align*}
for  $x, y, z, w \in \J$.
For any $g_1 \in G_{1\, \ks}$, since ch$(k) \neq 3$, there exists $t \in \ks$ 
such that
$t^3 = c(g_1)$.
Then $t^{-1} g_1  \in H_{1\, \ks}$. So 
\begin{align*}
g_1\left((x \times y)\times (z \times w)\right)  
&=
t(t^{-1}g_1)\left((x \times y)\times (z \times w)\right)\\
&=
t\left(t^{-1}g_1(x) \times t^{-1}g_1(y)\right) 
\times \left(t^{-1}g_1(z) \times t^{-1}g_1(w)\right)\\
&=
 c(g_1)^{-1}
\left(g_1(x) \times g_1(y)\right) \times \left(g_1(z) \times g_1(w)\right)
\end{align*}
for $x, y, z, w \in \J$.
It follows that 
\begin{equation*}
\renewcommand{\arraystretch}{1.1}
\begin{array}{rl}
\vspace{2pt}
g_1\left(m(x_1, x_2)\right) 
=& 6g_1((x_1 \times x_1) \times (x_2 \times x_2)) \\
&\ - 3D(x_1,x_2,x_2)g_1(x_1) -3D(x_1,x_1,x_2)g_1(x_2)\\
=& c(g)^{-1}(m\left(g_1(x_1),  g_1(x_2))\right).
\end{array}
\end{equation*}
Therefore, for any $(g_1, g_2) \in G$, we have 
\begin{equation}
\label{G-equiv-m}
m\left((g_1, g_2)(x_1v_1 +  x_2v_2)\right) =
 c(g_1)\det(g_2)^2g_1(m\left(x_1v_1 +  x_2v_2)\right).
\end{equation}
\end{proof}

For $x\in \Vss_k$, 
to define the isotope $\J_{m(x)}$ of $\J$,  $\det(m(x))$ must not be
equal to $0$.
So we need the following lemma. 
\begin{lem}
\label{lem:m(w)-compu}
For $w = w_1v_1 + w_2v_2 \in V$ in (\ref{eq:w-defn}), 
$
m(w) = e.
$
In particular, $\J_{m(w)} = \J$.
\end{lem}
\begin{proof}
By direct computation, we have
\begin{align*}
&w_1\times w_1 = - E_3, \;\, w_2\times w_2 = - E_1, \;\,
(w_1 \times w_1) \times (w_2\times w_2) = \frac 12 E_2, \\[-2pt]
&D(w_1,w_1,w_2) = \frac 13, \;\,
D(w_1,w_2,w_2) = -\frac 13,
\end{align*}
(see  (\ref{elements-defn}) for the definition of $E_1, E_2, E_3$).
So 
\begin{equation*}
\begin{array}{rl}
m(w) &= 6(w_1 \times w_1) \times (w_2\times w_2) - 
3D(w_1,w_2,w_2)w_1 - 3D(w_1,w_1,w_2)w_2\\
&=  3E_2 +w_1 - w_2\\
&= e.
\end{array}
\end{equation*}

It is known that the multiplication in $\J_e$ is the same as 
that of $\J$
(see \cite[p.158]{Springer}).
So $\J_{m(w)} = \J$.
\end{proof}
By the above lemma and
Proposition \ref{prop:G-eq}, 
the map $x \mapsto \det(m(x))$
is a non-trivial relative invariant polynomial of degree 12.
So $x \in  V^{\mathrm{ss}}$ if and only if $\det(m(x)) \neq 0$.
Thus, for all $x \in V^{\mathrm{ss}}_k$,
we can define the isotope $\J_{m(x)}$ of $\J$.

We next define a subspace $\at(x)$ of $\J$ 
which is a cubic \' etale $k$-subalgebra  
with respect to the multiplication 
in $\J_{m(x)}$ for $x = x_1v_1 + x_2v_2 \in V^{\mathrm{ss}}_k$. 
We shall show that the pair $(\J_{m(x)}, \at(x))$ is a $k$-form
of $(\J, \at)$ which is equivalent to $A(c_x)$ 
in Theorem \ref{thm:main-theorem}.
 
Let $\at(x)$ be the vector space spanned by $x_1, x_2$ and $m(x)$.
Since $m(w) = e$, it is easy to see that $\at(w) = \at$.
For any $g_1 \in G_1$, by Proposition \ref{prop:G-eq}, 
\begin{equation*}
\at((g_1, I_2)(x)) = g_1\left(\at(x)\right).
\end{equation*}
Moreover, it is easy to see that $\at(x)$ is $\GL(2)$-invariant.
Hence, for any $g \in G$, we have 
\begin{equation*}
\at(g(x)) = g_1\left(\at(x)\right).
\end{equation*}

We show that  $\at(x)$ is a cubic \' etale subalgebra 
of $\J_{m(x)}$.
For $g = (g_1,g_2) \in G_{\ks}$, 
let $\mu_g = c(g_1) \det(g_2)^2 g_1 \in G_{1\,\ks}$
(see (\ref{eq:mu-g-defn}) also).
For any $x = x_1v_1 + x_2v_2 \in V^{\mathrm{ss}}_k$,
there exists $g = (g_1, g_2) \in G_{\ks}$ such that 
$gw = x$.
Then, 
\begin{equation}
\label{eq:ax}
m(x) = m(gw) = c(g_1)\det(g_2)^2 g_1 (m(w)) = \mu_g(e).
\end{equation}
%
By the same argument as the one just before $\eqref{eq:Mj-defn}$,
the equation $\eqref{eq:ax}$ implies that 
the linear map 
$\mu_g : \J_{\ks} \to \J_{m(x)\, \ks}$ is an isomorphism.

Moreover, since 
\begin{equation*}
\at(x)_{\ks} = \at(gw)_{\ks} = g_1(\at(w)_{\ks}) = g_1(\at_{\ks})
= \mu_g(\at_{\ks}),
\end{equation*}
$\at(x)_{\ks}$ is a cubic \' etale subalgebra of $\J_{m(x)\, \ks}$.
Since the multiplication in $\J_{m(x)\, \ks}$ is defined over $k$, 
$\at(x)$ is a cubic \' etale subalgebra of $\J_{m(x)}$.
Moreover, $(\J_{m(x)}, \at(x))$ is a $k$-form of $(\J, \at)$. 

By the argument above, 
 the pair $( \J_{m(x)}, \at(x))$ corresponds 
to the cohomology class  $c_x \in \HGw$
 (see $\eqref{eq:cx-defn}$ for the definition of $c_x$). 
Therefore, we have the following theorem. 
\begin{thm}
The following map is bijective.
\label{thm:main-theorem2}
\begin{equation*}
G_k \backslash V_k^{\mathrm{ss}} \ni x \mapsto	(\J_{m(x)}, \at(x)) 
\in \mathrm{JIC}(k).
\end{equation*}
\end{thm}

Since the product structure on $\J_m$ is defined by homogeneous
polynomials  of degree $11$, the product structure on $\J_{m(x)}$
is defined by homogeneous polynomials  of degree $44$.  However,
it is possible to construct an equivariant map from $V$ to
${\mathrm{Hom}}_k(\J\otimes \J,\J)$ defined by homogeneous polynomials
of degree $8$. For this, see \cite{kato-yukie-jordan-equiv}.


\section{The split case}
\label{sec:split-case}

The purpose of this section is to prove 
Theorem \ref{thm:main-split-intro}. 
For that purpose we need some preparation. 

When we consider actions of ${\mathrm{Gal}}(\ks/k)$ 
on $\mathfrak S_n$ ($n$ is a positive integer) 
in this section, we always consider the 
trivial action of the Galois group. 
So if $\{h(\sigma)\}$ is a 1-cocycle with coefficients in 
$\mathfrak S_n$, then
${\rm{Gal}}(k^{\rm{sep}}/k)\ni 
\sigma\mapsto h(\sigma)^{-1}\in\mathfrak S_n$
is a homomorphism by Definition \ref{def:galcoh}. 
We identify $\text{H}^1(k,\mathfrak S_n)$ 
with conjugacy classes
of homomorphisms from 
${\rm{Gal}}(k^{\rm{sep}}/k)$ to $\mathfrak S_n$
in this manner. 

We use the space of pairs of ternary quadratic forms
in order to prove Theorem \ref{thm:main-split-intro}. 
The interpretation of rational orbits 
in this case was carried out in \cite{wryu}. 
We first review this case and point out 
properties which will be needed in 
this section.

Let $\cV_1$ be the space of ternary quadratic forms. 
Since $\ch (k)\not=2$, we may identify 
$\cV_1$ space of symmetric $3\times 3$ 
matrices.  Let $\cG_1=\GL(3)$. If $g\in \cG_1$ and 
$x\in \cV_1$, then we define $\rho_1(g) x = gx{}^t g$.
This defines a representation of $\cG_1$ on $\cV_1$. 
Let $\cG=\GL(3)\times \GL(2)$. 
Regarding ${\rm{Aff}}^2$ as the standard representation
of $\GL(2)$, 
$\cV=\cV_1\otimes {\rm{Aff}}^2$ is a representation of
$\cG$. We use the notation $\rho(g)$ for this representation. 
The representations $\rho_1,\rho$ can be 
extended to representations on $\J,V$ respectively
by the same formulas.  Note that elements of $k$
commute with elements of $\Oct$ and so
the order of multiplication does not matter.  
We later show (Corollary \ref{cor:enchou})
that $\rho_1(g_1)\in G_1$ (resp. $\rho(g)\in G$)
for $g_1\in \GL(3)$ (resp. $g\in \cG$).

We express elements of $\cV$ as 
$x=(x_1,x_2)$ ($x_1,x_2\in\cV_1$)
or $x=x_1v_1+x_2v_2$ using variables
$v=(v_1,v_2)$.  
If $x\in \cV_1$, the determinant 
defined by (\ref{eq:Jdet-defn}) 
coincides with the usual determinant. 
It is well-known that $(\cG,\cV)$ is a regular prehomogeneous vector 
space and if we express elements of 
$\cV$ as $x=x_1v_1+x_2v_2$, then the discriminant of 
the cubic form $F_x(v)$ (see (\ref{eq:Fx-defn}))  
is a relative invariant polynomial. 
Therefore, $x\in \cV$ is semi-stable 
with respect to the action of $\cG$
if and only if 
it is semi-stable as an element of $V$ 
with respect to the action of $G$.

If $x=(x_1,x_2)\in\cV_k$, then 
$x_1,x_2$ define two conics in 
$\Proj^2_{\overline k}$. 
Let $\Zero_{\Proj^2}(x)\subset \Proj^2_{\overline k}$ be 
the intersection of these two conics. 
This should not be confused with 
$\Zero(x)\subset \Proj^1_{\overline k}$ 
(see (\ref{eq:zero-p1})).  
If $g=(g_1,g_2)\in \cG_{\overline k}$, 
$x,y\in \cV^{\rm{ss}}_{\overline k}$
and $x=gy$, then it is easy to see that 
$\Zero_{\Proj^2}(x)=\Zero_{\Proj^2}(y)g_1^{-1}$.
In particular, if $g=(g_1,g_2)\in \cG_{x\, \overline k}$, 
then $\Zero_{\Proj^2}(x)=\Zero_{\Proj^2}(x)g_1^{-1}$ and so 
$\Zero_{\Proj^2}(x)=\Zero_{\Proj^2}(x)g_1$ also. 

The following proposition is proved in 
\cite[pp.290, Propositions 1.5, 1.6]{wryu}. 

\begin{prop}
\label{prop:quartic-semi-stable}
Suppose that $x\in \cV_k$. Then 
$x\in \cV^{\rm{ss}}_k$ if and only if $\Zero_{\Proj^2}(x)$ 
consists of four distinct points. 
Moreover, if $x\in \cV^{\rm{ss}}_k$ then 
the coordinates of points in $\Zero_{\Proj^2}(x)$ 
belong to $k^{\rm{sep}}$. 
\end{prop}
Note that the second statement is now a consequence of 
\cite{regularity-of-PVS}.

If $x\in \cV^{\rm{ss}}_k$, 
$\Zero_{\Proj^2}(x)=\{p_{x,1},\ldots,p_{x,4}\}$
and $g\in \cG_{x\,\overline k}$,  
then let $\gamma_x(g)\in \mathfrak S_4$ be the 
element such that 
\begin{equation}
\label{eq:pix-defn}
p_{x,j}g = p_{x,\gamma_x(g)(j)} 
\end{equation}
for $j=1,\ldots,4$.  Then 
$\gamma_x:\cG_{x\,\overline k}\ni g\mapsto \gamma_x(g)\in \mathfrak S_4$
is a homomorphism. 

%

The element $w$ in Section 1 of \cite{wryu} is not $w$ in this paper, 
but is the following element:
\begin{align*}
w'=\frac 12\left(
\begin{pmatrix}
0 & 0 & -1 \\
0 & 0 & 1 \\
-1 & 1 & 0 
\end{pmatrix},
\begin{pmatrix}
0 & 1 & 0 \\
1 & 0 & -1 \\
0 & -1 & 0
\end{pmatrix}
\right).
\end{align*}
Let 
\begin{equation*}
g = \left(
\frac 12
\begin{pmatrix}
1 & -1 & -1 \\
-1 & 1 & -1 \\
-1 & -1 & 1  
\end{pmatrix},
\begin{pmatrix}
1 & 0 \\
-1 & -1 
\end{pmatrix}
\right)\in \cG_k. 
\end{equation*}
Then $gw=w'$.

It is easy to see that  
$\Zero_{\Proj^2}(w)=\{p_{w,1},\ldots,p_{w,4}\}$
where 
\begin{equation*}
p_{w,1}=(1,-1,-1),p_{w,2}=(-1,1,-1),p_{w,3}=(-1,-1,1),
p_{w,4}=(1,1,1). 
\end{equation*}
We denote $\tau_1,\tau_2,\tau_3$ in \cite[p.289]{wryu}
by $\overline{\tau_1},\overline{\tau_2},\overline{\tau_3}$
in this paper. It is easy to see that 
$\overline{\tau_1},\overline{\tau_2},\overline{\tau_3}$
generate a subgroup of $\cG_{w'\, k}$ which 
is isomorphic to $\mathfrak S_4$ 
and maps isomorphically to 
$\cG_{w'\, k}/\cG_{w'\, k}^{\circ}$.
Moreover, $\gamma_{w'}$ induces an isomorphism
$\cG_{w'\, k}/\cG_{w'\, k}^{\circ}\cong \mathfrak S_4$.

By computation, 
\begin{equation}
\label{eq:grhotau}
\begin{aligned}
g^{-1} \overline{\tau_1} g & = 
\left(
\begin{pmatrix}
0 & 1 & 0 \\
1 & 0 & 0 \\
0 & 0 & 1 
\end{pmatrix},
\begin{pmatrix}
-1 & 0 \\
1 & 1
\end{pmatrix}\right), \\
g^{-1} \overline{\tau_2} g & = 
\left(
\begin{pmatrix}
1 & 0 & 0 \\
0 & 0 & 1 \\
0 & 1 & 0
\end{pmatrix},
\begin{pmatrix}
1 & 1 \\
0 & -1
\end{pmatrix}\right), \\
g^{-1} \overline{\tau_3} g & = 
\left(
\begin{pmatrix}
0 & -1 & 0 \\
-1 & 0 & 0 \\
0 & 0 & 1 
\end{pmatrix},
\begin{pmatrix}
-1 & 0 \\
1 & 1
\end{pmatrix}\right).
\end{aligned}
\end{equation}
Therefore,  
\begin{equation}
\label{eq:rhotau}
\rho(g^{-1} \overline{\tau_1} g)  
= \tau_1,\;
\rho(g^{-1} \overline{\tau_2} g)  
= \tau_2,\;
\rho(g^{-1} \overline{\tau_3} g)  
= \tau_1(-I_8, - I_8, I_8), 
\end{equation}
where $\tau_1, \tau_2 $ are the elements in $\eqref{eq:tau-defn}$ 
and $(I_8, - I_8, -I_8) \in \Spin(Q)$ (with the understanding
that elements of $\Spin(Q)$ act on diagonal entries of elements of
$\J$ trivially).
Note that $\tau_1, \tau_2, 
\tau_1(-I_8, - I_8, I_8)$ 
correspond to $(1, 2), (2,3), (1,2) \in \Sss$ 
by the homomorphism $\eta_w$ from $G_w$ to $\Sss$
(see $\eqref{eq:psix-defn}$ and $\eqref{eq:eta-cal}$).
These elements generate a subgroup of  
$G_{w\, k}$ which is isomorphic to $\mathfrak S_4$.
We identify $\mathfrak S_4$ with this subgroup.

For $x\in \cV^{\rm{ss}}_k$, let
$\Zero_{\Proj^2}(x)=\{p_{x,1},\ldots,p_{x,4}\}$. 
We define a map $\phi_x:\Gal \to \mathfrak S_4$
by 
\begin{equation}
\label{eq:p-sigma-defn}
p_{x,j}^{\sigma} = p_{x,\phi_x(\sigma)(j)}
\end{equation}
for $\sigma\in {\rm{Gal}}(k^{\rm{sep}}/k)$ and  
$j=1,\ldots,4$. It is easy to verify that 
$\phi_x$ is a homomorphism.

Suppose that $x=g_x w\in \cV_k^{\rm{ss}}$
where $g_x=(g_{x,1},g_{x,2})\in \cG_{k^{\rm{sep}}}$. 
We define another map 
$\psi_x:{\rm{Gal}}(k^{\rm{sep}}/k)\to \mathfrak S_4$
by 
\begin{equation}
\label{eq:psix-sigma-defn}
\psi_x(\sigma) = \gamma_w(g_x^{-1}g_x^{\sigma})^{-1} 
\end{equation}
for $\sigma\in {\rm{Gal}}(\ks/k)$. 
Since the map
${\rm{Gal}}(\ks/k)\ni \sigma\mapsto
\gamma_w(g_x^{-1}g_x^{\sigma})\in\mathfrak S_4$
is an anti-homomorphism, $\psi_x$ is a homomorphism.

By changing the numbering of 
$\Zero_{\Proj^2}(x)$, we may assume that 
$p_{x,j}=p_{w,j} g_{x,1}^{-1}$ ($j=1,\ldots,4$). 
Then 
\begin{equation*}
p_{x,j}^{\sigma}
= p_{w,j}(g_{x,1}^{\sigma})^{-1}
= p_{w,j}(g_{x,1}^{\sigma})^{-1}g_{x,1} g_{x,1}^{-1}
= p_{w,\psi_w(\sigma)(j)} g_{x,1}^{-1}
= p_{x,\psi_w(\sigma)(j)}.  
\end{equation*}
Therefore, $\phi_x$ and $\psi_x$ are conjugate. 
Hence, if we identify $\text{H}^1(k,\mathfrak S_4)$
with conjugacy classes of homomorphims from 
${\rm{Gal}}(k^{\rm{sep}}/k)$ to $\mathfrak S_4$, 
then the action of ${\rm{Gal}}(k^{\rm{sep}}/k)$ 
on $\{1,2,3,4\}$ can be identified with 
the action of ${\rm{Gal}}(k^{\rm{sep}}/k)$ 
on $\Zero_{\Proj^2}(x)$ up to conjugation.

The following theorem is proved in 
\cite[p.310, Theorem 5.3]{wryu}
(see \cite[p.120, Lemma (1.8)]{yukie-field-extension3} also). 
\begin{thm}
\label{thm:WY-quartic}
The set 
$\cG_k\backslash \cV^{\rm{ss}}_k$ 
is in bijective correspondence with
$\text{H}^1(k,\mathfrak S_4)$. 
If $x\in \cV^{\rm{ss}}_k$, then 
$\phi_x$ is the corresponding 
homomorphism from ${\rm{Gal}}(k^{\rm{sep}}/k)$
to $\mathfrak S_4$. 
\end{thm}

We next review known results regarding 
rational orbits of the space of binary cubic forms.

Let $\cW$ be the space of cubic forms
in two variables $v=(v_1,v_2)$ and 
$\cH = \GL(1)\times \GL(2)$. Since we are assuming
that ${\mathrm{ch}}(k)\not=3$, 
$\cW$ is an irreducible representation of 
$\cH$ where $\alpha \in \GL(1)$ acts by 
multiplication by $\alpha$. 
It is well-known and in fact very easy to prove 
that $(\cH,\cW)$ is a regular prehomogeneous
vector space, $w'' = F_w(v) = v_1v_2(v_1-v_2)\in \cW^{\rm{ss}}_k$,
$\cH_{w''}/\cH_{w''}^{\circ}\cong \mathfrak S_3$
and the set $\cH_k\backslash \cW_k^{\rm{ss}}$
is in bijective correspondence with 
$\text{H}^1(k,\mathfrak S_3)$.  Moreover, 
$x\in \cW$ is semi-stable if and only if
$x$ has three distinct factors.

If $x=x(v)\in \cW_k^{\rm{ss}}$, then we define 
\begin{equation}
\label{eq:zero-p1sym3}
\Zero_{\Proj^1}(x) = \{q\in \Proj^1_{\overline k}\ | \ x(q)=0\}
\end{equation}
and call it the zero set of $x$. 
Obviously, if $x\in V^{\mathrm{ss}}_k$, 
then $\Zero(x)=\Zero_{\Proj^1}(F_x)$. 
If $g=(g_1, g_2) \in \cH_{\oline k}$,
$x,y\in \cW^{\rm{ss}}_{\overline k}$ and $x = gy$, 
then it is easy to see that
$\Zero_{\Proj^1}(x) = \Zero_{\Proj^1}(y) g_2^{-1}$.
In particular, 
if $g_xw'' = x$ where $g_x = (g_{x, 1}, g_{x,2}) \in \cH_{\oline k}$,
then
$\mathrm{Zero}_{\Proj^1}(x) = \mathrm{Zero}_{\Proj^1}(w'')g_{x,2}^{-1}$.
If $\alpha\in\overline k$ and $(\alpha,1)\in\Zero_{\Proj^1}(x)$, 
then we may call $\alpha$ a root of $x$ also.

Let $\Zero_{\Proj^1}(x) =\{q_{x,1}, q_{x,2}, q_{x,3}\}$.
For $g = (g_1, g_2) \in \cH$, 
let $\theta_x(g) \in \Sss$ be the element such that
\begin{equation*}
q_{x,i}g_2 = q_{x, \theta_x(g)(i)}
\end{equation*}
for $i = 1,2,3$.
We define a map $\theta_x$ by
\begin{equation*}
\theta_x : \cH \ni g \mapsto \theta_x(g) \in \Sss.
\end{equation*}
By the same argument as in $\eqref{eq:q-eta-cal}$,
$\theta_x$ is a homomorphism.

If $x = g_x w''\in \cW_k^{\rm ss}$ ($g_x\in \cH_{w''\, \ks}$),
then ${\mathrm{Gal}}(\ks/k)\ni\sigma 
\mapsto \theta_{w''}(g_x^{-1}g_x^{\sigma})$
is the 1-cocycle which defines the element of 
${\mathrm H}^1(k,\Sss)$ which corresponds to $x$. 
Let $\lambda_x(\sigma) = \theta_{w''}(g_x^{-1}g_x^{\sigma})^{-1}$
for $\sigma\in {\mathrm{Gal}}(\ks/k)$. Then $\lambda_x$ 
is a homomorphism from ${\mathrm{Gal}}(\ks/k)$ to $\Sss$. ]

We define another map $\kappa_x : \Gal \to \mathfrak S_3$ by 
\begin{equation}
\label{eq:kappa-defn}
q_{x,i}^{\sigma} = q_{x, \kappa_x(\sigma)(i)}
\end{equation}
for $\sigma \in \Gal$ and $i = 1,2,3$.
Then $\kappa_x$ is also a homomorphism 
and as in the case of $\cV$, $\kappa_x$ 
and $\lambda_x$ are conjugate. 
So when we consider conjugacy classes, 
$\kappa_x$ and $\lambda_x$ can be identified.

We now start the proof of Theorem \ref{thm:main-split-intro}. 
%


\begin{proof}
[Proof of Theorem \ref{thm:main-split-intro}]
Let $h^1_{\eta_w}:{\mathrm{H}}^1(k,G_w)\to {\mathrm{H}}^1(k,\mathfrak S_3)$ 
be the map induced by $\eta_w$. This of course induces a map 
from the subset $\Ker({\mathrm{H}}^1(k,G_w)\to {\mathrm{H}}^1(k,G))$
of ${\mathrm{H}}^1(k,G_w)$ to ${\mathrm{H}}^1(k,\mathfrak S_3)$. 
We first prove that if this map 
from $\Ker({\mathrm{H}}^1(k,G_w)\to {\mathrm{H}}^1(k,G))$
to ${\mathrm{H}}^1(k,\mathfrak S_3)$ is a bijection,
then the element of ${\mathrm{H}}^1(k,\mathfrak S_3)$ which
corresponds to the orbit of $x\in V^{\mathrm{ss}}_k$
is the one obtained by the action of ${\mathrm{Gal}}(\ks/k)$
on the set $\Zero(x)$. Note that 
the set $\Ker({\mathrm{H}}^1(k,G_w)\to {\mathrm{H}}^1(k,G))$
is in bijective correspondence with $G_k\backslash V^{\mathrm{ss}}_k$.
 
We consider the homomorphism $G \to \cH$ defined as follows: 
\begin{equation*}
\Phi:G \ni (g_1, g_2) \mapsto (c(g_1), g_2) \in \cH.
\end{equation*}
It is easy to see that for any $x \in V$ and $g\in G$, 
$F_{gx} = \Phi(g)(F_x)$. Since $F_w=w''$,
$\Phi(G_w)\subset \cH_{w''}$.  
By choosing $q_{w'',1}=(0,1),q_{w'',2}=(1,1),q_{w'',3}=(1,0)$, 
the definitions of $\eta_w,\theta_{w''}$ imply that 
the following diagram 
\begin{equation}
\label{etaw-commutative}
\xymatrix{
G_w \ar[d]_{\Phi} \ar[r]^{\eta_w}&\Sss \\
\cH_{w''}
\ar[ur]^{\!\!\!\!\!\!\!\!\!\!\!\!\!\!\circlearrowleft}_{\theta_{w''}} \\
}
\end{equation}
is commutative.

Let $x=g_x w\in V^{\mathrm{ss}}_k$ where 
$g_x\in G_{\ks}$. Then 
$F_x = \Phi(g_x) w''$. Let $c_x$ be 
the element of ${\mathrm H}^1(k,G_w)$ 
as in (\ref{eq:cx-defn}). 
Then $F_x$ corresponds to the cohomology 
class in ${\mathrm H}^1(k,\cH_{w''})$  
defined by the 1-cocycle 
\begin{equation*}
{\mathrm{Gal}}(\ks/k)\ni\sigma \mapsto
\Phi(g_x)^{-1}\Phi(g_x)^{\sigma}
= \Phi(g_x^{-1}g_x^{\sigma})\in \cH_{w''\,\ks}.   
\end{equation*}
We define 
$h^1_{\Phi}:{\mathrm{H}}^1(k,G_w)\to {\mathrm{H}}^1(k,\cH_{w''})$, 
$h^1_{\theta_{w''}}:{\mathrm{H}}^1(k,\cH_{w''})\to {\mathrm{H}}^1(k,\mathfrak S_3)$
similarly as above. 
Then by the commutativity of the above diagram, 
\begin{equation*}
h^1_{\eta_w}(c_x)^{-1} =  h^1_{\theta_{w''}}(h^1_{\Phi}(c_x))^{-1} 
= \lambda_{F_x} = \kappa_{F_x}.  
\end{equation*}
Since the correspondence $x\mapsto h^1_{\eta_w}(c_x)^{-1}$ 
is bijective, so is the 
correspondence $x\mapsto \kappa_{F_x}$.  
So the orbits in $G_k\backslash V^{\mathrm{ss}}_k$
are determined by the action of the Galois 
group on $\Zero(x)$.

To prove the above bijectivity, 
we first consider the case where $k$ is a finite field. 
In this case, the well-known theorem of Lang \cite{Lang-Galois} 
says that the first Galois cohomology set of any smooth
connected algebraic group is trivial. So
${\mathrm{H}}^1(k,G)=\{1\}$. Therefore,
$G_k\backslash  V^{\rm{ss}}_k$ is in bijective correspondence with
${\mathrm{H}}^1(k,G_w)$. The map $\phi$ in the exact sequence (\ref{gw/gw0})
can be replaced by $\eta_w$ by the commutativity 
of the diagram (\ref{later on}). 
By \cite[p.120, Lemma (1.8)]{yukie-field-extension3}, 
$h^1_{\eta_w}$ is injective. Since the sequence (\ref{gw/gw0}) is split,
it is surjective. Hence, the map $h^1_{\eta_w}$
is bijective. 
This completes the proof of Theorem \ref{thm:main-split-intro}
for the case where $k$ is a finite field.


Therefore, we assume for the rest of this
section that 
$k$ is an infinite field and 
$\mathbb O$ is the split octonion $\widetilde {\mathbb O}$. 
The reason why we can prove the theorem 
when the octonion is split is that the
norm map $\mathbb O^{\times}\to k^{\times}$ is surjective.
Note that $\mathbb O^{\times}$
contains $\GL(2)_k$ as a subset, the restriction of 
the norm of $\mathbb O^{\times}$ to $\GL(2)_k$ 
coincides with the determinant and 
the determinant induces a surjective homomorphism 
$\GL(2)_k \to k^{\times}$ ($\GL(2)_k$ is a group).



What we need for later consideration about
the space $\cV$ is the following. 
First, the stabilizer of $4$ in $\mathfrak S_4$ is
obviously isomorphic to $\mathfrak S_3$.
Since this obvious subgroup plays an important role
in the subsequent proof, we name this subgroup as follows.  
\begin{mydef}
\label{defn:cS3}  
$\cS = \mathfrak S_3\subset \mathfrak S_4$.   
\end{mydef}

As we pointed out after (\ref{eq:rhotau}) that
$\cS$ is generated by $g^{-1}\overline{\tau_1}g$,  
$g^{-1}\overline{\tau_2}g$ of (\ref{eq:grhotau}).
As elements of $\mathfrak S_4$, these elements
are the transpositions $(1\, 2),(2\, 3)$.
The equation $\eqref{eq:eta-cal}$ implies that 
$(1\,2),(2\,3)\in\cS$ map to  
$(1\,2),(2\,3)\in\mathfrak S_3$ by  
$\eta_w$. Therefore, $\eta_w$
induces an isomorphism  
$\cS\cong \mathfrak S_3$.  

We shall prove the following proposition.
\begin{prop}
\label{prop:crucial-step}
Any orbit in $G_k\backslash \Vss_k$
has a representative in $\cV^{{\mathrm{ss}}}_k$. 
Moreover, such a representative
can be chosen so that it corresponds to a cohomology class in
$\text{H}^1(k,\mathfrak S_4)$ which comes from 
$\text{H}^1(k,\cS)$. 
\end{prop}
 
We show that Theorem \ref{thm:main-split-intro}
now follows from Proposition \ref{prop:crucial-step}. 
Suppose that Proposition \ref{prop:crucial-step} holds.

The natural maps $\cS\to\mathfrak S_4\to \cG_w\to G_w\to \cH_{w''}\to \Sss$
induce maps: 
\begin{equation*}
{\mathrm{H}}^1(k,\cS)
\to {\mathrm{H}}^1(k,\mathfrak S_4) 
\to {\mathrm{H}}^1(k,\cG_w) 
\to {\mathrm{H}}^1(k,G_w) 
\to {\mathrm{H}}^1(k,\cH_{w''}) 
\to {\mathrm{H}}^1(k,\mathfrak S_3).   
\end{equation*}
Note that the homomorphism $G_w\to \Sss$ 
can be regarded as $\eta_w$. 
Since the map
\begin{math}
{\mathrm{H}}^1(k,\cG_w)\to
{\mathrm{H}}^1(k,\cG_w/\cG_w^{\circ})
\cong {\mathrm{H}}^1(k,\mathfrak S_4)  
\end{math}
is bijective ($\cong$ means bijection)
and $\mathfrak S_4\to \cG_w$ is a section,
the map ${\mathrm{H}}^1(k,\mathfrak S_4) \to {\mathrm{H}}^1(k,\cG_w)$
is bijective. 
The map
\begin{math}
{\mathrm{H}}^1(k,\cG_w) 
\to {\mathrm{H}}^1(k,G_w)   
\end{math}
induces a map  
\begin{equation*}
{\mathrm{H}}^1(k,\cG_w)\cong \Ker({\mathrm{H}}^1(k,\cG_w)\to  {\mathrm{H}}^1(k,\cG))
\to \Ker({\mathrm{H}}^1(k,G_w)\to {\mathrm{H}}^1(k,G)).   
\end{equation*}

By assumption, the map
\begin{equation}
\label{eq:cStoGw}
{\mathrm{H}}^1(k,\cS) \to \Ker({\mathrm{H}}^1(k,G_w)\to {\mathrm{H}}^1(k,G))
\cong G_k\backslash V^{\mathrm{ss}}_k   
\end{equation}
is surjective. Since the composition 
\begin{equation*}
{\mathrm{H}}^1(k,\cS) \to {\mathrm{H}}^1(k,G_w)\to
{\mathrm{H}}^1(k,\mathfrak S_3) 
\end{equation*}
is bijective, (\ref{eq:cStoGw}) is injective and so bijective. 
Since ${\mathrm{H}}^1(k,\cS) \to {\mathrm{H}}^1(k,\mathfrak S_3)$
is bijective, 
\begin{equation*}
\Ker({\mathrm{H}}^1(k,G_w)\to {\mathrm{H}}^1(k,G))
\to {\mathrm{H}}^1(k,\mathfrak S_3)
\end{equation*}
is bijective and this map is induced by $\eta_w$.

We now prove Proposition \ref{prop:crucial-step}.
\begin{proof}
[Proof of Proposition \ref{prop:crucial-step}]
If $x\in \cV^{\mathrm{ss}}_k$ and 
$\mathrm{Zero}_{\Proj^2}(x)$ has a $k$-rational point, 
then the element of $\text{H}^1(k,\mathfrak S_4)$ 
which corresponds to $x$ comes from $\text{H}^1(k,\cS)$.
So what we have to prove is that any orbit in 
$G_k\backslash V^{\mathrm{ss}}_k$ has a representative
$x\in \cV^{\mathrm{ss}}_k$ such that 
$\mathrm{Zero}_{\Proj^2}(x)$ has a $k$-rational point.

We first describe some elements of $G_1$. 
This will be needed when we
show later that any orbit in $G_k\backslash \Vss_k$
has a representative in $\cV^{{\mathrm{ss}}}_k$.

For $\alpha\in\Oct^{\times}$, 
let 
$d_1(\alpha),d_2(\alpha),d_3(\alpha)$
be the $k$-linear maps $\J\to \J$ 
defined by 
\begin{equation}
\label{eq:di-defn}
\begin{aligned}
& d_1(\alpha):\J \ni 
\begin{pmatrix} 
s_1&x_3&\oline{x_2} \\
\oline{x_3}&s_2&x_1 \\
x_2&\oline{x_1}&s_3
\end{pmatrix} 
\mapsto 
\begin{pmatrix} 
\|\alpha\|s_1& x_3\oline{\alpha}& \oline{x_2}\alpha \\
\alpha\oline{x_3}&s_2&\alpha x_1\alpha /\|\alpha\|\\
\oline{\alpha} x_2&\oline{\alpha x_1\alpha}/\|\alpha\|&s_3
\end{pmatrix}\in \J, \\
& d_2(\alpha):\J \ni  
\begin{pmatrix} 
s_1 & x_3 & \oline{x_2} \\
\oline{x_3} & s_2 & x_1 \\
x_2 & \oline{x_1} & s_3
\end{pmatrix} 
\mapsto 
\begin{pmatrix} 
s_1 & \oline{\alpha}x_3 & \oline{\alpha x_2\alpha}/\|\alpha\|  \\
\oline{x_3}\alpha & \|\alpha\|s_2 & x_1\oline{\alpha}\\
\alpha x_2\alpha/\|\alpha\| & \alpha\oline{x_1} & s_3
\end{pmatrix}\in \J, \\
& d_3(\alpha):\J \ni 
\begin{pmatrix} 
s_1&x_3&\oline{x_2} \\
\oline{x_3}&s_2&x_1 \\
x_2&\oline{x_1}&s_3
\end{pmatrix} 
\mapsto 
\begin{pmatrix} 
s_1&\alpha x_3\alpha /\|\alpha\|& \alpha\oline{x_2} \\
\oline{\alpha x_3\alpha}/\|\alpha\|&s_2&\oline{\alpha}x_1\\
 x_2\oline{\alpha}&\oline{x_1}\alpha&\|\alpha\|s_3
\end{pmatrix}\in \J.
\end{aligned}
\end{equation}
Also for $\alpha_1,\alpha_2,\alpha_3\in k^{\times}$, 
let
\begin{equation}
\label{eq:dalpha-defn}
d(\alpha_1,\alpha_2,\alpha_3)
= \rho_1({\mathrm{diag}}(\alpha_1,\alpha_2,\alpha_3))
\in \GL(\J). 
\end{equation}

For $u\in \mathbb O$, let $n_{21}(u):\J\to \J$ be 
the $k$-linear map defined by 
\begin{equation}
\label{eq:nu-defn}
\J\ni X\mapsto 
\begin{pmatrix}
1 && \\
u & 1 & \\
&& 1  
\end{pmatrix}
X 
\begin{pmatrix}
1 & \overline u & \\
 & 1 & \\
&& 1  
\end{pmatrix}\in \J.
\end{equation}
We define $n_{ij}(u)$ for $i\not=j$  similarly.

\begin{lem}
\label{lem:niju-well-defined}
\begin{itemize}
\item[(1)] 
For all 
$\alpha_1,\alpha_2,\alpha_3\in k^{\times}$ 
and $u\in \mathbb O$, 
the right hand sides of (\ref{eq:dalpha-defn}), (\ref{eq:nu-defn}) 
do not depend on the order of multiplication and 
the resulting elements belong to $\J$. 
\item[(2)] 
For all $\alpha\in\Oct^{\times}$, 
the maps $d_i(\alpha)$ ($i=1,2,3$) 
belong to $G_{1\, k}=GE_{6\, k}$. 
Moreover, if $\|\alpha\|=1$, then 
the maps $d_i(\alpha)$ ($i=1,2,3$) 
belong to $H_{1\, k}=E_{6\, k}$.
\item[(3)] 
If $\alpha_1,\alpha_2,\alpha_3\in k^{\times}$, 
then the map $d(\alpha_1,\alpha_2,\alpha_3)$
belong to $G_{1\, k}=GE_{6\, k}$. Moreover,	   
$c(d(\alpha_1,\alpha_2,\alpha_3))=(\alpha_1\alpha_2\alpha_3)^2$
\item[(4)] 
The maps $n_{ij}(u)$ ($i\not=j$), 
belong to $H_{1\, k}=E_{6\, k}$. 
\end{itemize}
\end{lem}
\begin{proof}
The proof is easy for 
$d(\alpha_1,\alpha_2,\alpha_3)$. 
We only consider $d_1(\alpha)$ and 
$n_{21}(u)$ since other cases are similar. 

Let $\alpha\in \Oct^{\times}$. 
It is known that for
any $a, x, y \in \Oct$, 
\begin{equation}
\label{oct-formula}
a(xy)a = (ax)(ya), \quad a(xa) = (ax)a
\end{equation}
(see \cite[p.9, Proposition 1.4.1]{Springer}, 
\cite[p.10,  Lemma 1.4.2]{Springer}). 

We first consider $d_1(\alpha)$ ($\alpha\in\Oct^{\times}$). 
Let $X$ be as in (\ref{eq:X-defn}). 
By (\ref{eq:trace-relation}) and (\ref{oct-formula}), 
\begin{align*}
\tr\left((\alpha x_1\alpha)(\oline{\alpha} x_2)(x_3\oline{\alpha})\right) 
&=\tr\left((\alpha x_1\alpha)(\oline{\alpha} (x_2x_3)\oline{\alpha})\right) 
= \tr\left(\left((\alpha 
x_1\alpha)\oline{\alpha}\right)\left((x_2x_3)\oline{\alpha}\right)\right)\\
&= \tr\left(\left(((\alpha 
x_1)\alpha)\oline{\alpha}\right)\left((x_2x_3)\oline{\alpha}\right)\right)
= \|\alpha\|\tr\left(\left(\alpha 
x_1\right)\left((x_2x_3)\oline{\alpha}\right)\right) \\
& = \|\alpha\|\tr\left(\oline{\alpha}\left(\alpha 
x_1\right)(x_2x_3)\right)
= \|\alpha\|^2\tr\left( 
x_1(x_2x_3)\right).
\end{align*}
So, 
\begin{align*}
\det(d_1(\alpha)(X)) 
& = \|\alpha\|s_1s_2s_3 + \frac 1{\|\alpha\|}
\tr\left((\alpha x_1\alpha)(\oline{\alpha} x_2)(x_3\oline{\alpha})\right) \\
& \quad - \frac 1{\|\alpha\|} s_1 \|\alpha x_1\alpha\| - s_2\|\oline{\alpha} x_2\|
-s_3 \|x_3\oline{\alpha}\| \\
& = \|\alpha\|s_1s_2s_3 + \|\alpha\|
\tr(x_1x_2x_3) \\
& \quad- \|\alpha\| s_1 \|x_1\| -  \|\alpha\| s_2\|x_2\|
- \|\alpha\| s_3 \|x_3\| \\
& = \|\alpha\|\det(X). 
\end{align*}
Therefore, $d_1(\alpha)$ belongs to $GE_6$.
Moreover, $d_1(\alpha)$ belongs to $E_6$ if $\|\alpha\|=1$.

Next, we consider  $n_{21}(u)$ ($u\in\Oct$). 
Let $X$ be as in (\ref{eq:X-defn}). 
By taking the product of the first two matrices first, 
\begin{align*}
n_{21}(u)(X) 
& = 
\begin{pmatrix}
s_1 & x_3 & \overline {x_2} \\
\overline {x_3}+s_1 u & s_2 + ux_3 & x_1+u\overline{x_2} \\
x_2 & \overline{x_1} & s_3   
\end{pmatrix}
\begin{pmatrix}
1 & \overline u & \\
 & 1 & \\
&& 1  
\end{pmatrix} \\
& = \begin{pmatrix}
s_1 & x_3 + s_1\overline u & \overline{x_2} \\
\overline{x_3}+s_1 u & s_2 + \text{tr}(ux_3)+s_1\|u\| 
& x_1+u\overline{x_2} \\
x_2 & \overline{x_1}+x_2\overline u & s_3   
\end{pmatrix}.
\end{align*}
One can also verify that 
changing the order of multiplication 
does not change the result.  

By definition, 
\begin{align*}
\det(n_{21}(u)(X)) 
& = s_1(s_2 + \text{tr}(ux_3)+s_1\|u\|)s_3
+ \text{tr}(((x_1+u\overline{x_2})x_2)(x_3+s_1\overline u)) \\
& \quad - s_1\|x_1+u\overline{x_2}\|
- (s_2 + \text{tr}(ux_3)+s_1\|u\|)\|x_2\|
- s_3\|x_3+s_1\overline u\| \\
& = \det(X) +s_1s_3\text{tr}(ux_3)+s_1^2s_3\|u\| \\
& \quad + s_1\text{tr}((x_1x_2)\overline u) 
+ \text{tr}(((u\overline{x_2})x_2)x_3) 
+ s_1\text{tr}(((u\overline{x_2})x_2)\overline u)) \\
& \quad  - s_1\|u\| \|x_2\|
- s_1\text{tr}(x_1(x_2\overline u)) 
- (\text{tr}(ux_3)+s_1\|u\|)\|x_2\| \\
& \quad -s_1^2s_3\|u\| 
- s_1s_3\text{tr}(x_3u). 
\end{align*}
It is known that 
for any $x, y\in \Oct$,
\begin{equation*}
(xy)\overline y 
= x\|y\|, \; 
\overline y(yx) = x\|y\|. 
\end{equation*}
(see \cite[p.8, Lemma 1.3.3]{Springer}).
So, by this and the relations in (\ref{eq:trace-relation}), 
\begin{align*}
& \text{tr}((x_1x_2)\overline u) 
= \text{tr}(x_1(x_2\overline u)), \\ 
& \text{tr}(((u\overline{x_2})x_2)x_3) 
= \|x_2\|\text{tr}(ux_3)
= \|x_2\|\text{tr}(x_3u), \\
& \text{tr}(((u\overline x_2)x_2)\overline u))
= \|x_2\|\text{tr}(u\overline u)
= 2\|x_2\| \|u\|. 
\end{align*}
So all the terms except for $\det(X)$ cancel out and we obtain
\begin{math}
\det(n_{21}(u)(X)) = \det (X).  
\end{math}
\end{proof}

Note that for $X\in \J$ and 
$\alpha\in\Oct^{\times}$, 
the determinant of 
\begin{equation*}
\left(
{\rm{diag}}(\alpha,1,1)
X \right)
{\rm{diag}}(\overline \alpha,1,1)
\end{equation*}
may not be $\|\alpha\|\det(X)$. 

Since $\GL(3)$ is generated by elements of
the forms $d(\alpha_1,\alpha_2,\alpha_3)$ 
($\alpha_1,\alpha_2,\alpha_3\in k^{\times}$),
$n_{ij}(u)$ ($i\not=j$, $u\in k$), the 
following corollary follows. 
\begin{cor}
\label{cor:enchou}
For any $g_1 \in \GL(3)$ (resp. $g\in \cG$),
$\rho(g_1) \in G_1$ (resp. $\rho(g)\in G$). 
\end{cor}
Thus, if $x_1, x_2 \in \cV_k^{\rm ss}$ belong to the
same orbit in $\cG \backslash \cV_k^{\rm ss}$,
then they belong to the same 
orbit in $G_k \backslash \Vss_k$.

%
 
In the following, if $x=(x_1,x_2)\in V_k$, 
then we assume that  
\begin{equation}
\label{eq:xi-form}
x_i = \hji = 
\begin{pmatrix}
s_{i,1} & x_{i,3} & \overline {x_{i,2}} \\
\overline {x_{i,3}} & s_{i,2} & x_{i,1} \\
x_{i,2} & \overline {x_{i,1}} & s_{i,3}
\end{pmatrix} 
\end{equation}
where $s_{i,j}\in k,x_{i,j}\in \mathbb O$  
for $i=1,2,j=1,2,3$.

\begin{lem}
\label{lem:detzerohenkei}  
Suppose that $x=(x_1,x_2)\in V^{\rm{ss}}_k$
and that $\det(x_1)=0$. Then one can choose 
$g\in G_{k}$ so that after replacing 
$x$ by $gx$, $x$ is in the form:
\begin{equation}
\label{eq:x-form-lemma}
\left(
\begin{pmatrix}
0 & 1 & 0 \\
1 & 0 & 0 \\
0 & 0 & 0  
\end{pmatrix},
\begin{pmatrix}
-1 & 0 & 0 \\
0 & s_{2,2} & 0 \\
0 & 0 & 1 
\end{pmatrix}
\right). 
\end{equation}
\end{lem}
\begin{proof}
We first show that there exists $g\in G_1$ such that 
the $(1,1)$-entry of $gx_1$ is non-zero.

Let $W'$ be the subgroup of $\GL(3)_k$ 
consisting of permutation matrices. 
Then $W'$ is isomorphic to  $\mathfrak S_3$. 
Let $W=\rho_1(W')\subset G_1$. 
We put 
\begin{equation}
\label{eq:nu-prime-defn}
\nu' = 
\begin{pmatrix}
0 & 0 & 1 \\
0 & 1 & 0 \\
1 & 0 & 0  
\end{pmatrix},\;
\nu=\rho_1(\nu').  
\end{equation}
Then $\nu'\in W',\nu\in W$.

If there exists $1\leq j\leq 3$ such that  
$s_{1,j}\not=0$, then by applying an element of $W$, 
we may assume that $s_{1,1}\not=0$.
Suppose that all diagonal entries of $x_1$ are $0$. 
Since $x\in V_k^{\rm{ss}}$, $x_1\not=0$. 
So there exists $j$ such that $x_{1,j}\not=0$. 
By applying an element of $W$, we may assume that 
$x_{1,3}\not=0$. If $u\in\Oct$, then the $(1,1)$-entry of 
$n_{12}(u)(x_1)$ is $\text{tr}(u\oline{x_{1,3}})$. Since 
$x_{1,3}\not=0$, there exists $u\in \mathbb O$
such that $\text{tr}(u\oline{x_{1,3}})\not=0$. 
Therefore, we may assume that $s_{1,1}\not=0$.

By applying elements of the forms
$n_{21}(u),n_{31}(u)$ ($u\in\Oct$), 
we may assume that $x_{1,2}=x_{1,3}=0$. 
If $s_{1,2}=s_{1,3}=x_{1,1}=0$, then 
the weights of all non-zero coordinates 
of $x$ with respect to the one parameter subgroup
(which will be abbreviated as 1PS from now on)
\begin{equation*}
\lambda:\text{GL}(1)\ni \al\mapsto 
\left(
d(\alpha^2,\alpha^{-1},\alpha^{-1}), 
{\rm{diag}}(\al^{-3},\al^3)
\right)\in G_k
\end{equation*}
is positive (this $\GL(1)$ 
is $\GL(1)$ over $k$). 
This contradicts to 
the assumption that $x\in V^{\rm{ss}}_k$. 
This is well-known
(see \cite[p.49, Theorem 2.1]{mufoki})
but can be seen in this situation as follows also. 
If $\Delta(x)$ is the relative invariant polynomial
of degree $12$ as in Introduction, 
since $\text{Im}(\lambda)\subset E_6\times \text{SL}(2)$, 
$\Delta(\lambda(\al)x)=\Delta(x)$. However, 
since the weights of all non-zero coordinates 
with respect to $\lambda$ are positive, 
there exists a polynomial $P(x,\alpha)$ 
such that $\Delta(x)= \Delta(\lambda(\alpha)x)=\alpha P(x,\alpha)$. 
Since $\Delta(x)$ is a non-zero polynomial which does not
depend on $\alpha$, this is a contradiction.

So $(s_{1,2},s_{1,3},x_{1,1})\not=(0,0,0)$. 
By a similar argument as above, 
there exists $g\in G_{1\,k}$ such that 
$gx_1$ is in the form 
\begin{math}
\text{diag}(s_{1,1},s_{1,2},s_{1,3})  
\end{math} 
where $s_{1,1},s_{1,2}\not=0$. 
Since $\det(x_1)=0$, $s_{1,3}=0$.

If $s_{2,3}=0$ then 
the weights of all non-zero coordinates with respect to 
the 1PS 
\begin{equation*}
\lambda:\text{GL}(1)\ni \alpha\mapsto 
\left(
d(\alpha^2,\alpha^2,\alpha^{-4}), 
{\rm{diag}}(\al^{-3},\al^3)
\right)
\end{equation*}
is positive, which is a contradiction.   
Therefore, $s_{2,3}\not=0$.

Then by applying elements of the forms 
$n_{31}(u),n_{32}(u)$, 
we may assume that $x_{2,1}=x_{2,2}=0$. 
Note that $n_{31}(u),n_{32}(u)$ 
do not change $x_1$. 
Therefore, $x$ is in the form:
\begin{equation}
\label{eq:x-form2}
\left(
\begin{pmatrix}
s_{1,1} && \\
& s_{1,2} & \\
&& 0  
\end{pmatrix},
\begin{pmatrix}
s_{2,1} & x_{2,3} & 0 \\ 
\overline {x_{2,3}} & s_{2,2} & 0 \\
0 & 0 & s_{2,3}
\end{pmatrix} 
\right).  
\end{equation}

Since the norm map $\Oct^{\times}\to k^{\times}$ is 
surjective, by applying elements of the forms
$d_1(\alpha_1),d_2(\alpha_2),d_3(\alpha_3)$ 
($\alpha_1,\alpha_2,\alpha_3\in\Oct^{\times}$), 
we may assume that $s_{1,1}=1/2,s_{1,2}=-1/2,s_{2,3}=1$. 
Since
\begin{equation}
\label{eq:diagonal2antidiagonal}
\begin{pmatrix}
1 & -1 & 0 \\
1 & 1 & 0 \\
0 & 0 & 1  
\end{pmatrix}
\begin{pmatrix}
\frac 12  & 0 & 0 \\
0 & -\frac 12  & 0 \\
0 & 0 & 0 
\end{pmatrix}
\begin{pmatrix}
1 & 1 & 0 \\
-1 & 1 & 0 \\
0 & 0 & 1  
\end{pmatrix}
= 
\begin{pmatrix}
0 & 1 & 0 \\
1 & 0 & 0 \\
0 & 0 & 0 
\end{pmatrix},  
\end{equation}
we may assume that 
$x$ is in the form: 
\begin{equation}
\label{eq:detx0form}
x = (x_1,x_2) = \left(
\begin{pmatrix}
0 & 1 & 0 \\
1 & 0 & 0 \\
0 & 0 & 0 
\end{pmatrix},
\begin{pmatrix}
s_{2,1} & x_{2,3} & 0 \\ 
\overline {x_{2,3}} & s_{2,2} & 0 \\
0 & 0 & 1
\end{pmatrix} 
\right).  
\end{equation}

We first show that there exists $g\in G_{k}$
such that $gx$ is still in the form 
(\ref{eq:detx0form}) and that 
$s_{2,1}\not=0$.

Suppose that $s_{2,1}=0$. If 
$s_{2,2}\not=0$, then we only have to apply an 
element of $W$. 
So we assume  that $s_{2,1}=s_{2,2}=0$. 
In this situation, if 
$u\in\Oct$ and ${\rm{tr}}(u)=0$, 
then $n_{12}(u)x_1=x_1$ and 
the $(1,1)$-entry of $n_{12}(u) x_2$ 
is ${\rm{tr}}(u\overline {x_{2,3}})={\rm{tr}}(\overline {x_{2,3}}u)$. 
Let $\overline {x_{2,3}}=a+b$ where 
$a\in k$ and ${\rm{tr}}(b)=0$. 
Then ${\rm{tr}}(\overline {x_{2,3}}u)={\rm{tr}}(bu)$. 
The assumption implies that $b\not=0$. 
Suppose that ${\rm{tr}}(bu)=0$ for all $u\in\Oct$ such that
${\rm{tr}}(u)=0$. Since ${\rm{tr}}(bc)=0$ if $c\in k$, 
${\rm{tr}}(bu)=0$ for all $u\in\Oct$. 
Since $\Oct^2\ni (u_1,u_2)\mapsto {\rm{tr}}(u_1u_2)$
is a non-degenerate bilinear form, 
$b=0$, which is a contradiction.  Therefore, 
there exists $u\in \Oct$ such that 
${\rm{tr}}(u)=0$ and that ${\rm{tr}}(\overline{x_{23}}u)\not=0$. 
So we may assume that $s_{2,1}\not=0$. 

We next show that 
there exists $g\in G_{k}$
such that $gx$ is still in the form 
(\ref{eq:detx0form}) and that 
$x_{2,3}\in k$. 

\begin{lem}
\label{lem:belongtok}
Suppose that $s\in k^{\times}$ and $x\in \Oct$. 
Then there exists $u\in\Oct$ such that  
$\mathrm{tr}(u)=0$ and $x+su\in k$. 
\end{lem} 
\begin{proof}
It is enough to take $u=(1/2s)(\overline x-x)$. 
\end{proof}

If ${\rm{tr}}(u)=0$, then $n_{12}(u)x_1=x_1$ and 
the $(1,2)$-entry of $n_{12}(u)x_2$ is 
$x_{2,3}+s_{2,1}u$. Since $s_{2,1}\not=0$, 
there exists $u$ such that ${\mathrm{tr}}(u)=0$ and 
$x_{2,3}+s_{2,1}u\in k$ by Lemma \ref{lem:belongtok}. 
So we may assume further that $x_{2,3}\in k$. 
Then if we apply the element 
\begin{math}
\left(\begin{smallmatrix}
1 & 0  \\ -x_{2,3} & 1       
\end{smallmatrix}\right) \in \GL(2)_k, 
\end{math}
then $x_{2,3}$ changes to $0$.   

By applying $d_2(\alpha_2)d_1(\alpha_1)$, 
$x$ changes to 
\begin{equation*}
\left(
\begin{pmatrix}
0 & \overline {\alpha_2} \, \overline {\alpha_1} & 0 \\
\alpha_1 \alpha_2 & 0 & 0 \\
0 & 0 & 0 
\end{pmatrix},
\begin{pmatrix}
s_{2,1}\|\alpha_1\| & 0 & 0 \\ 
0 & s_{2,2}\|\alpha_2\| & 0 \\
0 & 0 & 1
\end{pmatrix} 
\right).  
\end{equation*}
We can choose $\alpha_1,\alpha_2$ so that 
$s_1\|\alpha_1\|=-1$ and $\alpha_1\alpha_2=1$. 
Then $x$ changes to 
\begin{equation*}
\left(
\begin{pmatrix}
0 & 1 & 0 \\
1 & 0 & 0 \\
0 & 0 & 0 
\end{pmatrix},
\begin{pmatrix}
-1 & 0 & 0 \\ 
0 & s_{2,2} & 0 \\
0 & 0 & 1
\end{pmatrix} 
\right).  
\end{equation*}
This completes the proof of the lemma. 
\end{proof}

Suppose that $F_x(v)$ has a rational factor. 
By applying an element of $\GL(2)_k$, 
we may assume that $F_x(1,0)=0$.  
Then $\det(x_1)=0$. 
By Lemma \ref{lem:detzerohenkei}, we may assume that
$x\in \cV^{\rm{ss}}_k$ and is in the form 
(\ref{eq:x-form-lemma}). 
Then $(1,0,1),(-1,0,1)\in \Zero_{\Proj^2}(x)$ are
rational points. 
Therefore, the element of ${\rm{H}}^1(k,\mathfrak S_4)$
which corresponds to $x$ with respect to the action of
$\cG_k$ comes from ${\rm{H}}^1(k,\mathfrak S_2)$
and so comes from ${\rm{H}}^1(k,\cS)$.

Finally we consider the case where 
$F_x$ is irreducible. 
We assume that $x_1,x_2$ are in the form 
(\ref{eq:xi-form}).  By assumption, $\det(x_1)\not=0$. 
By the same argument as in the proof of 
Lemma \ref{lem:detzerohenkei},  
we may assume that $x_1$ is a diagonal matrix. 
Since $\det(x_1)\not=0$, all diagonal entries are non-zero.
Also since the norm map $\Oct^{\times}\to k^{\times}$
is surjective, we may assume that 
\begin{math}
x_1 = {\rm{diag}}
(\tfrac 12,-2,-\tfrac 12).  
\end{math}
Let
\begin{equation}
\label{eq:Lambda-defn}
g = \begin{pmatrix}
1 & & 1 \\
& 1 & \\
1 & & -1  
\end{pmatrix}, \;
\Lambda = \begin{pmatrix}
0 & 0 & 1 \\
0 & -2 & 0 \\
1 & 0 & 0  
\end{pmatrix}.  
\end{equation}
Then 
\begin{math}
\rho_1(g) x_1
= \Lambda. 
\end{math}
So we may assume that $x_1=\Lambda$.

For $u \in \Oct$, let 
\begin{equation*}
A_1(u) =  \begin{pmatrix}
1 & 0 & 0 \\
u & 1 & 0 \\
\frac 14\|u\| & \frac 12 \oline u & 1
\end{pmatrix}, \quad
A_2(u) =  \begin{pmatrix}
1 & \frac 12 \overline u & \frac 14\|u\| \\
0 & 1 & u \\
0 & 0 & 1
\end{pmatrix}.
\end{equation*} 
Then one can easily verify that 
$(A_i(u) \Lambda) {}^t\overline {A_i(u)}
=A_i(u) (\Lambda {}^t\overline {A_i(u)})=\Lambda$ 
($i=1,2$). 

\begin{lem}
\label{lem:Au-defn}
\begin{itemize}
\item[(1)] 
For $Y\in\J$, 
$(A_i(u) Y) {}^t\oline{A_i(u)}=A_i(u) (Y {}^t\oline{A_i(u)})\in\J$
for $i=1,2$. 
\item[(2)] 
If we put $B_i(u)(Y)=A_i(u) Y {}^t\oline{A_i(u)}$, 
then $B_i(u):\J\to\J$ is a $k$-linear map 
which belongs to $E_6$ for $i=1,2$. 
\end{itemize} 
\end{lem}
\begin{proof}
Since the argument is similar, 
we only consider $A_1(u)$.

(1) \, Suppose that $Y=h(t_1,t_2,t_3,y_1,y_2,y_3)$ 
where $t_1,t_2,t_3\in k, y_1,y_2,y_3\in\Oct$. 
Then 
\begin{align*}
& (A_1(u) Y) {}^t\oline{A_1(u)} \\
& = 
\begin{pmatrix}
t_1 & y_3 & \oline{y_2} \\
\oline{y_3}+t_1u & t_2 + uy_3 
& y_1 + u\oline{y_2} \\
y_2+\frac 12 \oline{u} \, \oline{y_3} +\frac 14 t_1\|u\| 
& \oline{y_1} + \frac 12t_2\oline{u}  +\frac 14 \|u\| y_3 
& t_3 + \frac 12\oline{u}y_1 + \frac 14 \|u\|\, \oline{y_2} 
\end{pmatrix} \\
& \quad \times 
\begin{pmatrix}
1 & \oline u & \frac 14\|u\| \\
0 & 1 & \frac 12 u \\
0 & 0 & 1
\end{pmatrix} \; 
\text{($\times$ is the matrix multiplication)}\\
& = (A_{1,ij})
\end{align*}
where
\begin{align*}
A_{1,11} & = t_1, \\
A_{1,12} & = y_3+ t_1\oline{u}, \\
A_{1,13} & = \overline {y_2} +\frac 12 y_3u + \frac 14 t_1\|u\|, \\
A_{1,21} & = \oline{y_3} +t_1u, \\
A_{1,22} & = t_2 + {\rm{tr}}(uy_3) + t_1 \|u\|, \\
A_{1,23} & = y_1+u\overline{y_2}+\frac 12 t_2u
+\frac 12(uy_3)u+\frac 14\oline{y_3}\|u\|+\frac 14t_1\|u\|u, \\
A_{1,31} & = y_2 + \frac 12  \oline{u}\, \oline{y_3} + \frac 14 t_1\|u\|, \\
A_{1,32} & = \oline{y_1} + y_2\oline{u} + \frac 12 t_2 \oline{u}
+\frac 12 (\oline{u}\,\oline{y_3})\oline{u}
+\frac 14 y_3\|u\|+\frac 14t_1\|u\|\oline{u}, \\
A_{1,33} & = t_3 + \frac 12{\rm{tr}}(\oline{y_1}u)
+\frac 14t_2\|u\|+\frac 14{\rm{tr}}(y_2)\|u\|
+\frac 18{\rm{tr}}(y_3u)\|u\|
+\frac 1{16}t_1\|u\|^2.
\end{align*}

The only places where the order of multiplication matters
are $\frac 12(uy_3)u$ in $A_{23}$ and 
$\frac 12(\oline{u} \,\oline{y_3})\oline{u}$ in $A_{32}$.
However, by (\ref{oct-formula}), 
these are equal to $\frac 12u(y_3u)$,  
$\frac 12\oline{u} (\oline{y_3}\, \oline{u})$ respectively. 
So $(A_1(u) Y) {}^t\oline{A_1(u)}=A_1(u) (Y {}^t\oline{A_1(u)})$.
Also, 
\begin{equation}
\label{eq:ux3u}
\oline{(uy_3)u}=\oline{u}(\oline{y_3}\,\oline{u})
=(\oline{u}\,\oline{y_3})\oline{u}. 
\end{equation}
Therefore, 
$\oline{A_{1,23}}=A_{1,32}$. The conditions
$\oline{A_{1,12}}=A_{1,21}$, 
$\oline{A_{1,13}}=A_{1,31}$ and $A_{1,ii}\in k$ ($i=1,2,3$) 
are obvious. 
Thus $B_1(u)(Y) \in \J$.

The consideration is similar for 
$B_2(u)(Y)$, but 
we list entries of 
$B_2(u)(Y)$ since they will be needed later. 
Let $B_2(u)(Y)=(A_{2,ij})$. 
Then 
\begin{align*}
A_{2,11} & = t_1 + \frac 12{\rm{tr}}(y_3 u)
+\frac 14t_2\|u\|+\frac 14{\rm{tr}}(y_2)\|u\|
+\frac 18{\rm{tr}}(\overline{y_1}u)\|u\|
+\frac 1{16}t_3\|u\|^2, \\
A_{2,12} & = y_3 + \overline{y_2} \,\oline{u} + \frac 12 t_2 \oline{u}
+\frac 12 (\oline{u}\,y_1)\oline{u}
+\frac 14 \overline{y_1}\|u\|+\frac 14t_3\|u\|\overline u, \\ 
A_{2,13} &  = \oline{y_2} + \frac 12 \oline u y_1 
+ \frac 14t_3\|u\|, \\
A_{2,22} & = t_2 + \tr(u\oline{y_1}) + t_3\|u\|,\\
A_{2,23} & = y_1 + t_3 u, \\
A_{2,33} & = t_3, \\
A_{2,21} &= \oline{A_{2,12}}, \, A_{2,31} =  \oline{A_{2,13}}, \,
A_{2,32} = \oline{A_{2,23}}.
\end{align*}

(2) By easy computations, 
\begin{equation*}
A_1(u) = 
\begin{pmatrix}
1 & 0 & 0 \\
0 & 1 & 0 \\
-\frac 14 \|u\| & 0 & 1
\end{pmatrix}
\begin{pmatrix}
1 & 0 & 0 \\
0 & 1 & 0 \\
0 & \frac 12\oline u & 1
\end{pmatrix}
\begin{pmatrix}
1 & 0 & 0 \\
u & 1 & 0 \\
0 & 0 & 1
\end{pmatrix}
\end{equation*}
where the order of multiplication does not matter. 
One can verify by tedious computations using (\ref{eq:ux3u}) that 
$B_1(u) = n_{31}(-\frac 14 \|u\|) n_{32}(\frac 12\oline u)\, n_{21}(u)$
(this is not trivial since $\Oct$ is not associative). 
Therefore, $B_1(u)$ belongs to $G_{1\,k}$.
\end{proof}

We successively apply elements of $G_k$ to 
modify $x$.  

We eventually want to make $s_{2,1}=0$ and
$\|x_{2,3}\|\not=0$. However, for that purpose, 
we first make $s_{2,1}$ non-zero to eliminate some 
entries of $x_2$. 
If $s_{2,1}=0$ and $s_{2,3}\not=0$, 
then we only have to apply $\nu$ (see (\ref{eq:nu-prime-defn})). 
Note that $\nu x_1=x_1$. 
Suppose that 
$s_{2,1}=s_{2,3}=0$. 
If further $x_{2,3}=0$, then 
by applying 
\begin{math}
\left(\begin{smallmatrix}
1 & 0 \\
\frac 12s_{2,2} & 1        
\end{smallmatrix}\right)\in\GL(2)_k,  
\end{math}
$s_{2,2}$ becomes $0$. 
Then $\det(x_2)=0$, which is a contradiction. 
Therefore, $x_{2,3}\not=0$. 
If we choose $u\in\Oct$ such that $\|u\|=0$, 
then by applying  $B_2(u)$ to $x_2$, 
$s_{2,1}$ (which is $0$ now) becomes
$\frac 12{\rm{tr}}(x_{2,3} u)$ and $x_1$ 
does not change.

We can choose a basis of $\Oct$ as a $k$-vector space
so that according to the 
coordinate $u=(u_1,\ldots,u_8)$
of $u$ with respect to this basis, 
$\|u\|=u_1u_2+u_3u_4+u_5u_6+u_7u_8$ 
($\Oct$ is split). 
Since $\frac 12{\rm{tr}}(x_{2,3} u)$ 
is a non-zero $k$-linear function, 
it is in the form 
$c_1u_1+\cdots + c_8u_8$ ($c_1,\ldots,c_8\in k$). 
We may assume that $c_1\not=0$. 
We choose $u_2=1$. Then 
$\|u\|=0$ if and only if 
$u_1=-(u_3u_4+u_5u_6+u_7u_8)$ 
and so 
\begin{equation*}
\frac 12{\rm{tr}}(x_{2,3} u)
= -c_1(u_3u_4+u_5u_6+u_7u_8)+c_2+c_3u_3+\cdots+c_8u_8. 
\end{equation*}
This is a non-zero polynomial. Since $k$ 
is an infinite field, there exists 
$u\in\Oct$ such that 
$\|u\|=0$ and that 
$\frac 12{\rm{tr}}(x_{2,3} u)\not=0$. 
Therefore, we may assume that 
$s_{2,1}\not=0$ for the moment
(we later change $s_{2,1}$ to $0$ after 
eliminating some entries of $x_2$).

By applying an element of the form 
$B_1(u)$ ($u\in\Oct$), $x_1$ does not change and 
$x_{2,3}$ changes to 
$x_{2,3}+s_{2,1}\overline u$. Since $s_{2,1}\not=0$, 
there exists $u$ such that $x_{2,3}+s_{2,1}\overline u =0$.  
Therefore, we may assume that 
$x_{2,3}=0$ for the moment (we later change $x_{2,3}$ 
so that $\|x_{2,3}\|\not=0$).

By applying the element 
\begin{math}
\left(\begin{smallmatrix}
1 & 0 \\
\frac 12 s_{2,2} & 1        
\end{smallmatrix}\right) 
\in \GL(2)_k,
\end{math}
we may assume that $s_{2,2}=0$. 
So $x$ is in the form: 
\begin{equation*}
\left(
\begin{pmatrix}
 && 1 \\
& -2 & \\
1 &&   
\end{pmatrix},
\begin{pmatrix}
s_{2,1} & 0 & \overline x_{2,2} \\
0 & 0 & x_{2,1} \\
x_{2,2} & \overline{x_{2,1}} & s_{2,3}
\end{pmatrix}
\right).  
\end{equation*}
Then $\det(x_2)=s_{2,1}\|x_{2,1}\|\not=0$. 
Therefore, $x_{2,1}\in \Oct^{\times}$. 


We want to make $x_{2,1} \in k$ while preserving 
$x_1 = \Lambda$ (see (\ref{eq:Lambda-defn})).
So  we prove the following lemma.
\begin{lem}
\label{change-oct}
Let $y = h(t_1,t_2,t_3,y_1,y_2,y_3) \in \J$ 
where $t_i \in k, y_i \in \Oct$ 
$(i= 1,2,3)$. 
If $\|y_1\| \neq 0$, then 
there exists $g \in G_{1\Lambda\,k}$ ($ G_{1\Lambda}$ 
is the stabilizer of $\Lambda$ in $G_1$)
such that
the $(2,3)$ entry 
of $gy$ is $1$ and $g$ preserves subspaces $kE_i$, $\Oct_i$
(see (\ref{elements-defn}) and (\ref{eq:oct-i-defn}))
of $\J$ for $i = 1,2,3$.
\end{lem}
\begin{proof}
We construct such an element $g$ by combining the maps $d_i(\be)$
where $\be \in \Oct^{\times}, i = 1,2$.
For $\be\in \Oct^{\times}$, 
\begin{equation*}
d_1(-\|\be\|) d_2(\be)(\Lambda) = 
\begin{pmatrix}
 && -\oline{\be}^2 \\
& -2\|\be\| & \\
-\be^2 &&   
\end{pmatrix}.
\end{equation*}
If $\tr(\be) = 0$, 
then since $\be = -\oline \be$, $\be^2 = \oline{\be}^2 = - \|\be\|$
and so
$d_1(-\|\be\|) d_2(\be)(\Lambda)$ $= \|\be\|\Lambda$.
Hence, the element
\begin{math}
\|\be\|^{-1}I_{27}\,
d_1(-\|\be\|) d_2(\be) \in G_{1\,k}  
\end{math}
fixes $\Lambda$
for any $\be \in {\bbW} \cap \Oct^{\times}$ 
where ${\bbW}$ is the the orthogonal complement of $k\cdot1$ in $\Oct$.
For $\be \in \Oct^{\times}$, we denote the element
\begin{math}
\|\be\|^{-1}I_{27}\,
d_1(-\|\be\|) d_2(\be) \in G_{1\,k}  
\end{math}
by $\mathcal{D}(\be)$.
By the above argument, if $\tr(\be) = 0 \,( \be \in \bbW)$, 
then $\mathcal{D}(\be) \in G_{1\Lambda\, k}$.
Since $d_j(\al)$ where $\al \in \Oct^{\times}, j = 1,2$ 
preserves subspaces $kE_i$, $\Oct_i$
of $\J$ for $i = 1,2,3$,
$\mathcal{D}(\be)$ also preserves these subspaces.
By applying  $\mathcal{D}\left(\,\oline{\be}\,\right)$ where 
$\be \in \bbW \cap \Oct^{\times}$ to $y$, 
$y_1$ changes to $\|\be\|^{-1} y_1\be$ 
$(d_1(-\|\oline{\be}\|)$
does not change $y_1$)
and so $\tr(y_1)$ changes to $\|\be\|^{-1}\tr( y_1\be)$. 

If 
there exists $\be \in \bbW\cap \Oct^{\times}$
such that $\tr(\be y_1) = 0$,
then 
the $(2,3)$ entry of $\mathcal{D}\left(\|\be\|^{-1}\oline{y_1\be}\right)
\left(\mathcal{D}\left(\oline{\be}\right)(y)\right)$
is 
\begin{equation*}
\|(\|\be\|^{-1}(\,\oline{y_1\be}\,))\|^{-1}
(\|\be\|^{-1}(y_1\be))(\|\be\|^{-1}(\,\oline{y_1\be}\,))
= 1.
\end{equation*}
Note that by assumption, $\left\|\oline{y_1\be}\right\| \neq 0$.
Since $\mathcal{D}\left(\|\be\|^{-1}\oline{y_1\be}\right)
\mathcal{D}\left(\,\oline{\be}\,\right)$ belongs to $G_{1\Lambda\,k}$ and
 preserves subspaces $kE_i$, $\Oct_i$
of $\J$ where $i = 1,2,3$, this element is what we wanted.
So it is enough to show that there exists 
$\be \in \bbW\cap \Oct^{\times}$
such that $\tr(\be y_1) = 0$.

 For $a\in \Oct$, we denote the linear map 
$\Oct \ni x \to \tr(ax) \in k$ by $\tr_a$
and denote the restriction of $\tr_a$ 
on $\bbW$ by $\tr_a|_{\bbW}$.
Then $\mathrm{dim}(\Ker(\tr_{y_1}|_{\bbW})) \ge 6$.  
Since $\Oct$ is split, the dimension of the maximal totally isotropic subspace, 
which is called the Witt index, 
of $\Oct$ is $4$ ($= \frac 12$ dim $\Oct$).
Thus, there exists $\be \in \Ker(\tr_{y_1}|_{\bbW})$ 
such that $\|\be\| \neq 0$.
This completes the proof of the lemma.
\end{proof}
%

%

By the above lemma,
we may assume that $x$ is in the form: 
\begin{equation*}
\left(
\begin{pmatrix}
 && 1 \\
& -2 & \\
1 &&   
\end{pmatrix},
\begin{pmatrix}
s_{2,1} & 0 & \overline{x_{2,2}} \\
0& 0 & 1 \\
x_{2,2} & 1 & s_{2,3}
\end{pmatrix}
\right)
\end{equation*}
where $s_{2,1} \neq 0$.

If $u\in\Oct$ and ${\mathrm{tr}}(u)=0$, 
then $n_{31}(u)$ does not change $x_1$ and 
$x_{2,2}$ changes to $x_{2,2}+s_{2,1}u$. 
By Lemma \ref{lem:belongtok}, 
there exists $u$ such that 
${\mathrm{tr}}(u)=0$ and $x_{2,2}+s_{2,1}u\in k$. 
Therefore, we may assume that $x_{2,2}\in k$. 
By applying the element 
\begin{math}
\left(
\begin{smallmatrix}
1 & 0 \\
-x_{2,2}s_{2,1}^{-1} & s_{2,1}^{-1}
\end{smallmatrix}
\right), 
\end{math}
we may assume that
$x$ is in the form: 
\begin{equation*}
\left(
\begin{pmatrix}
 && 1 \\
& -2 & \\
1 &&   
\end{pmatrix},
\begin{pmatrix}
1 & 0 & 0 \\
0& s_{2,2} & c \\
0 & c & s_{2,3}
\end{pmatrix}
\right)
\end{equation*}
where $c\in k^{\times}$.

By applying $B_1(2)$ and then 
\begin{math}
\left(
\begin{smallmatrix}
1 & 0 \\
-1 & 1
\end{smallmatrix}
\right),  
\end{math}
$x_1$ does not change and 
$x_2$ becomes
%
\begin{equation*}
x_2' = \begin{pmatrix}
1 & 2 & 0 \\
2 & 6 + s_{2,2} & 2 + c + s_{2,2} \\
0 & 2 + c + s_{2,2} & 1+2c + s_{2,2} + s_{2,3}
\end{pmatrix}.
\end{equation*}
We denote the $(2,2), (2,3), (3,3)$ entries of $x_2'$
by $s'_{2,2}, c', s'_{2,3}$ respectively.

We choose an element  $u\in \Oct$ such that
$\tr(u) = -1, \|u\| = 0$.
For example, we can choose $u =  
\left(\frac 12 \left(
\begin{smallmatrix}
-1&1\\
1&-1
\end{smallmatrix}
\right),0\right) \in \mathrm{M}(2)\oplus \mathrm{M}(2)
= \Oct.
$
Then $u^2 = -u$.
By applying  $B_2(u)$, $x_1$ does not change and
$x_2'$ becomes
\begin{align*}
x_2'' = B_2(u)(x_2') 
& = 
\begin{pmatrix}
1 + \tr(u)&2 + \frac 12 s'_{2,2}\oline u + \frac 12 c' \oline u^2
&\frac 12 c' \oline u 
\\
2 + \frac 12 s'_{2,2}u + \frac 12 c' u^2
& s'_{2,2} + c'\tr(u) & c' + s'_{2,3}u 
\\
\frac 12 c' u & c' + s'_{2,3}\oline u & s'_{2,3}
\end{pmatrix}\\
& = 
\begin{pmatrix}
0 & 2 + \frac 12 (s'_{2,2}- c') \oline u
&\frac 12 c' \oline u 
\\
2 + \frac 12 (s'_{2,2}- c') u
&s'_{2,2} - c' & c'  + s'_{2,3}u 
\\
\frac 12 c' u &c' + s'_{2,3}\oline u & s'_{2,3}
\end{pmatrix}.
\end{align*}
We denote the $(1,2)$-entry of $x_2''$ by 
$x_{2,3}''$.  Then 
\begin{align*}
\|x_{2,3}''\| 
& = \left\|2 + \frac 12 (s'_{2,2}- c') \oline u \right\| 
= 4 + (s'_{2,2}- c') \tr(u) 
= 4 - s'_{2,2} + c'.
\end{align*}
Since $s'_{2,2} = 6 + s_{2,2}, c' = 2 + c + s_{2,2}$,
we have $\|x_{2,3}\| = c\not=0$.
Therefore, we may assume that $x$ is in the form:
\begin{equation*}
x = (x_1,x_2) = \left(
\begin{pmatrix}
0 & 0 & 1 \\
0 & -2 & 0 \\
1 & 0 & 0 
\end{pmatrix},
\begin{pmatrix}
0 & x_{2,3} & \overline {x_{2,2}} \\
\overline {x_{2,3}} & s_{2,2} & x_{2,1} \\
x_{2,2} & \overline {x_{2,1}} & s_{2,3}
\end{pmatrix}
\right) 
\end{equation*}
where $x_{2,3}\in\Oct^{\times}$.

We apply Lemma \ref{change-oct} to $\nu x_2$.
Note that $\nu\Lambda = \Lambda$.
By Lemma \ref{change-oct}, 
there exists $g\in G_{1\Lambda\, k}$
such that the $(2,3)$ and $(3,3)$ entries of $g(\nu x_2)$ 
are $1$ and $0$ respectively.
So, by applying  $\nu$ to
$g(\nu x_2)$, we may assume that $x$ is in the
form:
\begin{equation*}
\left(
\begin{pmatrix}
 && 1 \\
& -2 & \\
1 &&   
\end{pmatrix},
\begin{pmatrix}
0 & 1 & \overline{x_{2,2}} \\
1 & s_{2,2} & x_{2,1} \\
x_{2,2} & \overline{x_{2,1}} & s_{2,3}
\end{pmatrix}
\right).  
\end{equation*}

By applying  $B_1(-2\overline{x_{2,2}})$, 
we may assume that $x_{2,2}=0$. 
By applying an element of the form
$n_{31}(u)$ where $\text{tr}(u)=0$, 
Lemma \ref{lem:belongtok} implies that  
we may assume that $x_{2,1}\in k$. 
Therefore, $x$ is in the form: 
\begin{equation*}
x = \left(
\begin{pmatrix}
 && 1 \\
& -2 & \\
1 &&   
\end{pmatrix},
\begin{pmatrix}
0 & 1 & 0 \\
1 & s_{2,2} & x_{2,1} \\
0 & x_{2,1} & s_{2,3}
\end{pmatrix}
\right)
\end{equation*}
where $s_{2,2},x_{2,1},s_{2,3}\in k$. 
Then $x\in\cV^{\mathrm{ss}}_k$ and 
$\Zero_{\Proj^2}(x)$ has a rational 
point $(1,0,0)$. Therefore, 
the element of ${\rm{H}}^1(k,\mathfrak S_4)$
which corresponds to $x$ now comes from 
an element of ${\rm{H}}^1(k,\cS)$. 
This completes the proof of Proposition \ref{prop:crucial-step}.
\end{proof}

Since Proposition \ref{prop:crucial-step} is proved,  
the proof of Theorem \ref{thm:main-split-intro} is finished now. 
\end{proof}

%

\bibliographystyle{plain}
\bibliography{ref-rep7}

\end{document}